\documentclass[a4paper,11pt,dvips]{article}
\usepackage[latin1]{inputenc}
\usepackage[OT1,OT2,T1]{fontenc}
\usepackage{a4wide}
\usepackage{amsmath,amssymb}
\usepackage{mathrsfs}
\usepackage[ps,all]{xy}
\usepackage[vcentermath]{youngtab}
\newcommand\Rep{\mathord{\mathrm{Rep}}}
\newcommand\Sym{\mathord{\mathrm{Sym}}}
\newcommand\QSym{\mathord{\mathrm{QSym}}}
\newcommand\FQSym{\mathord{\mathbf{FQSym}}}

\newcommand\NTQSym[1]{\mathord{\mathit{QSym}^{(#1)}}}
\newcommand\T{\mathord{\mathsf{T}}}
\newcommand\TS{\mathord{\mathsf{TS}}}
\newcommand\NT{\mathord{\mathsf{NT}}}
\newcommand\MR{\mathord{\mathsf{MR}}}
\newcommand\ch{\mathop{\mathrm{ch}}\nolimits}
\newcommand\Irr{\mathord{\mathrm{Irr}}}
\newcommand\Hom{\mathop{\mathrm{Hom}}\nolimits}

\newcommand\id{\mathord{\mathrm{id}}}
\newcommand\sha{\mathbin{\text{\fontencoding{OT2}\fontfamily{cmss}%
  \selectfont sh}}}
\newcommand\Ind{\mathop{\mathrm{Ind}}\nolimits}
\newcommand\Res{\mathop{\mathrm{Res}}\nolimits}
\newcommand\rad{\mathop{\mathrm{rad}}\nolimits}
\newcommand\std{\mathop{\mathrm{std}}\nolimits}
\newcommand\ev{\mathord{\mathrm{ev}}}
\newcommand\sgn{\mathord{\mathrm{sgn}}}
\newcommand\rk{\mathop{\mathrm{rk}}\nolimits}
\newcommand\im{\mathop{\mathrm{im}\;}\nolimits}
\newcommand\sh{\mathop{\mathrm{sh}}\nolimits}
\newcommand\wt{\mathop{\mathrm{wt}}\nolimits}
\newcommand\Mat{\mathord{\mathrm{Mat}}}
\newcommand\colr{\mathord{\mathbf{cr}}}
\newcommand\seqr{\mathord{\mathbf{sr}}}
\newcommand\Part{\mathord{\mathrm{Part}}}
\newcommand\llangle{\langle\!\langle}
\newcommand\rrangle{\rangle\!\rangle}
\makeatletter
\def\mysmash{\relax
  \ifmmode\expandafter\mathpalette\expandafter\mymathsm@sh
  \else\expandafter\mymakesm@sh\fi}
\def\mymakesm@sh#1{\setbox\z@\hbox{\color@begingroup#1\color@endgroup}%
  \myfinsm@sh}
\def\mymathsm@sh#1#2{\setbox\z@\hbox{$\m@th#1{#2}$}\myfinsm@sh}
\def\myfinsm@sh{\wd\z@\z@\box\z@}
\makeatother
\renewcommand\labelitemi{$\ \bullet$}
\renewcommand\theenumi{(\roman{enumi})}
\renewcommand\labelenumi{\theenumi}

\renewcommand\enumerate{\list\labelenumi
  {\setlength\leftmargin{0pt}\setlength\labelwidth{0pt}
  \usecounter{enumi}\def\makelabel##1{\kern\labelsep{##1}}}}
\renewcommand\itemize{\list\labelitemi
  {\setlength\leftmargin{0pt}\setlength\labelwidth{0pt}
  \def\makelabel##1{\kern\labelsep{##1}}}}
\renewenvironment{description}{\list{}
  {\setlength\leftmargin{0pt}\setlength\labelwidth{0pt}
  \let\makelabel\descriptionlabel}}
  {\endlist}
\newcommand\altfootnotetext[1]{\insert\footins{\normalfont\footnotesize
  \interlinepenalty\interfootnotelinepenalty\splittopskip\footnotesep
  \splitmaxdepth\dp\strutbox\hsize\columnwidth\rule{0pt}\footnotesep
  \ignorespaces#1\par}}
\newtheorem{theorem}{Theorem}
\newtheorem{lemma}[theorem]{Lemma}
\newtheorem{proposition}[theorem]{Proposition}
\newtheorem{corollary}[theorem]{Corollary}
\newenvironment{proof}{\trivlist
  \item[\hskip \labelsep{\itshape Proof.}]\upshape}{\nobreak\noindent$
  \square$\endtrivlist}
\newenvironment{other}[1]{\refstepcounter{theorem}\trivlist
  \item[\hskip \labelsep{\itshape #1~\arabic{theorem}.}]\upshape}
  {\endtrivlist\bigbreak}
\SelectTips{cm}{10}
\UseTips
\begin{document}
\title{A Solomon descent theory for the wreath products $G\wr\mathfrak S_n$}
\author{Pierre Baumann and Christophe Hohlweg\thanks{Partially supported
by Canada Research Chairs.}}
\date{}
\maketitle
\altfootnotetext{MSC: Primary 16S99, Secondary 05E05 05E10 16S34 16W30
20B30 20E22.}
\altfootnotetext{Keywords: wreath products, Solomon descent algebra,
quasisymmetric functions.}
\begin{abstract}
We propose an analogue of Solomon's descent theory for the case of a
wreath product $G\wr\mathfrak S_n$, where $G$ is a finite abelian
group. Our construction mixes a number of ingredients:
Mantaci-Reutenauer algebras, Specht's theory for the representations
of wreath products, Okada's extension to wreath products of the
Robinson-Schensted correspondence, Poirier's quasisymmetric functions.
We insist on the functorial aspect of our definitions and explain the
relation of our results with previous work concerning the hyperoctaedral
group.
\end{abstract}

\section*{Introduction}
The problem studied in this article has its roots in a discovery by
Solomon in 1976. Let $(W,(s_i)_{i\in I})$ be a Coxeter system. For any
subset $J\subseteq I$, call $W_J$ the parabolic subgroup generated by
the elements $s_j$ with $j\in J$. In each left coset $wW_J$ of $W$
modulo $W_J$, there is a unique element of minimal length, called the
distinguished representative of that coset. We denote the set of these
distinguished representatives by $X_J$, and we form the sum
$x_J=\sum_{w\in X_J}w$ in the group ring $\mathbb ZW$. Finally we
denote by $\Sigma_W$ the $\mathbb Z$-submodule of $\mathbb ZW$ spanned
by all elements $x_J$.

Now let $R(W)$ be the character ring of $W$, and let $\varphi_J\in R(W)$
be the character of $W$ induced from the trivial character of $W_J$.
Given two subsets $J$ and $K$ of $I$, each double coset $C\in
W_J\backslash W/W_K$ contains a unique element $x$ of minimal length,
and a result of Tits, Kilmoyer~\cite{Kilmoyer69} and/or Solomon
\cite{Solomon76} asserts that the intersection $x^{-1}W_Jx\cap W_K$ is
the parabolic subgroup $W_{L(C)}$, where $L(C)=\{k\in K\mid\exists j\in
J,\ x^{-1}s_jx=s_k\}$. Joint to Mackey's tensor product theorem, this
yields the multiplication rule in the representation ring $R(W)$
$$\varphi_J\varphi_K=\sum_{L\subseteq I}a_{JKL}\,\varphi_L,\quad
\text{where}\quad a_{JKL}=\bigl|\{C\in W_J\backslash W/W_K\mid
L=L(C)\}\bigr|.$$

With these notations, Solomon's discovery \cite{Solomon76} is the
equality $x_Jx_K=\sum_{L\subseteq I}a_{JKL}x_L$ in the ring
$\mathbb ZW$. It implies that $\Sigma_W$ is a subring of $\mathbb ZW$
and it shows the existence a morphism of rings~$\theta_W:\Sigma_W\to R(W)$
such that $\theta_W(x_J)=\varphi_J$. This result means that (a part of)
the character theory of $W$ can be lifted to a subring of its group ring.
Additional details (for instance, a more precise description of the image
of $\theta_W$) can be found in the paper \cite{BBHT92} by F.~Bergeron,
N.~Bergeron, Howlett and Taylor.

It is natural to look for a similar theory for groups other than Coxeter
systems. The first examples that come to mind are finite groups of Lie
type and finite complex reflection groups. Among the latter, the groups
of type $G(r,1,n)$ are wreath products $(\mathbb Z/r\mathbb Z)\wr
\mathfrak S_n$ of a cyclic group $\mathbb Z/r\mathbb Z$ by the symmetric
group $\mathfrak S_n$. One is then led to investigate the case of a
general wreath product $G\wr\mathfrak S_n$. To build the theory, it is
necessary to have some knowledge about the representation theory of $G$
itself, and we assume in this paper that $G$ is abelian. One of our
main results explains how to construct a subring $\MR_n(\mathbb ZG)$
inside the group ring $\mathbb Z\bigl[G\wr\mathfrak S_n]$ and a
surjective ring homomorphism $\theta_G$ from $\MR_n(\mathbb ZG)$ onto
the representation ring $R(G\wr\mathfrak S_n)$ of the wreath product.
Here the notation $\MR$ refers to the names of Mantaci and Reutenauer;
indeed it turns out that the remarkable subring inside $\mathbb Z
\bigl[G\wr\mathfrak S_n]$ discovered in 1995 by these two authors
\cite{Mantaci-Reutenauer95} is adequate to our purpose.

A usually efficient method to tackle problems with the symmetric group
$\mathfrak S_n$ is to treat all $n$ at the same time. For instance,
Malvenuto and Reutenauer observed in 1995 \cite{Malvenuto-Reutenauer95}
that the direct sum $\mathscr F=\bigoplus_{n\geq0}\mathbb Z[\mathfrak
S_n]$ can be endowed \vspace{2pt}with the structure of a graded bialgebra
in such a way that the submodule $\Sigma=\bigoplus_{n\geq0}
\Sigma_{\mathfrak S_n}$ is a graded subbialgebra. A similar phenomenon
appears here: the direct sum $\mathscr F(\mathbb ZG)=\bigoplus_{n\geq0}
\mathbb Z\bigl[G\wr\mathfrak S_n\bigr]$ can be endowed with the structure
of a graded bialgebra, of which $\MR(\mathbb ZG)=\bigoplus_{n\geq0}
\MR_n(\mathbb ZG)$ is a subbialgebra. (A particular case of this
construction was previously considered by Aguiar and Mahajan; the paper
\cite{Aguiar-Bergeron-Nyman04} by Aguiar, N.~Bergeron and Nyman presents an
account of their result. Aguiar and his coauthors view the hyperoctaedral
group of order $2^nn!$ as the wreath product $\{\pm1\}\wr\mathfrak S_n$,
that is, as the group of signed permutations. Then they construct the
graded bialgebra $\mathscr F\bigl(\mathbb Z\bigl[\{\pm1\}\bigr]\bigr)$ and
its subbialgebra $\MR\bigl(\mathbb Z\bigl[\{\pm1\}\bigr]\bigr)$. Using the
morphism of group `forgetting the signs' from $\{\pm1\}\wr\mathfrak S_n$
onto $\mathfrak S_n$, they compare these graded bialgebras with Malvenuto
and Reutenauer's bialgebra $\mathscr F$ and its subbialgebra $\Sigma$.
Our construction and its functoriality generalize Aguiar and his coauthors'
results to the case of all wreath products $G\wr\mathfrak S_n$.) This
bialgebra structure on $\mathscr F(\mathbb ZG)$ will be the starting point
of our story; indeed we define a `free quasisymmetric algebra' $\mathscr
F(V)$ for any $\mathbb Z$-module $V$ and investigate its properties.

We now present the plan and the main results of this paper.

In Section~\ref{se:GenConstr}, we define the free quasisymmetric
algebra on a module $V$ over a commutative ground ring $\mathbb K$: this
is a graded module $\mathscr F(V)=\bigoplus_{n\geq0}\mathscr F_n(V)$,
which we endow with an `external product' and a coproduct to turn it
into a graded bialgebra (Theorem~\ref{th:FGIsBialgebra}). In the case
where $V$ is endowed with the structure of a coalgebra, $\mathscr F(V)$
contains a remarkable subbialgebra $\MR(V)$, the so-called
Mantaci-Reutenauer bialgebra, which is a free associative algebra as
soon as $V$ is a free module (Propositions~\ref{pr:SigmaIsFreeAlg} and
\ref{pr:SigmaWIsSubbialg}).

In Section~\ref{se:Duality}, we show that the functor $V\rightsquigarrow
\mathscr F(V)$ is compatible with the duality of $\mathbb K$-modules,
in the sense that any pairing between two $\mathbb K$-modules $V$ and
$W$ gives rise to a pairing of bialgebras between $\mathscr F(V)$ and
$\mathscr F(W)$ (Proposition~\ref{pr:FAndDuality}). In particular, the
bialgebra $\mathscr F(V)$ is self-dual as soon as the module $V$ is
endowed with a perfect pairing.

In Section~\ref{se:InternProd}, we investigate the case where the module
$V$ is a $\mathbb K$-algebra. Then $\mathscr F(V)$ can be endowed with an
`internal product', which turns each of the graded components $\mathscr
F_n(V)$ into an algebra. The interesting point here is the existence of
a splitting formula that describes the compatibility between this
internal product, the external product and the coproduct
(Theorem~\ref{th:SplitForm}). This formula is a generalization of the
splitting formula of Gelfand, Krob, Lascoux, Leclerc, Retakh and
Thibon~\cite{GKLLRT95}; it entails that the Mantaci-Reutenauer
bialgebra $\MR(V)$ is a subalgebra of $\mathscr F(V)$ for the internal
product whenever $V$ is endowed with the structure of a cocommutative
bialgebra (Corollaries~\ref{co:SubalgIntSigmaB} and~\ref{co:MantReutRule}).
In Section~\ref{ss:CaseGpAlg}, we consider for $V$ the case of the
group algebra $\mathbb K\Gamma$ of a finite group $\Gamma$ and justify
that the graded component $\mathscr F_n(\mathbb K\Gamma)$ is canonically
isomorphic to the group algebra $\mathbb K\bigl[\Gamma\wr\mathfrak
S_n\bigr]$, and that the graded component $\MR_n(\mathbb K\Gamma)=
\MR(\mathbb K\Gamma)\cap\mathscr F_n(\mathbb K\Gamma)$ coincides with
the subalgebra defined by Mantaci and Reutenauer in
\cite{Mantaci-Reutenauer95}.

In Section~\ref{se:SoloTheoWrProd}, we at last provide the link between
these constructions and a Solomon descent theory for wreath products.
We first recall Specht's classification of the irreducible complex
characters of a wreath product $G\wr\mathfrak S_n$ and Zelevinsky's
structure of a graded bialgebra on the direct sum $\Rep(G)=
\bigoplus_{n\geq0}R(G\wr\mathfrak S_n)$ for the induction product and
the restriction coproduct (Section~\ref{ss:CharacWreathProd}). We then
focus on the case where $G$ is abelian. We denote the dual group of
$G$ by $\Gamma$, we observe that the group ring $\mathbb Z\Gamma$ is
a cocommutative bialgebra, so that the Mantaci-Reutenauer bialgebra
$\MR(\mathbb Z\Gamma)$ is defined and is a subalgebra of $\mathscr F
(\mathbb Z\Gamma)$ for the internal product, and we define a map
$\theta_G:\MR(\mathbb Z\Gamma)\to\Rep(G)$. Then we show that $\theta_G$
is a surjective morphism of graded bialgebras, and that in each degree,
$\theta_G:\MR_n(\mathbb Z\Gamma)\to R(G\wr\mathfrak S_n)$ is a surjective
morphism of rings whose kernel is the Jacobson radical of
$\MR_n(\mathbb Z\Gamma)$ (Theorem~\ref{th:SoloHomo}). We also show
that $\theta_G$ enjoys a remarkable symmetry property analogous
to the symmetry property of Solomon's homomorphisms $\theta_W$ proved
by Jöllenbeck and Reutenauer \cite{Jollenbeck-Reutenauer01} and by
Blessenohl, Hohlweg and Schocker \cite{Blessenohl-Hohlweg-Schocker04}
(Theorem~\ref{th:SymmSoloHomo}). Finally we compare our results with
the work of Bonnafé and Hohlweg, who treated in \cite{Bonnafe-Hohlweg04}
the case of the hyperoctaedral group $\{\pm1\}\wr\mathfrak S_n$ using
methods from the theory of Coxeter groups (Section~\ref{ss:PartCaseZ/2Z}).

The questions about the bialgebras $\mathscr F(V)$ investigated in
Sections~\ref{se:GenConstr} to \ref{se:InternProd} are functorial
in the $\mathbb K$-module $V$. As usual, the most interesting point
in this assertion is the compatibility of the constructions with the
homomorphisms, namely here the $\mathbb K$-linear maps. On the contrary
the questions studied in Section~\ref{se:ColCombHopfAlg} require that
$V$ be a free $\mathbb K$-module and depend on the choice of a basis
$B$ of $V$. Such a basis $B$ can be viewed as the data of a structure
of a pointed coalgebra on $V$, which yields in turn a Mantaci-Reutenauer
subbialgebra $\MR(V)$ inside $\mathscr F(V)$. The choice of $B$ also
gives rise to a second subbialgebra $\mathscr Q(B)$, bigger than
$\MR(V)$, which we call the coplactic bialgebra. The definition of
$\mathscr Q(B)$ involves a combinatorial construction due to
Okada~\cite{Okada90}, which extends the well-known Robinson-Schensted
correspondence to `coloured' situations; at this point, we take the
opportunity to provide an analogue of Knuth relations for Okada's
correspondence (Proposition~\ref{pr:KnuthRSORel}). In the case where $B$
is a singleton set, the bialgebra $\mathscr Q(B)$ is one of the `algèbres
de Hopf de tableaux' of Poirier and Reutenauer \cite{Poirier-Reutenauer95}.
Extending the work of these authors, we define a surjective homomorphism
$\Theta_B$ of graded bialgebras from $\mathscr Q(B)$ onto a
bialgebra $\Lambda(B)$ of `coloured' symmetric functions
(Theorem~\ref{th:ThetaCoplac}). We then go back to the situation
investigated in Section~\ref{se:SoloTheoWrProd} and take the group
algebra $\mathbb Z\Gamma$ for $V$ and the group $\Gamma$ for $B$;
here $\Theta_\Gamma$ can be viewed as a lift of $\theta_G:\MR(\mathbb
Z\Gamma)\to\Rep(G)$ to $\mathscr Q(\Gamma)$ that yields a nice
description of the simple representations of all wreath products
$G\wr\mathfrak S_n$. We recover Jöllenbeck's construction of the Specht
modules \cite{Jollenbeck99} as the particular case where $G$ is the
group with one element; we refer the reader to Blessenohl and Schocker's
survey \cite{Blessenohl-Schocker02} for additional details about
Jöllenbeck's construction.

Finally we present in Section~\ref{se:ColQSymFun} a realization of the
bialgebra $\mathscr F(V)$ in terms of free quasisymmetric functions.
As in Section~\ref{se:ColCombHopfAlg}, the $\mathbb K$-module $V$ is
assumed to be free; we choose a basis $B$ of $V$ and endow $B$ with a
linear order. When $V$ has rank one, our free quasisymmetric functions
coincide with the usual ones \cite{Gessel84}. In higher rank however,
our free quasisymmetric functions are different from those defined by
Novelli and Thibon in~\cite{Novelli-Thibon04}. This disagreement has
its roots in the fact that Novelli and Thibon's construction and ours
were designed with different aims: roughly speaking, Novelli and
Thibon's goal was to find a noncommutative version of Poirier's
quasisymmetric functions \cite{Poirier98}; on the other side, we
view the dual algebra $\MR(V)^\vee$ as a quotient of $\mathscr F(V)$
and describe it in terms of commutative quasisymmetric functions.

At this point, we should mention that the assignment $(V,B)
\rightsquigarrow\MR(V)^\vee$ enjoys a certain functoriality property;
this property and the isomorphism between $\MR(\mathbb K)^\vee$ and
the graded bialgebra $\QSym$ of usual quasisymmetric functions yield
in turn homomorphisms of graded bialgebras from $\mathscr F(V)$ and
$\MR(V)^\vee$ to $\QSym$, which amounts to say that $\mathscr F(V)$
and $\MR(V)^\vee$ are `combinatorial Hopf algebras' in the sense of
Aguiar, N.~Bergeron and Sottile \cite{Aguiar-Bergeron-Sottile04}.

The authors wish to thank Jean-Christophe~Novelli and Jean-Yves~Thibon for
fruitful and instructive conversations, which took place on March 30, 2004
in Ottrott and on May 3, 2004 at the Institut Gaspard Monge (University of
Marne-la-Vallée). Their preprint~\cite{Novelli-Thibon04} influenced our
writing of Sections~\ref{se:GenConstr} and \ref{se:InternProd}. The main
part of this work was carried out when C.~H.\ was at the Institut de
Recherche Mathématique Avancée in Strasbourg.

We fix a commutative ground ring $\mathbb K$. Connected $\mathbb N$-graded
$\mathbb K$-bialgebras appear everywhere in the paper. Such bialgebras are
indeed automatically Hopf algebras, at least when $\mathbb K$ is a field.
However we will neither make use of this property nor attempt to work out
explicitly any antipode.

\section{Free quasisymmetric bialgebras}
\label{se:GenConstr}
In this section, we present our main objects of study, namely the
free quasisymmetric bialgebras and the generalized descent algebras,
among which the Novelli-Thibon bialgebras and the Mantaci-Reutenauer
bialgebras. Before that, we introduce some notations pertaining
permutations.

\subsection{Notations related to permutations}
\label{ss:NotaRelaPerm}
For each positive integer $n$, we denote the symmetric group of all
permutations of the set $\{1,2,\ldots,n\}$ by $\mathfrak S_n$. By
convention, $\mathfrak S_0$ is the group with one element. The unit
element of $\mathfrak S_n$ is denoted by $e_n$. The group algebra over
$\mathbb K$ of $\mathfrak S_n$ is denoted by $\mathbb K\mathfrak S_n$.
In practice, a permutation $\sigma\in\mathfrak S_n$ is written as the
word $\sigma(1)\sigma(2)\cdots\sigma(n)$ with letters in $\mathbb
Z_{>0}=\{1,2,\ldots\}$.

Let $\mathscr A$ be totally ordered set (an alphabet). The standardization
of a word $w=a_1a_2\cdots a_n$ of length $n$ with letters in $\mathscr A$
is the permutation $\sigma\in\mathfrak S_n$ with smallest number of
inversions such that the sequence
$$\bigl(a_{\sigma^{-1}(1)},a_{\sigma^{-1}(2)},\ldots,a_{\sigma^{-1}(n)}
\bigr)$$
is non-decreasing. In other words, the word $\sigma(1)\sigma(2)\cdots
\sigma(n)$ that represents $\sigma$ is obtained by putting the numbers
$1$, $2$, \dots, $n$ in the place of the letters $a_i$ of $w$; in this
process of substitution, the diverse occurrences of the smallest letter
of $\mathscr A$ get replaced first by the numbers $1$, $2$, etc.\ from
left to right; then we replace the occurrences of the second-smallest
element of $\mathscr A$ by the following numbers; and so on, up to the
exhaustion of all letters of $w$. An example clarifies this explanation:
given the alphabet $\mathscr A=\{a,b,c,\ldots\}$ with the usual order,
the standardization of the word $w=bcbaba$ is $\sigma=364152$.

A composition of a positive integer $n$ is a sequence $\mathbf c=(c_1,
c_2,\ldots,c_k)$ of positive integers which sum up to $n$. The usual
notation for that is to write $\mathbf c\models n$. Given two
compositions $\mathbf c=(c_1,c_2,\ldots,c_k)$ and $\mathbf d=(d_1,d_2,
\ldots,d_l)$ of the same integer $n$, we say that $\mathbf c$ is a
refinement of $\mathbf d$ and we write $\mathbf c\succcurlyeq\mathbf d$
if there holds
$$\{c_1,c_1+c_2,\ldots,c_1+c_2+\cdots+c_{k-1}\}\supseteq\{d_1,d_1+d_2,
\ldots,d_1+d_2+\cdots+d_{l-1}\}.$$
The relation $\preccurlyeq$ is a partial order on the set of compositions
of $n$. For instance, the following chain of inequalities hold among
compositions of $5$:
$$(5)\prec(4,1)\prec(1,3,1)\prec(1,2,1,1)\prec(1,1,1,1,1).$$

Let $\mathbf c=(c_1,c_2,\ldots,c_k)$ be a composition of $n$ and set
$t_i=c_1+c_2+\cdots+c_i$ for each $i$. Given a $k$-uple $(\sigma_1,
\sigma_2,\ldots,\sigma_k)\in\mathfrak S_{c_1}\times\mathfrak S_{c_2}
\times\cdots\times\mathfrak S_{c_k}$ of permutations, we define
$\sigma_1\times\sigma_2\times\cdots\times\sigma_k\in\mathfrak S_n$
as the permutation that maps an element $a$ belonging to the interval
$[t_{i-1}+1,t_i]$ onto $t_{i-1}+\sigma_i(a-t_{i-1})$. This assignment
defines an embedding $\mathfrak S_{c_1}\times\mathfrak S_{c_2}\times\cdots
\times\mathfrak S_{c_k}\hookrightarrow\mathfrak S_n$; we denote its
image by $\mathfrak S_{\mathbf c}$. Such a $\mathfrak S_{\mathbf c}$
is called a Young subgroup of $\mathfrak S_n$. We obtain for free an
embedding for the group algebras
$$\mathbb K\mathfrak S_{c_1}\otimes\mathbb K\mathfrak S_{c_2}\otimes
\cdots\otimes\mathbb K\mathfrak S_{c_k}\stackrel\simeq\longrightarrow
\mathbb K\mathfrak S_{\mathbf c}\subseteq\mathbb K\mathfrak S_n.$$
The map $\mathbf c\mapsto\mathfrak S_{\mathbf c}$ is an order reversing
bijection from the set of compositions of $n$, endowed with the refinement
order, onto the set of Young subgroups of $\mathfrak S_n$, endowed with
the inclusion order.

Let again $\mathbf c=(c_1,c_2,\ldots,c_k)$ be a composition of $n$ and
set $t_i=c_1+c_2+\cdots+c_i$. The subset
$$X_{\mathbf c}=\bigl\{\sigma\in\mathfrak S_n\bigm|\forall i,\ \sigma
\text{ is increasing on the interval }[t_{i-1}+1,t_i]\bigr\}$$
is a system of representatives of the left cosets of $\mathfrak
S_{\mathbf c}$ in $\mathfrak S_n$. Here are some examples:
$$X_{(2,2)}=\{1234,1324,1423,2314,2413,3412\},\quad X_{(n)}=\{\id\}
\quad\text{and}\quad X_{(\underbrace{1,1,\ldots,1}_{n\text{ times}})}
=\mathfrak S_n.$$
We define an element of the group ring $\mathbb K\mathfrak S_n$ by
setting $x_{\mathbf c}=\sum_{\sigma\in X_{\mathbf c}}\sigma$.

Let $\mathbf d=(d_1,d_2,\ldots,d_l)$ be a composition of an integer $n$.
Then a composition $\mathbf c$ of $n$ is a refinement of $\mathbf d$ if
and only if $\mathbf c$ can be obtained as the concatenation $\mathbf
f_1\mathbf f_2\cdots\mathbf f_l$ of a composition $\mathbf f_1$ of $d_1$,
a composition $\mathbf f_2$ of $d_2$, \dots, and a composition $\mathbf
f_l$ of $d_l$. If this holds, then the map
$$(\rho,\sigma_1,\sigma_2,\ldots,\sigma_l)\mapsto\rho\circ(\sigma_1
\times\sigma_2\times\cdots\times\sigma_l)$$
is a bijection from $X_{\mathbf d}\times X_{\mathbf f_1}\times
X_{\mathbf f_2}\times\cdots\times X_{\mathbf f_l}$ onto $X_{\mathbf c}$,
for $X_{\mathbf f_1}\times\cdots\times X_{\mathbf f_l}$ is a set of
minimal coset representatives of $\mathfrak S_{\mathbf c}$ in $\mathfrak
S_{\mathbf d}$. Therefore the equality
\begin{equation}
x_{\mathbf c}=x_{\mathbf d}\;(x_{\mathbf f_1}\otimes x_{\mathbf f_2}
\otimes\cdots\otimes x_{\mathbf f_l})
\label{eq:PropicXc}
\end{equation}
holds in the group ring $\mathbb K\mathfrak S_n$. As a particular case of
(\ref{eq:PropicXc}), we see that
\begin{equation}
x_{(n,n',n'')}=x_{(n,n'+n'')}\bigl(x_{(n)}\otimes x_{(n',n'')}\bigr)=
x_{(n+n',n'')}\bigl(x_{(n,n')}\otimes x_{(n'')}\bigr)
\label{eq:AssoExtProd}
\end{equation}
holds true for any three positive integers $n$, $n'$ and $n''$.

Let $\sigma\in\mathfrak S_n$. One may partition the word $\sigma(1)
\sigma(2)\cdots\sigma(n)$ that represents $\sigma$ into its longest
increasing subwords; the composition of $n$ formed by the successive
lengths of these subwords is called the descent composition of $\sigma$
and is denoted by $D(\sigma)$. For instance, the descent composition of
$\sigma=51243$ is $D(\sigma)=(1,3,1)$. Then for any composition
$\mathbf c$ of $n$, the assertions $\sigma\in X_{\mathbf c}$ and
$D(\sigma)\preccurlyeq\mathbf c$ are equivalent.

\subsection{Definition of the free quasisymmetric bialgebra
$\mathscr F(V)$}
\label{ss:FreeBialgFV}
Let $V$ be a $\mathbb K$-module. The group $\mathfrak S_n$ acts on the
$n$-th tensor power $V^{\otimes n}$; the submodule of invariants, that
is, the space of symmetric tensors, is denoted by $\TS^n(V)$. We may
form the tensor product of $V^{\otimes n}$ by $k\mathfrak S_n$. To
distinguish this tensor product from those used to build the tensor power
$V^{\otimes n}$, we denote it with a sharp symbol. We denote the result
$(V^{\otimes n})\#(\mathbb K\mathfrak S_n)$ by $\mathscr F_n(V)$. The
actions defined by
$$\pi\cdot\bigl[(v_1\otimes v_2\otimes\cdots\otimes v_n)\#\sigma\bigr]=
\bigl[(v_{\pi^{-1}(1)}\otimes v_{\pi^{-1}(2)}\otimes\cdots\otimes
v_{\pi^{-1}(n)})\#(\pi\sigma)\bigr]$$
and
$$\bigl[(v_1\otimes v_2\otimes\cdots\otimes v_n)\#\sigma\bigr]\cdot\pi=
\bigl[(v_1\otimes v_2\otimes\cdots\otimes v_n)\#(\sigma\pi)\bigr]$$
endow $\mathscr F_n(V)$ with the structure of a $\mathbb K\mathfrak
S_n$-bimodule, where $(v_1,v_2,\ldots,v_n)\in V^n$ and $\pi\in\mathfrak
S_n$. For instance, $\mathscr F_n(\mathbb K)$ is the (left and right)
regular $\mathbb K\mathfrak S_n$-module.

Our aim now is to endow the space $\mathscr F(V)=\bigoplus_{n\geq0}
\mathscr F_n(V)$ with the structure of a graded bialgebra. We define
the product of two elements $\alpha\in\mathscr F_n(V)$ and
$\alpha'\in\mathscr F_{n'}(V)$ of the form
$$\alpha=\bigl[(v_1\otimes v_2\otimes\cdots\otimes v_n)\#\sigma\bigr]
\quad\text{and}\quad\alpha'=\bigl[(v'_1\otimes v'_2\otimes\cdots\otimes
v'_{n'})\#\sigma'\bigr]$$
by the formula
$$\alpha*\alpha'=x_{(n,n')}\cdot\bigl[(v_1\otimes v_2\otimes\cdots\otimes
v_n\otimes v'_1\otimes v'_2\otimes\cdots\otimes v'_{n'})\#(\sigma\times
\sigma')\bigr].$$
(This formula can be made more concrete by noting that $x_{(n,n')}\;
(\sigma\times\sigma')$ is the sum in the group algebra $\mathbb K\mathfrak
S_{n+n'}$ of all permutations $\pi$ such that $\sigma$ is the
standardization of the word $\pi(1)\pi(2)\cdots\pi(n)$ and $\sigma'$ is
the standardization of the word $\pi(n+1)\pi(n+2)\cdots\pi(n+n')$.)
We extend this definition by multilinearity to an operation defined on
the whole space $\mathscr F(V)$ and call this latter the external product.

We define the coproduct of an element $\alpha=\bigl[(v_1\otimes v_2\otimes
\cdots\otimes v_n)\#\sigma\bigr]$ of $\mathscr F_n(V)$ as
$$\Delta\bigl((v_1\otimes v_2\otimes\cdots\otimes v_n)\#\sigma\bigr)=
\sum_{n'=0}^n\bigl[(v_1\otimes v_2\otimes\cdots\otimes v_{n'})\#\pi_{n'}
\bigr]\otimes\bigl[(v_{n'+1}\otimes v_{n'+2}\otimes\cdots\otimes v_n)\#
\pi'_{n-n'}\bigr],$$
where $\pi_{n'}\in\mathfrak S_{n'}$ is the inverse of the standardization
of the word $\sigma^{-1}(1)\;\sigma^{-1}(2)\;\cdots\;\sigma^{-1}(n')$
and $\pi'_{n-n'}\in\mathfrak S_{n-n'}$ is the inverse of the
standardization of the word $\sigma^{-1}(n'+1)\;\sigma^{-1}(n'+2)\;
\cdots\;\sigma^{-1}(n)$. In other words, $\pi_{n'}$ and $\pi'_{n-n'}$ are
such that the two sequences of letters $(\pi_{n'}(1),\pi_{n'}(2),\ldots,
\pi_{n'}(n'))$ and $(n'+\pi'_{n-n'}(1),n'+\pi'_{n-n'}(2),\ldots,n'+
\pi'_{n-n'}(n-n'))$ appear in this order in the word $\sigma(1)\sigma(2)
\cdots\sigma(n)$. We call the map $\Delta:\mathscr F(V)\to\mathscr F(V)
\otimes\mathscr F(V)$ the coproduct of $\mathscr F(V)$.

We define the unit of $\mathscr F(V)$ as the injection of the graded
component $\mathscr F_0(V)=\mathbb K$ into $\mathscr F(V)$; we define
the counit of $\mathscr F(V)$ as the projection of $\mathscr F(V)$
onto $\mathscr F_0(V)=\mathbb K$.

We now give an example to illustrate these definitions. Given six
elements $v_1$, $v_2$, $v_3$, $v_4$, $v'_1$, $v'_2$ in $V$, the product
of $\alpha=\bigl[(v_1\otimes v_2)\#e_2\bigr]$ and $\alpha'=
\bigl[(v'_2\otimes v'_1)\#21\bigr]=(21)\cdot\bigl[(v'_1\otimes v'_2)
\#e_2\bigr]$ is
\begin{align*}
\alpha*\alpha'&=(1243+1342+1432+2341+2431+3421)\cdot\bigl[(v_1\otimes
v_2\otimes v'_1\otimes v'_2)\#e_4\bigr]\\[4pt]
&=\bigl[(v_1\otimes v_2\otimes v'_2\otimes v'_1)\#1243\bigr]+
\bigl[(v_1\otimes v'_2\otimes v_2\otimes v'_1)\#1342\bigr]+\\[2pt]
&\quad\ \bigl[(v_1\otimes v'_2\otimes v'_1\otimes v_2)\#1432\bigr]+
\bigl[(v'_2\otimes v_1\otimes v_2\otimes v'_1)\#2341\bigr]+\\[2pt]
&\quad\ \bigl[(v'_2\otimes v_1\otimes v'_1\otimes v_2)\#2431\bigr]+
\bigl[(v'_2\otimes v'_1\otimes v_1\otimes v_2)\#3421\bigr],
\end{align*}
and the coproduct of $\alpha=\bigl[(v_3\otimes v_1\otimes
v_2\otimes v_4)\#2314\bigr]=(2314)\cdot\bigl[(v_1\otimes v_2\otimes
v_3\otimes v_4)\#e_4\bigr]$ is
\begin{align*}
\Delta(\alpha)&=\bigl[()\#e_0\bigr]\otimes\alpha+\bigl[(v_3)\#1\bigr]
\otimes\bigl[(v_1\otimes v_2\otimes v_4)\#123\bigr]+\\[2pt]
&\quad\ \bigl[(v_3\otimes v_1)\#21\bigr]\otimes\bigl[(v_2\otimes v_4)
\#12\bigr]+\\[2pt]
&\quad\ \bigl[(v_3\otimes v_1\otimes v_2)\#231\bigr]\otimes
\bigl[(v_4)\#1\bigr]+\alpha\otimes\bigl[()\#e_0\bigr]\\[4pt]
&=\bigl[()\#e_0\bigr]\otimes\alpha+\bigl[(v_3)\#e_1\bigr]\otimes
\bigl[(v_1\otimes v_2\otimes v_4)\#e_3\bigr]+\\[2pt]
&\quad\ +(21)\cdot\bigl[(v_1\otimes v_3)\#e_2\bigr]\otimes
\bigl[(v_2\otimes v_4)\#e_2\bigr]\\[2pt]
&\quad\ +(231)\cdot\bigl[(v_1\otimes v_2\otimes v_3)\#e_3\bigr]\otimes
\bigl[(v_4)\#e_1\bigr]+\alpha\otimes\bigl[()\#e_0\bigr].
\end{align*}

\begin{theorem}
\label{th:FGIsBialgebra}
The unit, the counit, and the operations $*$ and $\Delta$ endow $\mathscr
F(G)$ with the structure of a graded bialgebra.
\end{theorem}
\begin{proof}
It is clear that the four operations respect the graduation. The
associativity of $*$ follows immediately from Equation
(\ref{eq:AssoExtProd}). A moment's thought suffices to check the
coassociativity of $\Delta$ and the axioms for the unit and the counit.
It remains to show the pentagon axiom, which asks that $\Delta$
be multiplicative with respect to the product $*$.

Following Malvenuto and Reutenauer's method \cite{Malvenuto-Reutenauer95},
we first recall a classical construction in the theory of Hopf algebras.
Let $\mathscr A$ be a set, let $\langle\mathscr A\rangle$ denote the set
of words on $\mathscr A$, and let $\mathbb K\langle\mathscr A\rangle$ be
the free $\mathbb K$-module with basis $\langle\mathscr A\rangle$. The
shuffle product of two words $w$ and $w'$ of length $n$ and $n'$
respectively is the sum
$$w\sha w'=\sum_{\rho\in X_{(n,n')}}b_{\rho^{-1}(1)}b_{\rho^{-1}(2)}
\cdots b_{\rho^{-1}(n+n')},$$
where the word $b_1b_2\cdots b_{n+n'}$ is the concatenation of the
words $w$ and $w'$. This operation $\sha$ is then extended bilinearly
to a product on $\mathbb K\langle\mathscr A\rangle$. The deconcatenation
is the coproduct $\delta$ on $\mathbb K\langle\mathscr A\rangle$ such that
$$\delta(w)=\sum_{n'=0}^na_1a_2\cdots a_{n'}\otimes a_{n'+1}a_{n'+2}
\cdots a_n$$
for any word $w=a_1a_2\cdots a_n$. It is known that the operations
$\sha$ and $\delta$ endow $\mathbb K\langle\mathscr A\rangle$ with the
structure of a bialgebra (see Proposition~1.9 in \cite{Reutenauer93}
for a proof).

We are now ready to show the pentagon axiom in the case where the
$\mathbb K$-module $V$ is free. We take a basis $B$ of $V$ and we set
$\mathscr A=\mathbb Z_{>0}\times B$. We observe that the elements
$(b_1\otimes b_2\otimes\cdots\otimes b_n)\#\sigma$ form a basis of
$\mathscr F_n(V)$, where $(b_1,b_2,\ldots,b_n)\in B^n$ and
$\sigma\in\mathfrak S_n$. We may thus define linear maps
$j_k:\mathscr F(G)\to\mathbb K\langle\mathscr A\rangle$ (depending
on the choice of a non-negative integer $k$) by mapping an element
$\alpha=\bigl[(b_1\otimes b_2\otimes\cdots\otimes b_n)\#\sigma\bigr]$
to $j_k(\alpha)=a_1a_2\cdots a_n$, where $a_i=(k+\sigma^{-1}(i),b_i)$.
In the other direction, we define a linear map
$s:\mathbb K\langle\mathscr A\rangle\to\mathscr F(V)$ as
follows: given a word $w=a_1a_2\cdots a_n$ with letters in $\mathscr A$,
we write $a_i=(p_i,b_i)$ and set $s(w)=(b_1\otimes b_2\otimes\cdots\otimes
b_n)\#\sigma$, where $\sigma$ is the inverse of the standardization of
the word $p_1p_2\cdots p_n$.

One easily checks that $s\circ j_k=\id_{\mathscr F(G)}$ and that
$(s\otimes s)\circ\delta=\Delta\circ s$. Moreover, let $w=a_1a_2\cdots a_n$
and $w'=a'_1a'_2\cdots a'_{n'}$ be two words with letters in $\mathscr A$.
If we write $a_i=(p_i,b_i)$ and $a'_i=(p'_i,b'_i)$, then $s(w\sha w')=
s(w)*s(w')$ as soon as every integer $p_i$ is strictly smaller than
every integer $p'_i$.

We now take $\alpha\in\mathscr F_n(G)$ and $\alpha'\in\mathscr F_{n'}(G)$.
We compute:
\begin{align*}
\Delta(\alpha*\alpha')
&=\Delta\Bigl[s\bigl(j_0(\alpha)\bigr)*s\bigl(j_n(\alpha')\bigr)
\Bigr]\\[4pt]
&=\bigl(\Delta\circ s\bigr)\bigl(j_0(\alpha)\sha j_n(\alpha')
\bigr)\\[4pt]
&=(s\otimes s)\Bigl[\delta\bigl(j_0(\alpha)\sha j_n(\alpha')\bigr)
\Bigr]\\[4pt]
&=(s\otimes s)\Bigl[\delta\bigl(j_0(\alpha)\bigr)\sha\delta\bigl(j_n
(\alpha')\bigr)\Bigr]\\[4pt]
&=\bigl[(s\otimes s)\circ\delta\circ j_0(\alpha)\bigr]*\bigl[
(s\otimes s)\circ\delta\circ j_n(\alpha')\bigr]\\[4pt]
&=\bigl[\Delta\circ s\circ j_0(\alpha)\bigr]*\bigl[\Delta\circ s
\circ j_n(\alpha')\bigr]\displaybreak[0]\\[4pt]
&=\Delta(\alpha)*\Delta(\alpha').
\end{align*}

This relation proves the pentagon axiom for $\mathscr F(V)$ in the case
where $V$ is a free $\mathbb K$-module. In the general case, we may find
a free $\mathbb K$-module $\tilde V$ and a surjective morphism of
$\mathbb K$-modules $f:\tilde V\to V$. Then $f$ induces a surjective
map from $\mathscr F(\tilde V)$ onto $\mathscr F(V)$ which is a
morphism of algebras and of coalgebras. Since the operations $*$ and
$\Delta$ on $\tilde V$ satisfy the pentagon axiom, their analogues on
$V$ satisfy also the pentagon axiom. This completes the proof of the
theorem.
\end{proof}

We note that the assignment $V\rightsquigarrow\mathscr F(V)$ is a
covariant functor from the category of $\mathbb K$-modules to the
category of $\mathbb N$-graded bialgebras over $\mathbb K$.

The algebras $\mathscr F(V)$ were also indirectly defined by Novelli
and Thibon; in \cite{Novelli-Thibon04}, they denote our $\mathscr
F(\mathbb K^l)$ by $\FQSym^{(l)}$ and state that it is a free
associative algebra, whence the name `free quasisymmetric bialgebras.'

\begin{other}{Remark}
Given a $\mathbb K$-module $V$, one can endow the direct sum
$\bigoplus_{n\geq0}V^{\otimes n}$ with two structures of a graded
bialgebra: the tensor algebra, denoted by $\T(V)$, and the cotensor
algebra, sometimes denoted by $\T^c(V)$. (The bialgebra $\mathbb K
\langle\mathscr A\rangle$ used in the proof of Theorem
\ref{th:FGIsBialgebra} is indeed the cotensor algebra on the free
$\mathbb K$-module $\mathbb K\mathscr A$ with basis $\mathscr A$.)
One checks easily that the maps
$$\iota:\T(V)\to\mathscr F(V),\ v_1\otimes v_2\otimes\cdots\otimes
v_n\mapsto\sum_{\sigma\in\mathfrak S_n}\sigma\cdot(v_1\otimes v_2\otimes
\cdots\otimes v_n\#e_n)$$
and
$$p:\mathscr F(V)\to\T^c(V),\ (v_1\otimes v_2\otimes\cdots\otimes
v_n\#\sigma)\mapsto v_1\otimes v_2\otimes\cdots\otimes v_n$$
are morphisms of graded bialgebras. Moreover the composition $p\circ\iota$
is the symmetrization map
$$\T(V)\to\T^c(V),\ v_1\otimes v_2\otimes\cdots\otimes v_n\mapsto
\sum_{\sigma\in\mathfrak S_n}v_{\sigma(1)}\otimes v_{\sigma(2)}\otimes
\cdots\otimes v_{\sigma(n)}.$$
For details and applications of this construction, we refer the reader
to \cite{Nichols78} and \cite{Rosso00}.
\end{other}

\subsection{The descent subbialgebras $\Sigma(W)$}
\label{ss:SubbialgSigmaW}
In this section, we investigate a class of graded subalgebras of
$\mathscr F(V)$, called the descent algebras. We find a criterion for
a descent algebra to be a subbialgebra of $\mathscr F(V)$ and give a
couple of examples.

We fix here a $\mathbb K$-module $V$. To any graded submodule
$W=\bigoplus_{n\geq0}W_n$ of the tensor algebra
$\T(V)=\bigoplus_{n\geq0}V^{\otimes n}$, we associate the subalgebra
$\Sigma(W)$ of $\mathscr F(V)$ generated by all elements of the form
$(t\#e_n)$ with $t\in W_n$. We call such a subalgebra $\Sigma(W)$ a
descent algebra. A descent algebra is necessarily graded, for it is
generated by homogeneous elements.

\begin{proposition}
\label{pr:SigmaIsFreeAlg}
Assume that $V$ is flat and that each module $W_n$ is free of finite
rank. For each $n\geq1$, pick a basis $B_n$ of $W_n$. Then $\Sigma(W)$
is the free associative algebra on the elements $(b\#e_n)$, where
$n\geq1$ and $b\in B_n$.
\end{proposition}
\begin{proof}
By the way of contradiction, we assume that there exists a finite
family $(\mathbf u_i)_{i\in I}$ consisting of distinct finite sequences
$\mathbf u_i=\bigl(\bigl(c^{(i)}_1,b^{(i)}_1\bigr),\bigl(c^{(i)}_2,
b^{(i)}_2\bigr),\ldots,\bigl(c^{(i)}_{k_i},b^{(i)}_{k_i}\bigr)\bigr)$
of elements in $\bigcup_{n\geq1}\bigl(\{n\}\times B_n\bigr)$ and a finite
family $(\lambda_i)_{i\in I}$ of elements of $\mathbb K\setminus\{0\}$
such that
\begin{equation}
\sum_{i\in I}\lambda_i\Bigl[\Bigl(b^{(i)}_1\#e_{c^{(i)}_1}\Bigr)*
\Bigl(b^{(i)}_2\#e_{c^{(i)}_2}\Bigr)*\cdots*\Bigl(b^{(i)}_{k_i}\#
e_{c^{(i)}_{k_i}}\Bigr)\Bigr]=0.
\label{eq:PfPrSIFAa}
\end{equation}
Using the graduation, we may suppose without loss of generality that
all the sequences $\mathbf c_i=(c^{(i)}_1,c^{(i)}_2,\ldots,c^{(i)}_{k_i})$
are compositions of the same integer $n$. Then (\ref{eq:PfPrSIFAa})
yields
\begin{equation}
\sum_{i\in I}\lambda_i\,x_{\mathbf c_i}\cdot\bigl[(b^{(i)}_1\otimes
b^{(i)}_2\otimes\cdots\otimes b^{(i)}_{k_i})\#e_n\bigr]=0.
\label{eq:PfPrSIFAb}
\end{equation}

We choose a maximal element $\mathbf c=(c_1,c_2,\ldots,c_k)$ among
the set $\{\mathbf c_i\mid i\in I\}$ with respect to the refinement
order, we set $J=\{i\in I\mid\mathbf c_i=\mathbf c\}$, and we choose
a permutation $\sigma\in\mathfrak S_n$ whose descent composition is
$\mathbf c$. Then for any $i\in I$,
$$\sigma\in X_{\mathbf c_i}\Longleftrightarrow\mathbf c\preccurlyeq
\mathbf c_i\Longleftrightarrow i\in J.$$
Taking the image of (\ref{eq:PfPrSIFAb}) by the linear map $p:\mathscr
F_n(V)\to V^{\otimes n}$ defined by
$$p\bigl((v_1\otimes v_2\otimes\cdots\otimes v_n)\#\rho\bigr)=
\begin{cases}v_{\rho(1)}\otimes v_{\rho(2)}\otimes\cdots\otimes
v_{\rho(n)}&\text{if $\rho=\sigma$,}\\0&\text{otherwise,}\end{cases}$$
we obtain
\begin{equation}
\sum_{i\in J}\lambda_i\;b^{(i)}_1\otimes b^{(i)}_2\otimes\cdots\otimes
b^{(i)}_k=0.
\label{eq:PfPrSIFAc}
\end{equation}
By assumption however, the sequences $(b^{(i)}_1,b^{(i)}_2,\ldots,
b^{(i)}_k)$ are distinct when $i$ runs over $J$. Therefore the elements
$b^{(i)}_1\otimes b^{(i)}_2\otimes\cdots\otimes b^{(i)}_k$ are linearly
independent in $W_{c_1}\otimes W_{c_2}\otimes\cdots\otimes W_{c_k}$,
for $B_{c_1}\otimes B_{c_2}\otimes \cdots\otimes B_{c_k}$ is a basis of
this module. Since $V$ and the $W_{c_i}$ are flat modules, the images
of the elements $b^{(i)}_1\otimes b^{(i)}_2\otimes\cdots\otimes b^{(i)}_k$
in $V^{\otimes n}$ are linearly independent. We then reach a contradiction
with Equation~(\ref{eq:PfPrSIFAc}).
\end{proof}

Before we look for a condition on $W$ that would ensures that $\Sigma(W)$
is a subbialgebra of $\mathscr F(V)$, we introduce a piece of notation
that will be needed later, especially in Section~\ref{ss:SplitForm}. Let
$\mathbf c=(c_1,c_2,\ldots,c_k)$ be a composition (possibly with parts
equal to zero)\footnote{It is convenient in this context to allow
compositions to have parts equal to zero. We could use a special
terminology, following for example Reutenauer who coined in
\cite{Reutenauer93} the word pseudocomposition for that purpose. To limit
the advent of new words, we will however simply say `composition
(possibly with parts equal to zero).'} of $n$. Since $V^{\otimes n}=V^{
\otimes c_1}\otimes V^{\otimes c_2}\otimes\cdots\otimes V^{\otimes c_k}$,
each tensor $t\in V^{\otimes n}$ can be written as a linear combination
of products $t_1\otimes t_2\otimes\cdots\otimes t_k$, where $t_i\in
V^{\otimes c_i}$ for each $i$. We denote such a decomposition by
$t=\sum_{(t)}t^{(\mathbf c)}_1\otimes t^{(\mathbf c)}_2\otimes\cdots\otimes
t^{(\mathbf c)}_k$. In this equation, the symbol $t^{(\mathbf c)}_i$ is
meant as a place-holder for the actual elements $t_i$. With this notation,
the coproduct of an element of the form $t\#e_n$ is
\begin{equation}
\Delta(t\#e_n)=\sum_{n'=0}^n\Bigl[t^{((n',n-n'))}_1\#e_{n'}\Bigr]
\otimes\Bigl[t^{((n',n-n'))}_2\#e_{n-n'}\Bigr].
\label{eq:CoprodGenDA}
\end{equation}

Let us now return to our study of the descent algebras. We introduce
the following condition on a graded submodule $W=\bigoplus_{n\geq0}W_n$
of $\T(V)$:
\begin{description}
\item[(A)]There holds $W_n\subseteq W_{c_1}\otimes W_{c_2}\otimes\cdots
\otimes W_{c_k}$ for any composition (possibly with parts equal to zero)
$\mathbf c=(c_1,c_2,\ldots,c_k)$ of a positive integer $n$.\footnote{We
abusively confuse $W_{c_1}\otimes W_{c_2}\otimes\cdots\otimes W_{c_k}$
with its image in $V^{\otimes c_1}\otimes V^{\otimes c_2}\otimes\cdots
\otimes V^{\otimes c_k}=V^{\otimes n}$. Of course no ambiguity arises
when $\mathbb K$ is a field or $V$ is torsion-free module over a p.i.d.}
\end{description}
In other words, for any composition $\mathbf c=(c_1,c_2,\ldots,c_k)$
of a positive integer $n$ and any $t\in W_n$, we may assume that in the
writing $t=\sum_{(t)}t^{(\mathbf c)}_1\otimes t^{(\mathbf c)}_2
\otimes\cdots\otimes t^{(\mathbf c)}_k$, all the elements of
$V^{\otimes c_i}$ represented by the place-holder $t^{(\mathbf c)}_i$
can be picked in $W_{c_i}$. We can now find a sufficient condition for
$\Sigma(W)$ to be a subbialgebra of $\mathscr F(V)$.

\begin{proposition}
\label{pr:SigmaWIsSubbialg}
If $W$ satisfies Condition~(A), then $\Sigma(W)$ is a graded
subbialgebra of $\mathscr F(V)$.
\end{proposition}
\begin{proof}
We have already seen that $\Sigma(W)$ is a graded subalgebra of
$\mathscr F(V)$. It remains to prove the inclusion
$$\Sigma(W)\subseteq\{x\in\mathscr F(V)\mid\Delta(x)\in\Sigma(W)
\otimes\Sigma(W)\}.$$

The set $E$ on the right of the symbol $\subseteq$ above is a
subalgebra of $\mathscr F(V)$, because $\Delta$ is a morphism of
algebras and $\Sigma(W)\otimes\Sigma(W)$ is a subalgebra. Moreover,
Equation~(\ref{eq:CoprodGenDA}) shows that if $W$ satisfies
Condition~(B), then $E$ contains all the elements $t\#e_n$ with
$t\in W_n$. Since these elements generate $\Sigma(W)$ as an algebra,
it follows that $E$ contains $\Sigma(W)$.
\end{proof}

Besides the trivial choice $W=\T(V)$, there are two main examples. The
first one occurs with $W=\TS(V)$, the space of all symmetric tensors on
$V$.\footnote{\label{fn:AssNTCondA}Condition~(A) holds for $W=\TS(V)$ as
soon as $V$ is projective or $\mathbb K$ is a field or a Dedekind ring.
We do not know if these restrictions can be lifted.} We call the
corresponding subbialgebra $\Sigma(W)$ the Novelli-Thibon bialgebra and
we denote it by $\NT(V)$. One may notice that the assignment
$V\rightsquigarrow\NT(V)$ is functorial.

The second interesting example concerns the case where $V$ is the
underlying space of a coalgebra. We first fix two rather standard
notations that are convenient for dealing with coalgebras; we will
use them not only in the presentation below, but also later in
Section~\ref{ss:SplitForm} with the comultiplicative structure of
$\mathscr F(V)$. Let $C$ be a coalgebra with its coassociative coproduct
$\delta$ and its counit $\varepsilon$. We define the iterated coproducts
$\delta_n:C\to C^{\otimes n}$ by setting $\delta_0=\varepsilon$,
$\delta_1=\id_C$, $\delta_2=\delta$, and
$$\delta_n=\bigl(\delta\otimes(\id_C)^{\otimes n-2}\bigr)\circ\bigl(
\delta\otimes(\id_C)^{\otimes n-3}\bigr)\circ\cdots\circ\delta$$
for all $n\geq3$. The Sweedler notation proposes to write the image of
an element $v\in C$ by $\delta_n$ as
$$\delta_n(v)=\sum_{(v)}v_{(1)}\otimes v_{(2)}\otimes\cdots\otimes
v_{(n)};$$
in this writing, the symbol $v_{(i)}$ is a place-holder for an
actual element of $C$ which varies from one term to the other.

Now we assume that the module $V$ on which the free quasisymmetric
algebra $\mathscr F(V)$ is constructed is endowed with a structure of
a coalgebra, with a coproduct $\delta$ and a counit~$\varepsilon$. In this
case, we may consider the image $W_n$ of the iterated coproduct
$\delta_n:V\to V^{\otimes n}$ and we may set $W=\bigoplus_{n\geq0}W_n$.
For any composition (possibly with parts equal to zero) $\mathbf c=(c_1,
c_2,\ldots,c_k)$ of $n$ and any element $v\in V$, the coassociativity of
$\delta$ implies
\begin{equation}
\delta_n(v)=\sum_{(v)}\ \underbrace{\delta_{c_1}(v_{(1)})}_{(\delta_n(v)
)^{(\mathbf c)}_1}\otimes\underbrace{\delta_{c_2}(v_{(2)})}_{(
\delta_n(v))^{(\mathbf c)}_2}\otimes\cdots\otimes\underbrace{
\delta_{c_k}(v_{(k)})}_{(\delta_n(v))^{(\mathbf c)}_k},
\label{eq:CondAAlgMR}
\end{equation}
which shows that Condition~(B) holds. Therefore $\Sigma(W)$ is a
subbialgebra of $\mathscr F(V)$. We call it the Mantaci-Reutenauer
bialgebra of the coalgebra $V$ and we denote it by $\MR(V)$. The
assignment $V\rightsquigarrow\MR(V)$ is a covariant functor from the
category of $\mathbb K$-coalgebras to the category of $\mathbb N$-graded
bialgebras over $\mathbb K$. As we will see in Section~\ref{ss:SplitForm},
this construction is mainly useful when $V$ is a projective $\mathbb
K$-module and the coproduct of $V$ is cocommutative; in this case,
$\MR(V)$ is a subbialgebra of $\NT(V)$.

For convenience, we introduce the following special notation for the
generators of the Mantaci-Reutenauer bialgebra $\MR(V)$: given any
positive integer $n$ and any element $v\in V$, we set $y_{n,v}=
\bigl[\delta_n(v)\#e_n\bigr]$. Equations~(\ref{eq:CoprodGenDA})
and (\ref{eq:CondAAlgMR}) entail that the coproduct of $y_{n,v}$
is given by
\begin{equation}
\Delta(y_{n,v})=\sum_{(v)}\sum_{n'=0}^ny_{n',v_{(1)}}\otimes
y_{n-n',v_{(2)}}.
\label{eq:CoprodY}
\end{equation}
Moreover, Proposition~\ref{pr:SigmaIsFreeAlg} implies that if $V$ is a
free $\mathbb K$-module, then the associative algebra $\MR(V)$ is freely
generated by the elements $y_{n,v}$, where $n\geq1$ and $v$ is chosen in
a basis of $V$.

\section{Duality}
\label{se:Duality}
The main result of this section says that the dual bialgebra $\mathscr
F(V)^\vee$ of the free quasisymmetric bialgebra on $V$ is the free
quasisymmetric bialgebra $\mathscr F(V^\vee)$ on the dual module $V^\vee$.
This result is neither deep nor difficult, but has many interesting
consequences, as we will see in Sections~\ref{se:SoloTheoWrProd} and
\ref{se:ColCombHopfAlg}. We begin by a general and easy discussion of
duality for $\mathbb K$-modules and $\mathbb K$-bialgebras.

\subsection{Perfect pairings}
\label{ss:PerfPair}
We define the duality functor $?^\vee$ as the contravariant endofunctor
$\Hom_{\mathbb K}(?,\mathbb K)$ of the category of $\mathbb K$-modules.
In particular, this functor maps a morphism $f:M\to N$ to its transpose
$f^\vee:N^\vee\to M^\vee$. Restricted to the full subcategory consisting
of finitely generated projective $\mathbb K$-modules, the duality
functor is an anti-equivalence of categories.

Given two $\mathbb K$-modules $M$ and $N$, there is a canonical
isomorphism $(M\oplus N)^\vee\cong M^\vee\oplus N^\vee$ and a canonical
map $N^\vee\otimes M^\vee\to(M\otimes N)^\vee$; the latter is an
isomorphism as soon as $M$ or $N$ is finitely generated and projective.
Given a $\mathbb K$-module $M$, there is a canonical homomorphism
$M\to M^{\vee\vee}$, which is an isomorphism if $M$ is finitely
generated and projective.

Let $H$ be a $\mathbb K$-bialgebra whose underlying space is finitely
generated and projective. Then the dual $H^\vee$ of $H$ is also a
bialgebra: the multiplication, the coproduct, the unit and the counit
of $H^\vee$ are the transpose of the coproduct, the multiplication,
the counit and the unit of $H$, respectively.

A pairing between two $\mathbb K$-modules $M$ and $N$ is a bilinear
form $\varpi:M\times N\to\mathbb K$. It gives rive to two linear maps
$\varpi^\flat:\bigl(M\to N^\vee,\ x\mapsto\varpi(x,?)\bigr)$ and
$\varpi^\#:\bigl(N\to M^\vee,\ y\mapsto\varpi(?,y)\bigr)$. The pairing
$\varpi$ is called perfect if the maps $\varpi^\flat$ and $\varpi^\#$
are isomorphisms. A pairing on a $\mathbb K$-module $M$ is a pairing
between $M$ and itself; such a pairing $\varpi$ is called symmetric
if $\varpi^\flat=\varpi^\#$.

In the case where the $\mathbb K$-modules $M$ and $N$ are finitely
generated and projective, we may identify $M$ and $N$ with their
respective biduals, and for any pairing $\varpi$ between $M$ and $N$,
it holds $\varpi^\#=(\varpi^\flat)^\vee$. If moreover $M$ and $N$ are
bialgebras, then $M^\vee$ and $N^\vee$ are also bialgebras; in this
situation, a pairing $\varpi$ between $M$ and $N$ such that
$\varpi^\flat$ and $\varpi^\#$ are morphisms of bialgebras is called
a pairing of bialgebras.

The above constructions concerning biduality or bialgebras are only
valid with finitely generated projective modules. We can however
relax the requirement of finite generation by working with $\mathbb
N$-graded modules. In this situation, we must adapt the definition
for the dual module: the dual of $M=\bigoplus_{n\geq0}M_n$ is the
graded module $M^\vee=\bigoplus_{n\geq0}(M_n)^\vee$, whose graded
components are the dual modules in the previous sense of the graded
components of $M$. We must also make the further assumptions that the
morphisms preserve the graduation and that pairings make graded
components of different degrees orthogonal to each other. Then
everything works as before, and biduality and duality of bialgebras
go smoothly as soon as the modules are projective with finitely
generated homogeneous components.

\subsection{Duality and the functor $\mathscr F$}
\label{ss:DualFunctF}
The following proposition examines the relationship between the functor
$\mathscr F$ and duality.

\begin{proposition}
\label{pr:FAndDuality}
There is a natural transformation from the contravariant functor
$\mathscr F(?^\vee)$ to the contravariant functor $\mathscr F(?)^\vee$,
which is an isomorphism when the domain of these functors is restricted to
the full subcategory of finitely generated projective $\mathbb K$-modules.
\end{proposition}
In other words, for any $\mathbb K$-module $V$, we can define a morphism
of graded algebras $c_V:\mathscr F(V^\vee)\stackrel\simeq\longrightarrow
\mathscr F(V)^\vee$, the construction being such that the assignment
$V\rightsquigarrow c_V$ is natural in $V$, and that $c_V$ is an
isomorphism of bialgebras if $V$ is finitely generated and projective.
\begin{proof}
Let $V$ be a $\mathbb K$-module. With the help of the canonical duality
bracket $\langle?,?\rangle:V\times V^\vee\to\mathbb K$ between $V$ and
$V^\vee$, we define for each $n\geq0$ a pairing $\langle?,?\rangle_n$
between $\mathscr F_n(V)$ and $\mathscr F_n(V^\vee)$ by the following
formula:
\begin{equation}
\left\langle\bigl[(v_1\otimes v_2\otimes\cdots\otimes v_n)\#\sigma\bigr],
\bigl[(f_1\otimes f_2\otimes\cdots\otimes f_n)\#\pi\bigr]\right\rangle_n=
\begin{cases}\prod_{i=1}^n\langle v_{\sigma(i)},f_i\rangle,&
\text{if $\sigma=\pi^{-1}$,}\\0&\text{otherwise,}\end{cases}
\label{eq:DefPairingFV}
\end{equation}
where $(v_1,v_2,\ldots,v_n)\in V^n$, $(f_1,f_2,\ldots,f_n)\in(V^\vee)^n$,
and $\sigma$ and $\pi$ are elements of $\mathfrak S_n$. If $V$ is
assumed to be finitely generated and projective, the canonical duality
between $V$ and $V^\vee$ is perfect and extends to a perfect pairing
between $V^{\otimes n}$ and $(V^\vee)^{\otimes n}$, which implies that
the pairing $\langle?,?\rangle_n$ is perfect.

We combine these pieces to define a pairing $\langle?,?\rangle_{\mathrm
{tot}}$ between $\mathscr F(V)$ and $\mathscr F(V^\vee)$ by setting
$$\langle\alpha,\xi\rangle_{\mathrm{tot}}=\sum_{n\geq0}\langle
\alpha_n,\xi_n\rangle_n$$
for all $\alpha=\sum_{n\geq0}\alpha_n$ and $\xi=\sum_{n\geq0}\xi_n$,
where $\alpha_n\in\mathscr F_n(V)$ and $\xi_n\in\mathscr F_n(V^\vee)$.
The map
$$c_V:\mathscr F(V^\vee)\to\mathscr F(V)^\vee,\ x\mapsto\langle?,x
\rangle_{\mathrm{tot}}$$
is a morphism of $\mathbb K$-modules; it is even an isomorphism if $V$
is finitely generated and projective.

A straightforward verification shows that the product $*$ and the
coproduct $\Delta$ of $\mathscr F(V)$ are adjoint to the coproduct
$\Delta$ and to the product $*$ of $\mathscr F(V^\vee)$ with respect
to the pairing $\langle?,?\rangle_{\mathrm{tot}}$. Together with a
similar statement about the unit and the counits, this implies that
$c_V$ is a morphism of algebras, and even of bialgebras if $\mathscr
F(V)$ is projective with finitely generated homogeneous components.
One checks also easily the commutativity of the diagram
$$\xymatrix@C+=50pt{\mathscr F(V^\vee)\ar[d]_{c_V}&\mathscr F(W^\vee)
\ar[d]^{c_W}\ar[l]_{\mathscr F(f^\vee)}\\\mathscr F(V)^\vee&\mathscr
F(W)^\vee\ar[l]_{\mathscr F(f)^\vee}}$$
for any $\mathbb K$-linear map $f:V\to W$ of $\mathbb K$-modules. This
means that the assignment $V\rightsquigarrow c_V$ is a natural
transformation from $\mathscr F(?^\vee)$ to $\mathscr F(?)^\vee$,
which completes the proof.
\end{proof}

Using the precise definition of the maps $c_V$ given in the proof of
Proposition~\ref{pr:FAndDuality}, one may check the following additional
property: the two compositions
$$\mathscr F(V)\longrightarrow\mathscr F(V^{\vee\vee})
\xrightarrow{c_{(V^\vee)}}\mathscr F(V^\vee)^\vee\quad\text{and}\quad
\mathscr F(V)\longrightarrow\mathscr F(V)^{\vee\vee}
\xrightarrow{(c_V)^\vee}\mathscr F(V^\vee)^\vee$$
are equal. Abusing the notations, we will write the above equality
as $c_{(V^\vee)}=(c_V)^\vee$.

Now suppose that $\varpi$ is a pairing between two $\mathbb K$-modules
$V$ and $W$. We can then define a pairing $\varpi_{\mathrm{tot}}$
between $\mathscr F(V)$ and $\mathscr F(W)$ by the equality
$\varpi_{\mathrm{tot}}{}^\flat=c_W\circ\mathscr F(\varpi^\flat)$; in
other words, we set
$$\varpi_{\mathrm{tot}}(x,y)=\bigl(c_W\circ\mathscr F(\varpi^\flat)
\bigr)(x)(y),$$
where $x\in\mathscr F(V)$ and $y\in\mathscr F(W)$. Then
$$\varpi_{\mathrm{tot}}{}^\#=(\varpi_{\mathrm{tot}}{}^\flat)^\vee
=\mathscr F(\varpi^\flat)^\vee\circ(c_W)^\vee=\mathscr
F(\varpi^\flat)^\vee\circ c_{(W^\vee)}=c_V\circ\mathscr
F\bigl((\varpi^\flat)^\vee\bigr)=c_V\circ\mathscr F(\varpi^\#).$$
The equalities $\varpi_{\mathrm{tot}}{}^\flat=c_W\circ\mathscr
F(\varpi^\flat)$ and $\varpi_{\mathrm{tot}}{}^\#=c_V\circ\mathscr
F(\varpi^\#)$ show that $\varpi_{\mathrm{tot}}$ is a pairing of
bialgebras. Moreover if $\varpi$ is perfect, then so is
$\varpi_{\mathrm{tot}}$. In the case $V=W$, one can also see that
the symmetry of $\varpi$ entails that of $\varpi_{\mathrm{tot}}$.

\subsection{Orthogonals and polars}
\label{ss:OrthoPolars}
Let $M$ be a finitely generated projective $\mathbb K$-module. We view
it as an `ambient' space and identify it with its bidual $M^{\vee\vee}$.
We define the orthogonal of a submodule $S$ of $M$ as the submodule
$S^\perp=\{f\in M^\vee\mid f\bigr|_S=0\}$ of $M^\vee$. Then $S^\perp$
is canonically isomorphic to $(M/S)^\vee$. Likewise, the orthogonal of
a submodule $T$ of $M^\vee$ is a submodule $T^\perp$ of $M$.

Let $\mathscr S$ be the set of all submodules $S$ of $M$ such that $M/S$
is projective, or in other words, that are direct summands of $M$.
If $S\in\mathscr S$, then both $S$ and $M/S$ are finitely generated
projective $\mathbb K$-modules. Likewise, let $\mathscr T$ be the set of
all submodules $T$ of $M^\vee$ that are direct summands of $M^\vee$. We
endow both $\mathscr S$ and $\mathscr T$ with the partial order given by
the inclusion of submodules. The following results are well-known in
this context:
\begin{itemize}
\item The maps $\bigl(\mathscr S\to\mathscr T,\ S\mapsto S^\perp\bigr)$
and $\bigl(\mathscr T\to\mathscr S,\ T\mapsto T^\perp\bigr)$ are
mutually inverse, order decreasing bijections.
\item For any $S\in\mathscr S$, there is a canonical isomorphism
$S^\vee\cong M^\vee/S^\perp$. Moreover for each submodule $S'\subseteq S$,
there is a canonical isomorphism $(S/S')^\vee\cong
S^{\prime\perp}/S^\perp$.
\item Let $S$ and $S'$ be two elements in $\mathscr S$. We always have
$(S+S')^\perp=S^\perp\cap S^{\prime\perp}$ and $S^\perp+S^{\prime\perp}
\subseteq(S\cap S')^\perp$. If moreover $S+S'$ belongs to $\mathscr S$,
then so does $S\cap S'$, and the equality $(S\cap S')^\perp=S^\perp
+S^{\prime\perp}$ holds.
\item Assume that $M$ is endowed with the structure of a bialgebra. Then
a submodule $S\in\mathscr S$ is a subbialgebra of $M$ if and only if
$S^\perp$ is a biideal of $M^\vee$, and a submodule $T\in\mathscr T$ is
a subbialgebra of $M^\vee$ if and only if $T^\perp$ is a biideal of $M$.
\end{itemize}
Given two submodules $S\in\mathscr S$ and $T\in\mathscr T$, we have
then sequences of canonical maps
\begin{align}
\begin{split}
T/(S^\perp\cap T)\cong(S^\perp+T)/S^\perp=(S^\perp+T^{\perp\perp})/
S^\perp\hookrightarrow(S\cap T^\perp)^\perp/S^\perp\cong\bigl(
S/(S\cap T^\perp)\bigr){}^\vee,\\[4pt]
S/(S\cap T^\perp)\cong(S+T^\perp)/T^\perp=(S^{\perp\perp}+T^\perp)/
T^\perp\hookrightarrow(S^\perp\cap T)^\perp/T^\perp\cong\bigl(
T/(S^\perp\cap T)\bigr){}^\vee.
\end{split}\label{eq:DualSubquot}
\end{align}
In other words, there is a canonical pairing between $S/(S\cap T^\perp)$
and $T/(S^\perp\cap T)$, which is perfect as soon as $(S+T^\perp)\in
\mathscr S$ and $(S^\perp+T)\in\mathscr T$.

We assume now that the module $M$ is endowed with a symmetric and
perfect pairing $\varpi$. Then to any submodule $S$ of $M$ we can
associate its polar $P^\circ=\bigl(\varpi^\flat\bigr){}^{-1}
\bigl(S^\perp\bigr)$ with respect to $\varpi$. Using $\varpi^\flat$,
one can deduce properties for polar submodules analogous to the
properties for orthogonals recalled above.

One can also adapt these results to the case where the projective module
$M$ is not finitely generated, provided it is graded with finitely
generated homogeneous components.

This material will prove useful in Sections~\ref{ss:SoloHomo} and
\ref{se:ColCombHopfAlg}, where we will meet instances of the following
situation. Here $V$ is a finitely generated projective $\mathbb K$-module,
endowed with a symmetric and perfect pairing $\varpi$. Then
$\mathscr F(V)$ is a projective $\mathbb K$-module, graded with
finitely generated homogeneous components, and endowed with the perfect
and symmetric pairing $\varpi_{\mathrm{tot}}$. Let moreover $S$ be a
graded subbialgebra of $\mathscr F(V)$, assumed to be a direct summand
of the graded $\mathbb K$-module $\mathscr F(V)$. We have then the
following commutative diagram of graded bialgebras,
\begin{equation}
\raisebox{35pt}{\xymatrix@!0@C=50pt@R=35pt{&\mathscr F(V)\ar[rr]^\simeq
\ar@{->>}[dr]&&\mathscr F(V)^\vee\ar@{->>}[dr]&\\S\ \ar@{^{(}->}[ur]
\ar@{->>}[dr]&&\mathscr F(V)/S^\circ\ar[rr]^\simeq&&S^\vee.\\&S/(S\cap
S^\circ)\;\ar@{^{(}->}[ur]\ar@{^{(}->}[rr]&&\bigl(S/(S\cap S^\circ)\bigr)
{}^\vee\ar@{^{(}->}[ur]&}}
\label{eq:AutodualDiagr}
\end{equation}
Here the horizontal arrows are induced by $\varpi_{\mathrm{tot}}{}^\flat$;
the one at the bottom line is the pairing on $S/(S\cap S^\circ)$ defined
by the sequences~(\ref{eq:DualSubquot}) with the choice
$T=\varpi_{\mathrm{tot}}{}^\flat(S)$.

To conclude this section, we show that the framework above is general
enough to accomodate the case of a Mantaci-Reutenauer bialgebra,
viewed as a submodule in a free quasisymmetric bialgebra.
\begin{proposition}
\label{pr:MRIsDirSum}
For any $\mathbb K$-coalgebra $V$, the submodule $\MR(V)$ is a direct
summand of $\mathscr F(V)$.
\end{proposition}
\begin{proof}
We will show two facts:
\begin{description}
\item[a)]The submodule $\Sigma(\T(V))$ has a graded complement in
$\mathscr F(V)$.
\item[b)]The Mantaci-Reutenauer bialgebra $\MR(V)$ has a graded
complement in $\Sigma(\T(V))$.
\end{description}

\noindent Let $n$ be a positive integer. The map
$$\mathbb K\mathfrak S_n\otimes V^{\otimes n}\to\mathscr F_n(V),\
\sigma\otimes t\mapsto\sigma\cdot(t\#e_n),$$
where $\sigma\in\mathfrak S_n$ and $t\in V^{\otimes n}$, is an
isomorphism of $\mathbb K$-modules. The submodule $\mathscr F_n(V)\cap
\Sigma(\T(V))$ is spanned by elements of the form
$$x_{\mathbf c}\cdot(t\#e_n)=\sum_{\substack{\sigma\in\mathfrak
S_n\\[2pt]D(\sigma)\preccurlyeq\mathbf c}}\sigma\cdot(t\#e_n),$$
where $\mathbf c$ is a composition of $n$. Therefore the submodule of
$\mathscr F_n(V)$ spanned over $\mathbb K$ by the elements $(\sigma-\pi)
\cdot(t\#e_n)$, where $t\in V^{\otimes n}$ and $\sigma$ and $\pi$ are
two permutations in $\mathfrak S_n$ with $D(\sigma)=D(\pi)$, is
complementary to $\mathscr F_n(V)\cap\Sigma(\T(V))$. This proves
Claim~a).

Let us now denote the coproduct and the counit of $V$ by $\delta$ and
$\varepsilon$, respectively. Let $n$ be a positive integer. We denote
the image of the iterated coproduct $\delta_n:V\to V^{\otimes n}$ by
$W_n$. The short exact sequence
$$0\to V\stackrel{\delta_n}\longrightarrow V^{\otimes n}\to V^{\otimes
n}/W_n\to0$$
splits, because the map $V^{\otimes n}\xrightarrow{\varepsilon^{\otimes
n-1}\otimes\id_V}\mathbb K^{\otimes n-1}\otimes V\cong V$ is a retraction
of $\delta_n$. Therefore we can find a complementary submodule $Z_n$
of $W_n$ in $V^{\otimes n}$. Given a composition $\mathbf c=(c_1,c_2,
\ldots,c_k)$ of $n$, we set
$$W_{\mathbf c}=W_{c_1}\otimes W_{c_2}\otimes\cdots\otimes W_{c_k}
\quad\text{and}\quad Z_{\mathbf c}=\sum_{i=1}^kV_{c_1}\otimes\cdots
V_{c_{i-1}}\otimes Z_{c_i}\otimes V_{c_{i+1}}\otimes\cdots\otimes
V_{c_k},$$
so that $V^{\otimes n}=W_{\mathbf c}\oplus Z_{\mathbf c}$. Then
$$\mathscr F_n(V)\cap\Sigma(\T(V))=\bigoplus_{\mathbf c\models
n}\bigl[x_{\mathbf c}\cdot(V^{\otimes n}\#e_n)\bigr]=\underbrace{
\bigoplus_{\mathbf c\models n}\bigl[x_{\mathbf c}\cdot(W_{\mathbf c}
\#e_n)\bigr]}_{\mathscr F_n(V)\cap\MR(V)}\oplus\bigoplus_{\mathbf
c\models n}\bigl[x_{\mathbf c}\cdot(Z_{\mathbf c}\#e_n)\bigr],$$
which shows Claim~b) and completes the proof.
\end{proof}

\section{The internal product}
\label{se:InternProd}
In this section, we consider the case where $V$ is the underlying space
of an algebra $A$. This affords a new structure on $\mathscr F(A)$, called
the internal product. We study ways to construct subalgebras of
$\mathscr F(A)$ for the internal product and clarify the situation
that arises when $A$ is a symmetric algebra, that is, an algebra
endowed with a symmetric, associative and perfect pairing.

\subsection{The twisted group ring $\mathscr F_n(A)$}
\label{ss:TwistedGpRing}
So let $A$ be a $\mathbb K$-algebra. Then the group $\mathfrak S_n$ acts
on the tensor power $A^{\otimes n}$ by automorphisms of algebra, which
allows to construct a twisted group ring, which we denote by $(A^{\otimes
n})\#(\mathbb K\mathfrak S_n)$. (In the language of Hopf algebras, one
says that $A^{\otimes n}$ is a $\mathbb K\mathfrak S_n$-module algebra,
and then the twisted group ring $(A^{\otimes n})\#(\mathbb K\mathfrak S_n)$
is viewed as a particular case of the smash product construction; see for
instance \cite{Montgomery93}.) This twisted group ring is our $\mathbb
K\mathfrak S_n$-bimodule $\mathscr F_n(A)$ endowed additionally with the
structure of an algebra. The associative product is given by the rule
\begin{multline}
\bigl[(a_1\otimes a_2\otimes\cdots\otimes a_n)\#\sigma\bigr]\cdot
\bigl[(b_1\otimes b_2\otimes\cdots\otimes b_n)\#\tau\bigr]\\
=\bigl[\bigl(a_1b_{\sigma^{-1}(1)}\otimes a_2b_{\sigma^{-1}(2)}\otimes
\cdots\otimes a_nb_{\sigma^{-1}(n)}\bigr)\#\bigl(\sigma\tau\bigr)\bigr]
\label{eq:ProdTwistedGpRing}
\end{multline}
and the unit is $1^{\otimes n}\#e_n$. The structure map $\bigl(\mathbb
K\to A^{\otimes n},\ \lambda\mapsto\lambda1^{\otimes n}\bigr)$ gives
rise to an embedding of the group algebra $\mathbb K\mathfrak S_n=
\mathscr F_n(\mathbb K)$ into $\mathscr F_n(A)$, which allows to
represent the two-sided action of $\mathbb K\mathfrak S_n$ on $\mathscr
F_n(A)$ with the help of the product of $\mathscr F_n(A)$.

It is convenient to extend this product to the whole $\mathscr F(A)$
by linearity: if $\alpha=\sum_{n\geq0}\alpha_n$ and $\alpha'=\sum_{n\geq0}
\alpha'_n$ with $\alpha_n$ and $\alpha'_n$ in $\mathscr F_n(A)$, we define
$\alpha\cdot\alpha'=\sum_{n\geq0}\alpha_n\cdot\alpha'_n$. This `internal
product' as it is called lacks a unit element.

More generally, given two $\mathbb K$-modules $V$ and $W$, the composition
$$(V^{\otimes n}\#\mathbb K\mathfrak S_n)\otimes(W^{\otimes n}\#\mathbb
K\mathfrak S_n)\twoheadrightarrow(V^{\otimes n}\#\mathbb K\mathfrak
S_n)\otimes_{\mathbb K\mathfrak S_n}(W^{\otimes n}\#\mathbb K\mathfrak
S_n)\stackrel\simeq\longrightarrow(V\otimes W)^{\otimes n}\#\mathbb K
\mathfrak S_n$$
defines a canonical morphism of $\mathbb K\mathfrak S_n$-modules from
$\mathscr F_n(V)\otimes\mathscr F_n(W)$ into $\mathscr F_n(V\otimes W)$.
Taking the direct sum over all $n\geq0$, one can define an `internal
product' $\mathscr F(V)\otimes\mathscr F(W)\to\mathscr F(V\otimes W)$
which is natural in $(V,W)$. Given a third $\mathbb K$-module $X$
and a linear map $m:V\otimes W\to X$, we obtain an internal product
$\mathscr F(V)\otimes\mathscr F(W)\to\mathscr F(X)$ by composition
with $\mathscr F(m)$. We will not pursue this way for want of
application, but it is worth noticing that even the apparently
simple case where $V$ or $W$ is the ground ring $\mathbb K$ is not
empty. We leave it to the reader to generalize the results of
Section~\ref{ss:SplitForm} to this wider context.

To conclude this section, we introduce two pieces of terminology that
will prove convenient in Section~\ref{ss:FrobStruct}. Let $S$ and $M$
be two graded submodules of $\mathscr F(A)$, and set $S_n=S\cap\mathscr
F_n(A)$ and $M_n=M\cap\mathscr F_n(A)$. We say that $S$ is a subalgebra
of $\mathscr F(A)$ for the internal product if each $S_n$ is a
subalgebra of $\mathscr F_n(A)$. In this case, we say further that $M$
is a left (respectively, right) internal $S$-submodule of $\mathscr
F(A)$ if $S\cdot M\subseteq M$ (respectively, $M\cdot S\subseteq M$).

\subsection{Double cosets in the symmetric group}
\label{ss:DblCosetSymGp}
In this section, we translate to the case of the symmetric group a
theorem of Solomon valid in the more general context of Coxeter
groups. The result will prove crucial in the proof of the splitting
formula in Section~\ref{ss:SplitForm}.

To begin with, let $\bigl(W,(s_i)_{i\in I}\bigr)$ be a Coxeter system.
Given a subset $J\subseteq I$, the parabolic subgroup $W_J$ is the
subgroup of $W$ generated by the elements $s_j$ with $j\in J$. In each
left coset $wW_J$, there is a unique element with minimal length, called
the distinguished representative of that coset. We denote by $X_J$ the
set of distinguished representatives of the left cosets modulo $W_J$.
Given a second subset $K\subseteq I$, there is likewise a unique element
with minimal length in each double coset $W_JwW_K$, unsurprisingly
called the distinguished representative of the double coset. The set of
distinguished representatives of the double cosets modulo $W_J$ and
$W_K$ is $(X_J)^{-1}\cap X_K$. The following statement is a rephrasing of
Theorem~2 of \cite{Solomon76}.

\begin{theorem}
\label{th:ThmSolomon}
Given a double coset $C\in W_J\backslash W/W_K$, we set
$$L_C=\{k\in K\mid\exists j\in J,\ x^{-1}s_jx=s_k\},$$
where $x\in C\cap(X_J^{-1}\cap X_K)$ is the distinguished representative
of $C$. Then $X_{L_C}$ is the disjoint union of the sets $X_Jw$, where
$w\in C\cap X_K$.
\end{theorem}

We now translate this proposition in a combinatorial language more
adapted to the case of the symmetric group $\mathfrak S_n$. Let
$\mathbf c=(c_1,c_2,\ldots,c_k)$ and $\mathbf d=(d_1,d_2,\ldots,d_l)$
be two compositions of $n$, and set $t_i=c_1+c_2+\cdots+c_i$ and
$u_j=d_1+d_2+\cdots+d_j$. We denote by $\mathscr M_{\mathbf c,\mathbf d}$
the set of all matrices $M=(m_{ij})$ with non-negative integral entries
in $k$ rows and $l$ columns and with row-sum $\mathbf c$ and column-sum
$\mathbf d$, that is,
$$c_i=\sum_{j=1}^lm_{ij}\quad\text{for all $i$ and}\quad
d_j=\sum_{i=1}^km_{ij}\quad\text{for all $j$.}$$
There is a well-known bijection from $\mathscr M_{\mathbf c,\mathbf d}$
onto the double quotient $\mathfrak S_{\mathbf c}\backslash\mathfrak
S_n/\mathfrak S_{\mathbf d}$ that maps a matrix $M=(m_{ij})$ to the
double coset
$$C(M)=\bigl\{\sigma\in\mathfrak S_n\bigm|\forall(i,j),\ m_{ij}=
\bigl|\,[t_{i-1}+1,t_i]\cap\sigma([u_{j-1}+1,u_j])\,\bigr|\,\bigr\}.$$
Finally, we associate to a matrix $M\in\mathscr M_{\mathbf c,\mathbf
d}$ its column-reading composition
$$\colr(M)=(m_{11},m_{21},\ldots,m_{k1},m_{12},m_{22},\ldots,m_{k2},
\ldots,m_{1l},m_{2l},\ldots,m_{kl}).$$
With these notations, Theorem~\ref{th:ThmSolomon} translates to the
following statement.

\begin{corollary}
\label{co:CombiSymmGp}
For any matrix $M\in\mathscr M_{\mathbf c,\mathbf d}$, the set
$X_{\colr(M)}$ is the disjoint union of the sets $X_{\mathbf c}\sigma$,
where $\sigma\in C(M)\cap X_{\mathbf d}$.
\end{corollary}
\begin{proof}
We set $I=\{1,2,\ldots,n-1\}$. For each $i\in I$, we call $s_i$ be
the transposition in $\mathfrak S_n$ that exchanges $i$ and $i+1$.
Endowed with the family $(s_i)_{i\in I}$, the group $\mathfrak S_n$
becomes a Coxeter system $W$.

Set $J=I\setminus\{t_1,t_2,\cdots,t_{k-1}\}$ and $K=I\setminus\{u_1,
u_2,\ldots,u_{l-1}\}$. Then the Young subgroups $\mathfrak S_{\mathbf c}$
and $\mathfrak S_{\mathbf d}$ coincide with the parabolic subgroups
$W_J$ and $W_K$, respectively; moreover the sets $X_{\mathbf c}$ and
$X_{\mathbf d}$ are the sets of distinguished representatives $X_J$
and $X_K$.

We now fix a matrix $M\in\mathscr M_{\mathbf c,\mathbf d}$. We define
a permutation $\rho\in\mathfrak S_n$ by the following rule: for each
$a\in[1,n]$, we determine the index $j\in[1,l]$ such that
$a\in[u_{j-1}+1,u_j]$ and then the index $i\in[1,k]$ such that
$a-u_{j-1}\in[m_{1j}+m_{2j}+\cdots+m_{i-1,j}+1,m_{1j}+m_{2j}+
\cdots+m_{ij}]$, and we set
$$\rho(a)=a-(u_{j-1}+m_{1j}+m_{2j}+\cdots+m_{i-1,j})+(t_{i-1}+
m_{i1}+m_{i2}+\cdots+m_{i,j-1}).$$
One checks without difficulty that $\rho\in C(M)\cap(X_{\mathbf
c})^{-1}\cap X_{\mathbf d}$, which implies that $\rho$ is the
distinguished representative of the double coset $C(M)\in\mathfrak
S_{\mathbf c}\backslash\mathfrak S_n/\mathfrak S_{\mathbf d}$.

Moreover, let $a\in[1,n]$ and determine the indices $i$ and $j$ as
above. One checks easily that
$$a-u_{j-1}\in[m_{1j}+m_{2j}+\cdots+m_{i-1,j}+1,m_{1j}+m_{2j}+\cdots+
m_{ij}-1]\Longleftrightarrow\left\{\begin{aligned}&a\in K,\\&\rho(a)
\in J,\\&\rho(a+1)=\rho(a)+1.\end{aligned}\right.$$
(The ket point here is to observe that if $a=u_{j-1}+m_{1j}+m_{2j}+
\cdots+m_{ij}$ and $a\in K$, then the inequalities $\rho(a+1)>t_i\geq
\rho(a)$ hold.)

For any $j\in J$, the permutation $\rho^{-1}s_j\rho$ is the transposition
that exchanges $\rho^{-1}(j)$ and $\rho^{-1}(j+1)$, with necessarily
$\rho^{-1}(j)<\rho^{-1}(j+1)$ because $\rho^{-1}\in X_J$. The definition
$$L_{C(M)}=\{k\in K\mid\exists j\in J,\ \rho^{-1}s_j\rho=s_k\}$$
translates therefore to the equality
$$L_{C(M)}=\bigcup_{j=1}^l\bigcup_{i=1}^k[u_{j-1}+m_{1j}+m_{2j}+
\cdots+m_{i-1,j}+1,u_{j-1}+m_{1j}+m_{2j}+\cdots+m_{ij}-1],$$
or, in other words, to
$$L_{C(M)}=I\setminus\{f_1,f_1+f_2,\ldots,f_1+f_2+\cdots+f_{m-1}\}$$
if the parts of $\colr(M)$ form the sequence $(f_1,f_2,\cdots,f_m)$.
This implies that the sets $X_{L(C)}$ and $X_{\colr(M)}$ coincide.

This completes the dictionary that allows to deduce the corollary
from Theorem~\ref{th:ThmSolomon}.
\end{proof}

\subsection{The splitting formula}
\label{ss:SplitForm}
The splitting formula, due to Gelfand, Krob, Lascoux, Leclerc, Retakh
and Thibon \cite{GKLLRT95} in the case of $\mathscr F(\mathbb K)$ and
to Novelli and Thibon \cite{Novelli-Thibon04} in the general case, is
the tool that enables to show that certain graded subbialgebras
$\Sigma(W)$ of Section~\ref{ss:SubbialgSigmaW} are subalgebras
of $\mathscr F(A)$ for the internal product. We begin with a lemma.

\begin{lemma}
\label{le:SplitLemma}
Let $V$ be a projective $\mathbb K$-module,\footnote{The assuption that
$V$ is projective guarantee the existence of decompositions
$a_i=\sum_{(a_i)}a^{(M)}_{i1}\otimes a^{(M)}_{i2}\otimes\cdots\otimes
a^{(M)}_{il}$ below, as mentioned in the footnote~\ref{fn:AssNTCondA}.}
let $n$ be a positive integer, let $\mathbf c=(c_1,c_2,\ldots,c_k)$ and
$\mathbf d$ be two compositions of $n$, and for each $i\in[1,k]$, let
$a_i\in\TS^{c_i}(V)$ be a symmetric tensor of degree $c_i$. The $i$-th
line of a matrix $M=(m_{ij})$ in $\mathscr M_{\mathbf c,\mathbf d}$ can
be seen as a composition (possibly with parts equal to zero) of $c_i$.
According to the decomposition
$$\TS^{c_i}(V)\subseteq\TS^{m_{i1}}(V)\otimes\TS^{m_{i2}}(V)\otimes
\cdots\otimes\TS^{m_{il}}(V),$$
we write $a_i$ as a linear combination $\sum_{(a_i)}a^{(M)}_{i1}\otimes
a^{(M)}_{i2}\otimes\cdots\otimes a^{(M)}_{il}$ with $a^{(M)}_{ij}\in
\TS^{m_{ij}}(V)$. Then in the $\mathbb K\mathfrak S_n$-bimodule
$\mathscr F_n(V)$, there holds
\begin{multline}
x_{\mathbf c}\cdot\bigl[(a_1\otimes a_2\otimes\cdots\otimes a_k)\#
x_{\mathbf d}\bigr]\\[2pt]
\begin{aligned}[b]
\smash{=\sum_{M\in\mathscr M_{\mathbf c,\mathbf d}}\
\sum_{(a_1),\,(a_2),\,\ldots,\,(a_k)}}\
x_{\colr(M)}\cdot\Bigl[\Bigl(a^{(M)}_{11}\otimes
a^{(M)}_{21}&\otimes\cdots\otimes a^{(M)}_{k1}\\
\otimes a^{(M)}_{12}&\otimes a^{(M)}_{22}\otimes\cdots\otimes
a^{(M)}_{k2}\otimes\cdots\\
&\otimes a^{(M)}_{1l}\otimes a^{(M)}_{2l}\otimes\cdots\otimes
a^{(M)}_{kl}\Bigr)\#e_n\Bigr].
\end{aligned}
\label{eq:SplitLemma}
\end{multline}
\end{lemma}
\begin{proof}
We set $t_i=c_1+c_2+\cdots+c_i$ and $u_j=d_1+d_2+\cdots+u_j$. We take
$M\in\mathscr M_{\mathbf c,\mathbf d}$ and $\sigma\in\mathfrak S_n$.
If $\sigma$ belongs to the double coset $C(M)$, then for each $j$, the
set $\sigma([u_{j-1}+1,u_j])$ has $m_{1j}$ elements in $[1,t_1]$, $m_{2j}$
elements in $[t_1+1,t_2]$, \dots, $m_{kj}$ elements in $[t_{k-1}+1,t_k]$.
On the other hand, if $\sigma$ belongs to $X_{\mathbf d}$, then it is an
increasing map on the interval $[u_{j-1}+1,u_j]$. Therefore, if $\sigma$
belongs to $C(M)\cap X_{\mathbf d}$, the $m_{ij}$ elements of the set
$$\sigma([u_{j-1}+m_{1j}+m_{2j}+\cdots+m_{i-1,j}+1,u_{j-1}+m_{1j}+
m_{2j}+\cdots+m_{ij}])$$
belong to $[t_{i-1}+1,t_i]$, so that
\begin{equation*}
(a_1\otimes a_2\otimes\cdots\otimes a_k)\#\sigma\
\begin{aligned}[t]
\smash{=\ \sum_{(a_1),\,(a_2),\,\ldots,\,(a_k)}}\
\sigma\cdot\Bigl[\Bigl(a^{(M)}_{11}\otimes
a^{(M)}_{21}&\otimes\cdots\otimes a^{(M)}_{k1}\\
\otimes a^{(M)}_{12}&\otimes a^{(M)}_{22}\otimes\cdots\otimes
a^{(M)}_{k2}\otimes\cdots\\
&\otimes a^{(M)}_{1l}\otimes a^{(M)}_{2l}\otimes\cdots\otimes
a^{(M)}_{kl}\Bigr)\#e_n\Bigr],
\end{aligned}
\end{equation*}
because each $a_i$ is symmetric. Using the notations of
Section~\ref{ss:DblCosetSymGp}, we decompose $X_{\mathbf d}$ as the
disjoint union $\coprod_{M\in\mathscr M_{\mathbf c,\mathbf d}}\bigl(
C(M)\cap X_{\mathbf d}\bigr)$. Then
\begin{align*}
x_{\mathbf c}\cdot\bigl[(a_1\otimes a_2\otimes\cdots\otimes&a_k)\#
x_{\mathbf d}\bigr]\\
&=\sum_{M\in\mathscr M_{\mathbf c,\mathbf d}}\ \sum_{\sigma\in C(M)\cap
X_{\mathbf d}}x_{\mathbf c}\cdot\bigl[(a_1\otimes a_2\otimes\cdots\otimes
a_k)\#\sigma\bigr]\\
&=\sum_{M\in\mathscr M_{\mathbf c,\mathbf d}}\ \sum_{\sigma\in C(M)\cap
X_{\mathbf d}}\sum_{(a_i)}(x_{\mathbf c}\sigma)\cdot\Bigl[\Bigl(
a^{(M)}_{11}\otimes a^{(M)}_{21}\otimes\cdots\otimes a^{(M)}_{kl}\Bigr)
\#e_n\Bigr]\\
&=\sum_{M\in\mathscr M_{\mathbf c,\mathbf d}}\sum_{(a_i)}
x_{\colr(M)}\cdot\Bigl[\Bigl(a^{(M)}_{11}\otimes a^{(M)}_{21}
\otimes\cdots\otimes a^{(M)}_{kl}\Bigr)\#e_n\Bigr],
\end{align*}
the last equality coming from Corollary~\ref{co:CombiSymmGp}. This
calculation proves Lemma~\ref{le:SplitLemma}.
\end{proof}

In the remainder of this section, the letter $A$ denotes a $\mathbb
K$-algebra, whose underlying module is projective. We now state and
prove the splitting formula.
\begin{theorem}
\label{th:SplitForm}
Let $y$ be an element in $\NT(A)$ and $z_1$, $z_2$, \ldots, $z_l$ be
elements in $\mathscr F(A)$. Then
\begin{equation}
y\cdot(z_1*z_2*\cdots*z_l)=\sum_{(y)}(y_{(1)}\cdot z_1)*
(y_{(2)}\cdot z_2)*\dots*(y_{(l)}\cdot z_l).
\label{eq:SplitForm}
\end{equation}
\end{theorem}
\begin{proof}
By linearity, it is sufficient to prove Formula~(\ref{eq:SplitForm})
for elements $y$ of the form
$$y=(a_1\#e_{c_1})*(a_2\#e_{c_2})*\cdots*(a_k\#e_{c_k})=x_{\mathbf c}
\cdot\bigl[(a_1\otimes a_2\otimes\cdots\otimes a_k)\#e_n\bigr],$$
where $\mathbf c=(c_1,c_2,\ldots,c_k)$ is a composition of a positive
integer $n$ and where $a_1$, $a_2$, \dots, $a_k$ are symmetric tensors
on $A$ of degree $c_1$, $c_2$, \dots, $c_k$, respectively. By Formula
(\ref{eq:CoprodGenDA}), the $l$-th iterated coproduct of the element
$(a_i\#e_i)$ is
$$\Delta_l(a_i\#e_{c_i})=\sum_{\mathbf f}\ \sum_{(a_i)}\
\bigl((a^{}_i)^{(\mathbf f)}_1\#e_{f_1}\bigr)\otimes\bigl(
(a_i)^{(\mathbf f)}_2\#e_{f_2}\bigr)\otimes\cdots\otimes\bigl(
(a_i)^{(\mathbf f)}_l\#e_{f_l}\bigr),$$
where the first sum runs over all compositions $\mathbf f=(f_1,f_2,
\ldots,f_l)$ of $c_i$ in $l$ parts (possibly equal to zero).
Multiplying these expressions for $i=1,2,\ldots,k$ and expanding, we
obtain
\begin{align}
\Delta_l(y)=\sum_{\mathbf g}\;\sum_{M\in\mathscr M_{\mathbf c,\mathbf g}}
\quad\sum_{(a_1),\,(a_2),\,\ldots,\,(a_k)}\quad
&\Bigl[\Bigl(a^{(M)}_{11}\#e_{m_{11}}\Bigr)*\Bigl(a^{(M)}_{21}\#e_{m_{21}}
\Bigr)*\cdots*\Bigl(a^{(M)}_{k1}\#e_{m_{k1}}\Bigr)\Bigr]\notag\\
\otimes&\Bigl[\Bigl(a^{(M)}_{12}\#e_{m_{12}}\Bigr)*\Bigl(a^{(M)}_{22}
\#e_{m_{22}}\Bigr)*\cdots*\Bigl(a^{(M)}_{k2}\#e_{m_{k2}}\Bigr)\Bigr]
\notag\\[4pt]
\otimes&\cdots\notag\\[2pt]
\otimes&\Bigl[\Bigl(a^{(M)}_{1l}\#e_{m_{1l}}\,\Bigr)*\Bigl(a^{(M)}_{2l}
\#e_{m_{2l}}\,\Bigr)*\cdots*\Bigl(a^{(M)}_{kl}\#e_{m_{kl}}\Bigr)\Bigr],
\label{eq:PfSplitForma}
\end{align}
where the first sum runs over all compositions (possibly with zero
parts) $\mathbf g$ of $n$ in $l$ parts and where for each matrix
$M\in\mathscr M_{\mathbf c,\mathbf g}$, the tensors $a_i$ are decomposed
as in the statement of Lemma~\ref{le:SplitLemma}.

Now let $\mathbf d=(d_1,d_2,\ldots,d_l)$ be a composition (possibly
with parts equal to zero) of $n$ and consider the equality proved in
Lemma~\ref{le:SplitLemma}. The left-hand side of (\ref{eq:SplitLemma})
is equal to
\begin{equation}
x_{\mathbf c}\cdot\bigl[(a_1\otimes a_2\otimes\cdots\otimes a_k)\#
e_n\bigr]\cdot(1^{\otimes n}\#x_{\mathbf d})=y\cdot(1^{\otimes n}\#
x_{\mathbf d}).
\label{eq:PfSplitFormb}
\end{equation}
On the other hand, Equation~(\ref{eq:PfSplitForma}) joint to
Formula~(\ref{eq:PropicXc}) shows that the right-hand side of
(\ref{eq:SplitLemma}) is equal to
\begin{equation}
\sum_{(y)}\bigl(y_{(1)}\cdot(1^{\otimes d_1}\#e_{d_1})\bigr)*
\bigl(y_{(2)}\cdot(1^{\otimes d_2}\#e_{d_2})\bigr)*\dots*
\bigl(y_{(l)}\cdot(1^{\otimes d_l}\#e_{d_l})\bigr).
\label{eq:PfSplitFormc}
\end{equation}
We conclude that the quantities (\ref{eq:PfSplitFormb}) and
(\ref{eq:PfSplitFormc}) are equal.

By linearity, we may assume that the elements $z_j$ are of the form
$z_j=(b_j\#\sigma_j)$, where $b_j\in A^{\otimes d_j}$ and
$\sigma_j\in\mathfrak S_{d_j}$. For degree reasons, both sides of
(\ref{eq:SplitForm}) vanish unless $n=d_1+d_2+\cdots+d_l$. We may
therefore assume without loss of generality that $\mathbf
d=(d_1,d_2,\ldots,d_l)$ is a composition of $n$. We now multiply
both (\ref{eq:PfSplitFormb}) and (\ref{eq:PfSplitFormc}) on the
right by $\bigl[(b_1\otimes b_2\otimes\cdots\otimes b_l)\#
(\sigma_1\times\sigma_2\times\cdots\times\sigma_l)\bigr]$, using
the internal product. These multiplications yield the left-hand and
the right-hand side of (\ref{eq:SplitForm}), respectively. The
theorem follows.
\end{proof}

As a first application of this formula, we consider the two following
conditions for a graded submodule $B=\bigoplus_{n\geq0}B_n$ of $\T(A)$.
\begin{description}
\item[(B)]Each $B_n$ is a subalgebra of $A^{\otimes n}$.
\item[(C)]Each space $B_n$ consists of symmetric tensors, that is,
$B\subseteq\TS(A)$.
\end{description}

\begin{corollary}
\label{co:SubalgIntSigmaB}
For any graded submodule $B$ of $\T(A)$ satisfying Conditions~(A),
(B) and (C), the descent bialgebra $\Sigma(B)$ is a subalgebra
of $\mathscr F(A)$ for the internal product.
\end{corollary}
\begin{proof}
We have to prove that for any elements $y$ and $z$ in $\Sigma(B)$,
the product $y\cdot z$ belongs to $\Sigma(B)$. We first consider
the case where $z$ is of the form $(b\#e_n)$, where $b\in B_n$. The
homogeneous components of $y$ whose degree are different from $n$ do not
contribute to the product $y\cdot z$; they can therefore be put aside.
We then write $y$ as a linear combination of products $(a_1\#e_{c_1})*
(a_2\#e_{c_2})*\cdots*(a_k\#e_{c_k})$, where $\mathbf c=(c_1,c_2,\ldots,
c_k)$ is a composition of $n$ and $a_1\in B_{c_1}$, $a_2\in B_{c_2}$,
\dots, $a_k\in B_{c_k}$. By Condition~(A), we may find a decomposition
$b=\sum_{(b)}b^{(\mathbf c)}_1\otimes b^{(\mathbf c)}_2\otimes\cdots
\otimes b^{(\mathbf c)}_k$ for each composition $\mathbf c$ that arises
in the expression of $y$, where the elements represented by the
place-holder $b^{(\mathbf c)}_i$ belong to $B_{c_i}$. Therefore
$y\cdot(b\#e_n)$ is a linear combination of elements of the form
\begin{align*}
\bigl(x_{\mathbf c}\cdot\bigl[(a_1\otimes a_2\otimes\cdots\otimes a_k)
\#e_n\bigr]\bigl)\cdot\bigl[(b_1\otimes b_2&\otimes\cdots\otimes b_k)
\#e_n\bigr]\\[3pt]&=x_{\mathbf c}\cdot\bigl[((a_1b_1)\otimes(a_2b_2)
\otimes\cdots\otimes(a_kb_k))\#e_n\bigr]\\[3pt]&=(a_1b_1\#e_{c_1})*
(a_2b_2\#e_{c_2})*\cdots*(a_kb_k\#e_{c_k}).
\end{align*}
Since each element $a_ib_i$ appearing here belongs to $B_{c_i}$ by
Condition~(B), $y\cdot(b\#e_n)$ is in $\Sigma(B)$.

In the general case, we may write $z$ as a linear combination of products
$z_1*z_2*\cdots*z_l$, where each $z_j$ is of the form $b_j\#e_{d_j}$, where
$d_j$ is a positive integer and $b_j\in B_{d_j}$. We apply the splitting
formula~(\ref{eq:SplitForm}). Since $\Sigma(B)$ is a subcoalgebra of
$\mathscr F(A)$ (Proposition~\ref{pr:SigmaWIsSubbialg}), we may require
that in the decomposition $\Delta_l(y)=\sum_{(y)}y_{(1)}\otimes
y_{(2)}\otimes\cdots\otimes y_{(l)}$ used, all elements represented by
the placeholders $y_{(j)}$ belong to $\Sigma(B)$. By the first case,
each product $y_{(j)}\cdot z_j$ belongs to $\Sigma(B)$, which entails
that $y\cdot(z_1*z_2*\cdots*z_l)$ belongs to $\Sigma(B)$. We conclude
that $y\cdot z$ belongs to $\Sigma(B)$.
\end{proof}

Again there are two main examples to which Corollary
\ref{co:SubalgIntSigmaB} can be applied. The first one is the case of
the Novelli-Thibon algebra: the sequence $B_n=\TS^n(A)$ satisfies
Conditions (A), (B) and (C), so the submodule $\NT(A)$ is a subalgebra
of $\mathscr F(A)$ for the internal product.

The second example arises when $A$ is a cocommutative bialgebra. Each
iterated coproduct $\delta_n:A\to A^{\otimes n}$ is a morphism of
algebras, therefore its image is a subalgebra $B_n$ of $A^{\otimes n}$.
The submodule $B=\bigoplus_{n\geq0}B_n$ therefore satisfies Condition
(B). It also satisfies Conditions~(A) and (C), because the coproduct of
$A$ is coassociative and cocommutative. It thus follows from
Corollary~\ref{co:SubalgIntSigmaB} that the submodule $\Sigma(B)$,
which is of course the Mantaci-Reutenauer bialgebra $\MR(A)$, is a
subalgebra of $\mathscr F(A)$ for the internal product.

The following corollary gives the rule to compute internal products in
a Mantaci-Reute\-nauer algebra. It generalizes Corollary~6.8 and
Theorem~6.9 of~\cite{Mantaci-Reutenauer95}.

\begin{corollary}
\label{co:MantReutRule}
Let $A$ be a cocommutative bialgebra with coproduct $\delta$, let $n$ be
a positive integer, let $\mathbf c=(c_1,c_2,\ldots,c_k)$ and $\mathbf
d=(d_1,d_2,\ldots d_l)$ be two compositions of $n$, and let $a_1$, $a_2$,
\dots, $a_k$, $b_1$, $b_2$, \dots, $b_l$ be elements of $A$. Then
$$(y_{c_1,a_1}*y_{c_2,a_2}*\cdots*y_{c_k,a_k})\cdot
(y_{d_1,b_1}*y_{d_2,b_2}*\cdots*y_{d_l,b_l})=\sum_{(a_i),(b_j)}\
\sum_{M\in\mathscr M_{\mathbf c,\mathbf d}}\left(\prod_{j=1}^l\;
\prod_{i=1}^k\;y_{m_{ij},a_{i(j)}b_{j(i)}}\right),$$
where the first summation symbol on the right comes from the Sweedler
notation for writing the iterated coproducts $\delta_l(a_i)$ and
$\delta_k(b_j)$, where the two successive symbols $\prod$ stand for the
external product, and where the factors of this external product are
formed by reading column by column the entries of the matrix $M=(m_{ij})$.
\end{corollary}
\begin{proof}
An easy induction based on Formula~(\ref{eq:CoprodY}) implies that
$$\Delta_l(y_{c_i,a_i})=\sum_{(a_i)}\ \sum_{\mathbf f}\
y_{f_1,a_{i(1)}}\otimes y_{f_2,a_{i(2)}}\otimes\cdots\otimes
y_{f_l,a_{i(l)}},$$
where the second sum runs over all compositions (possibly with zero
parts) $\mathbf f=(f_1,f_2,\ldots,f_l)$ of $c_i$ in $l$ parts.
Setting $y=y_{c_1,a_1}*y_{c_2,a_2}*\cdots*y_{c_k,a_k}$, it follows that
\begin{align}
\Delta_l(y)&=\Delta_l(y_{c_1,a_1})*\Delta_l(y_{c_2,a_2})*\cdots*
\Delta_l(y_{c_k,a_k})\notag\\
&=\sum_{(a_1),\,(a_2),\,\ldots,\,(a_k)}\;\sum_M\;\left(\prod_{i=1}^k
y_{m_{i1},a_{i(1)}}\right)\otimes\left(\prod_{i=1}^ky_{m_{i2},a_{i(2)}}
\right)\otimes\cdots\otimes\left(\prod_{i=1}^ky_{m_{il},a_{i(l)}}\right),
\label{eq:CoMRRa}
\end{align}
where the second sum is over all matrices $M$ with non-negative integral
entries in $k$ rows and $l$ columns and with row-sum $\mathbf c$.

We now use the splitting formula
\begin{equation}
y\cdot(y_{d_1,b_1}*y_{d_2,b_2}*\cdots*y_{d_l,b_l})=\sum_{(y)}
(y_{(1)}\cdot y_{d_1,b_1})*(y_{(2)}\cdot y_{d_2,b_2})*\cdots*
(y_{(l)}\cdot y_{d_l,b_l})
\label{eq:CoMRRb}
\end{equation}
and substitute in it the expression for the iterated coproduct
$\Delta_l(y)$ found in (\ref{eq:CoMRRa}). For degree reasons, each
term in (\ref{eq:CoMRRa}) that yields a non-zero contribution to the
right-hand side of (\ref{eq:CoMRRb}) corresponds to a matrix $M$ whose
column sum is equal to $\mathbf d$, so that we may restrict the sum
to the matrices $M$ in $\mathscr M_{\mathbf c,\mathbf d}$. The result
of the substitution is a sum of products; in each product, the $j$-th
factor is
\begin{align*}
y_{(j)}\cdot&\;y_{d_j,b_j}\\[2pt]
&=\bigl[(y_{m_{1j},a_{1(j)}})*(y_{m_{2j},a_{2(j)}})*\cdots*(y_{m_{kj},
a_{k(j)}})\bigr]\cdot y_{d_j,b_j}\\[2pt]
&=x_{(m_{1j},m_{2j},\ldots,m_{kj})}\cdot\bigl[\bigl(\delta_{m_{1j}}
(a_{1(j)})\otimes\delta_{m_{2j}}(a_{2(j)})\otimes\cdots\otimes
\delta_{m_{kj}}(a_{k(j)})\bigr)\#e_{d_j}\bigr]\cdot\bigl[\delta_{d_j}
(b_j)\#e_{d_j}\bigr]\\[2pt]
&=\sum_{(b_j)}x_{(m_{1j},m_{2j},\ldots,m_{kj})}\cdot\bigl[\bigl(
\delta_{m_{1j}}(a_{1(j)}b_{j(1)})\otimes\delta_{m_{2j}}(a_{2(j)}b_{j(2)})
\otimes\cdots\otimes\delta_{m_{kj}}(a_{k(j)}b_{j(k)})\bigr)\#e_{d_j}
\bigr]\\[2pt]
&=\sum_{(b_j)}y_{m_{1j},a_{1(j)}b_{j(1)}}*y_{m_{2j},a_{2(j)}b_{j(2)}}*
\cdots*y_{m_{kj},a_{k(j)}b_{j(k)}}.
\end{align*}
The corollary follows immediately.
\end{proof}

\subsection{Frobenius structures}
\label{ss:FrobStruct}
In this section, we put together the structures defined in
Sections~\ref{ss:DualFunctF} and \ref{ss:TwistedGpRing}.

We begin by recalling some terminology. Let $A$ be an associative
$\mathbb K$-algebra with unit. A pairing $\varpi$ on $A$ is said
associative if $\varpi(ab,c)=\varpi(a,bc)$ for all $(a,b,c)\in A^3$.
A trace form on $A$ is a linear map $\tau:A\to\mathbb K$ such that
$\tau(ab)=\tau(ba)$ for all $(a,b)\in A^2$. The data of a symmetric
and associative pairing is equivalent to the data of a trace form:
to the trace form $\tau$ corresponds the pairing $\varpi:(a,b)
\mapsto\tau(ab)$, and conversely $\tau$ is given by $\tau=\varpi^\flat(1)$.
One says that an algebra $A$ is a Frobenius algebra if it can be endowed
with an associative and perfect pairing; if one can choose this pairing
symmetric, then one calls $A$ a symmetric algebra.

Now let $A$ be a such a symmetric algebra, endowed with a symmetric,
associative and perfect pairing $\varpi$. Then the graded bialgebra
$\mathscr F(A)$ is endowed with the symmetric and perfect pairing
$\varpi_{\mathrm{tot}}$ (Section~\ref{ss:DualFunctF}) and each graded
piece $\mathscr F_n(A)$ is an associative algebra for the internal
product $\cdot$ (Section~\ref{ss:TwistedGpRing}).

\begin{proposition}
\label{pr:FnAIsSymm}
For any degree $n$, the pairing $\varpi_{\mathrm{tot}}\bigr|_{\mathscr
F_n(A)\times\mathscr F_n(A)}$ is associative and endows $\mathscr F_n(A)$
with the structure of a symmetric algebra.
\end{proposition}
\begin{proof}
We denote the linear form $\varpi^\flat(1)$ by $\tau$ and define a
linear form $\tau_{\mathrm{tot}}:\mathscr F(A)\to\mathbb K$ by setting
$$\tau_{\mathrm{tot}}\bigl((a_1\otimes a_2\otimes\cdots\otimes
a_n)\#\sigma\bigr)=\begin{cases}\prod_{i=1}^n\tau(a_i)&\text{if
$\sigma=e_n$,}\\0&\text{otherwise,}\end{cases}$$
for any $n\in\mathbb N$, any $(a_1,a_2,\ldots,a_n)\in A^n$ and any
$\sigma\in\mathfrak S_n$. A straightforward verification based on
Formula~(\ref{eq:DefPairingFV}) shows that $\varpi_{\mathrm{tot}}(x,y)=
\tau_{\mathrm{tot}}(x\cdot y)$ for any $(x,y)\in\mathscr F(A)^2$. It
follows in particular that the pairing $\varpi_{\mathrm{tot}}
\bigr|_{\mathscr F_n(A)\times\mathscr F_n(A)}$ is associative. Since
this pairing is also symmetric and perfect, the algebra $\mathscr
F_n(A)$ is a symmetric algebra.
\end{proof}

We add to these ingredients the data of a graded subbialgebra $S$
of $\mathscr F(A)$, assumed to be a subalgebra of it for the internal
product. The polar $S^\circ$ of $S$ satisfies
$$\varpi_{\mathrm{tot}}(S,S\cdot S^\circ)=\varpi_{\mathrm{tot}}(S\cdot
S,S^\circ)\subseteq\varpi_{\mathrm{tot}}(S,S^\circ)=0,$$
so that $S\cdot S^\circ\subseteq S^\circ$. A similar argument shows
the inclusion $S^\circ\cdot S\subseteq S^\circ$, and we conclude that
$S^\circ$ is a two-sided internal $S$-submodule of $\mathscr F(A)$.
Assuming that $A$ is a projective $\mathbb K$-module and that $S$ is
a direct summand of the graded $\mathbb K$-module $\mathscr F(A)$, we
construct the diagram~(\ref{eq:AutodualDiagr}), with $\mathscr F(V)$
replaced by $\mathscr F(A)$; beside being a diagram of graded bialgebras,
it is then a diagram of two-sided internal $S$-submodules.

\subsection{The case of a group algebra}
\label{ss:CaseGpAlg}
Group algebras are at the same time cocommutative bialgebras and symmetric
algebras. They give therefore examples to which the constructions of
Sections~\ref{ss:SplitForm} and \ref{ss:FrobStruct} can be applied. We
study this situation here.

So let $\Gamma$ be a finite group. We endow the algebra $\mathbb K\Gamma$
with the pairing $\varpi$ defined by
$$\varpi(\gamma,\delta)=\begin{cases}1&\text{if $\gamma=\delta^{-1}$,}\\
0&\text{otherwise.}\end{cases}$$
This pairing is associative, symmetric and perfect; the corresponding
trace form $\tau=\varpi^\flat(1)$ is the linear form that maps an
element $\gamma\in\Gamma$ to $1$ if $\gamma$ is the unit and to $0$
otherwise. (One may observe that the familiar trace map of $\mathbb
K\Gamma$, i.e.\ the regular character of $\Gamma$, is a scalar multiple
of $\tau$.)

We now construct the graded bialgebra $\mathscr F(\mathbb K\Gamma)$
and endow it with the pairing of bialgebras $\varpi_{\mathrm{tot}}$.
By Proposition~\ref{pr:FnAIsSymm}, each graded component
$\mathscr F_n(\mathbb K\Gamma)$ is a symmetric algebra for the
pairing $\varpi_{\mathrm{tot}}\bigr|_{\mathscr F_n(\mathbb K\Gamma)
\times\mathscr F_n(\mathbb K\Gamma)}$. This property can also be
explained in the following way.

Let us first recall that the wreath product $\Gamma\wr\mathfrak S_n$
is the semidirect product $\Gamma^n\rtimes\mathfrak S_n$ for the usual
permutation action of $\mathfrak S_n$ on $\Gamma^n$. Thus an element
$\Gamma\wr\mathfrak S_n$ can always be written as the product of an
element of $\mathfrak S_n$ and an element of $\Gamma^n$, and the
commutation rule between these two kinds of elements is
$$\sigma\cdot(\gamma_1,\gamma_2,\ldots,\gamma_n)=(\gamma_{\sigma^{-1}(1)},
\gamma_{\sigma^{-1}(2)},\ldots,\gamma_{\sigma^{-1}(n)})\cdot\sigma.$$
A comparison with Equation~(\ref{eq:ProdTwistedGpRing}) which defines
the product in the twisted group ring\linebreak$\mathscr F_n(\mathbb
K\Gamma)=(\mathbb K\Gamma)^{\otimes n}\#(\mathbb K\mathfrak S_n)$
shows the existence of an isomorphism of algebras
$$\mathbb K\bigl[\Gamma\wr\mathfrak S_n\bigr]\to\mathscr F_n(\mathbb
K\Gamma),\ \bigl[(\gamma_1,\gamma_2,\ldots,\gamma_n)\cdot\sigma\bigr]
\mapsto\bigl[(\gamma_1\otimes\gamma_2\otimes\cdots\otimes\gamma_n)\#
\sigma\bigr].$$
Now the group algebra $\mathbb K\bigl[\Gamma\wr\mathfrak S_n\bigr]$ has
a standard structure of a symmetric algebra, whose trace form $\tau_n$
is given by
$$\tau_n\bigl[(\gamma_1,\gamma_2,\ldots,\gamma_n)\cdot\sigma\bigr]=
\begin{cases}1&\text{if $\gamma_1$, $\gamma_2$, \dots, $\gamma_n$ are
equal to the unit of $\Gamma$ and $\sigma=e_n$,}\\0&\text{otherwise.}
\end{cases}$$
Under the previous isomorphism, this trace form coincides with the
linear form $\tau_{\mathrm{tot}}\bigr|_{\mathscr F_n(\mathbb K\Gamma)}$
used in the proof of Proposition \ref{pr:FnAIsSymm}. We conclude that
the pairing $\varpi_{\mathrm{tot}}\bigr|_{\mathscr F_n(\mathbb K\Gamma)
\times\mathscr F_n(\mathbb K\Gamma)}$ on $\mathscr F_n(\mathbb K\Gamma)$
corresponds to the usual associative, symmetric and perfect pairing
on the group algebra $\mathbb K\bigl[\Gamma\wr\mathfrak S_n\bigr]$.

\section{A Solomon descent theory for the wreath products $G\wr\mathfrak S_n$}
\label{se:SoloTheoWrProd}
In this section, we study a particular case of the following problem,
inspired by Solomon's article~\cite{Solomon76}: given a finite group $H$,
is it possible to find a subalgebra of the group algebra $\mathbb KH$ of
which the representation ring of $H$ is a quotient? More precisely, we use
the theory developed in in the previous sections to give a positive answer
in the case where the group $H$ is the wreath product $G\wr\mathfrak S_n$
of the symmetric group with a finite abelian group $G$.

\subsection{Representation rings}
\label{ss:ReprRings}
We first set up the notation we plan to use concerning representation
rings. Let $H$ be a finite group. We denote the algebra of complex-valued
functions on $H$ by $\mathbb C^H$. The $\mathbb Z$-submodule of
$\mathbb C^H$ spanned by the characters of $H$ is called the ring
of complex linear representations of $H$ and is denoted by $R(H)$.
The involutive map $f\mapsto f^*$ which sends a function in $\mathbb C^H$
to its complex-conjugate leaves $R(H)$ stable. The assignment $H
\rightsquigarrow R(H)$ is a contravariant functor from the category
of finite groups to the category of commutative rings with involution.

Elements of $R(H)$ are usually called virtual characters. The set
$\Irr(H)$ of irreducible characters of $H$ is a basis of the
$\mathbb Z$-module $R(H)$. A virtual character is called effective if
all its coordinates with respect to the basis $\Irr(H)$ are positive.

The linear form on $\mathbb C^H$ that maps a function $f$ to the complex
number $\frac1{|H|}\sum_{h\in H}f(h)$ restricts to a $\mathbb Z$-valued
additive form $\varphi$ on $R(H)$, which is called the fundamental linear
form on $R(H)$. Its value at an irreducible character $\zeta\in\Irr(H)$ is
$1$ if $\zeta$ is the trivial character of $H$ and $0$ otherwise. We
define the fundamental bilinear form $\beta:R(H)\times R(H)\to\mathbb Z$
by $\beta(f,g)=\varphi(fg)$ for any $(f,g)\in R(H)^2$. The usual inner
product of characters is the bilinear form $(f,g)\mapsto\beta(f,g^*)$.
Given two irreducible characters $\zeta$ and $\psi$ in $\Irr(H)$, the
number $\beta(\zeta,\psi)$ is thus $1$ if $\zeta=\psi^*$ and $0$ otherwise.
As a consequence, the fundamental bilinear form $\beta$ is an associative,
symmetric and perfect pairing; endowed with it, $R(H)$ becomes a symmetric
commutative algebra.

We conclude this section with a proposition which is probably well-known.
\begin{proposition}
\label{pr:RepRingWoRad}
The representation ring $R(H)$ has trivial Jacobson radical.
\end{proposition}
\begin{proof}
Let $L$ be a number field big enough to contain all the roots of unity
of order $|H|$ in $\mathbb C$, and let $\mathscr O$ be the integral
closure of $\mathbb Z$ in $L$. Let $X$ the set of maximal ideals of
$\mathscr O$. Since $\mathscr O$ is a Dedekind ring, there holds
$$\bigcap_{\mathfrak m\in X}\mathfrak m=\{0\}.$$
For any $h\in H$, the evaluation $\zeta(h)$ of a virtual character
$\zeta\in R(H)$ at $h$ belongs to $\mathscr O$. The image of the
evaluation map $\ev_h:\zeta\mapsto\zeta(h)$ is therefore a subring of
$\mathscr O$, over which $\mathscr O$ is integral. This implies that for
any $\mathfrak m\in X$, the intersection $(\im\ev_h)\cap\mathfrak m$ is
a maximal ideal in $(\im\ev_h)$, and thus that the inverse image
$\ev_h^{-1}(\mathfrak m)$ is a maximal ideal of $R(H)$. The desired
result now follows from the equality
$$\bigcap_{h\in H}\bigcap_{\mathfrak m\in X}\ev_h^{-1}(\mathfrak m)=
\bigcap_{h\in H}\Biggl[\ev_h^{-1}\Biggl(\bigcap_{\mathfrak m\in X}
\mathfrak m\Biggr)\Biggr]=\bigcap_{h\in H}\ker\ev_h=\{0\},$$
because the Jacobson radical of $R(H)$ is the intersection of all its
maximal (left) ideals.
\end{proof}

\subsection{The characters of the wreath products $G\wr\mathfrak S_n$}
\label{ss:CharacWreathProd}
Let $G$ be a finite group, not necessarily abelian. We present in this
section Specht's results about the characters of the wreath products
$G\wr\mathfrak S_n$. Our presentation follows the appendix of
\cite{Macdonald80}, Appendix~B of Chap.~I in \cite{Macdonald95} and
and \S7 in \cite{Zelevinsky81}, to which we refer the reader for the
proofs.

The wreath product $G\wr\mathfrak S_n$ is the semidirect product
$G^n\rtimes\mathfrak S_n$ for the usual permutation action of
$\mathfrak S_n$ on $G^n$. (By convention, the notation $G\wr\mathfrak
S_0$ denotes the group with one element.) An element of $G\wr\mathfrak
S_n$ can always be written in two ways as the product of an element of
$\mathfrak S_n$ and an element of $G^n$, namely
$$\sigma\cdot(g_1,g_2,\ldots,g_n)=(g_{\sigma^{-1}(1)},g_{\sigma^{-1}(2)},
\ldots,g_{\sigma^{-1}(n)})\cdot\sigma.$$

Given a $\mathbb CG$-module $V$, we construct a complex representation
$\eta_n(V)$ of $G\wr\mathfrak S_n$ on the space $V^{\otimes n}$ by
letting a product $(g_1,\ldots,g_n)\cdot\sigma$ act on a pure tensor
$v_1\otimes\cdots\otimes v_n\in V^{\otimes n}$ by
$$((g_1,g_2,\ldots,g_n)\cdot\sigma)\cdot(v_1\otimes v_2\otimes\cdots
\otimes v_n)=(g_1\cdot v_{\sigma^{-1}(1)})\otimes(g_2\cdot
v_{\sigma^{-1}(2)})\otimes\cdots\otimes(g_n\cdot v_{\sigma^{-1}(n)}).$$
The character of $\eta_n(V)$ does not depend actually on $V$ but only of
its character; if $\zeta$ denotes the latter, then we will denote the
former by $\eta_n(\zeta)$. Two particular cases are worth mentioning.
\begin{itemize}
\item If $\gamma$ is a linear character of $G$, that is, a character of
degree $1$, then $\eta_n(\gamma)$ is the linear character $\bigl((g_1,
g_2,\ldots,g_n)\cdot\sigma\bigr)\mapsto\gamma(g_1g_2\cdots g_n)$ of
$G\wr\mathfrak S_n$.
\item If $\zeta$ is the regular character of $G$, then $\eta_n(\zeta)$
is the character induced from the trivial representation of the subgroup
$\mathfrak S_n$ to $G\wr\mathfrak S_n$.
\end{itemize}

Let $\mathbf c=(c_1,c_2,\ldots,c_k)$ be a composition of $n$. The Young
subgroup $\mathfrak S_{\mathbf c}$ of $\mathfrak S_n$ acts on $G^n$, and
the semidirect product $G^n\rtimes\mathfrak S_{\mathbf c}$ can be seen
as the subgroup of $G\wr\mathfrak S_n$ generated by $G^n$ and $\mathfrak
S_{\mathbf c}$. By analogy, we denote it by $G\wr\mathfrak S_{\mathbf c}$
and we call it a Young subgroup of $G\wr\mathfrak S_n$. The natural
isomorphism $\mathfrak S_{c_1}\times\mathfrak S_{c_2}\times\cdots\times
\mathfrak S_{c_k}\cong\mathfrak S_{\mathbf c}$ gives rise to an
isomorphism $(G\wr\mathfrak S_{c_1})\times(G\wr\mathfrak S_{c_2})\times
\cdots\times(G\wr\mathfrak S_{c_k})\cong(G\wr\mathfrak S_{\mathbf c})$.

A partition is an infinite non-increasing sequence $\lambda=(\lambda_1,
\lambda_2,\ldots)$ of non-negative integers, all of whose terms but a
finite number vanish. As usual, we denote the sum of the parts of
$\lambda$ by $|\lambda|$; if $|\lambda|=n$, then we say that $\lambda$
is a partition of $n$. To a partition $\lambda$ of $n$, we associate in
the usual way an irreducible complex representation $S_\lambda$ of
$\mathfrak S_n$, the so-called Specht module. Thus for instance the
characters of $S_{(n)}$ and $S_{(1,1,\ldots,1)}$ (with $n$ terms equal
to $1$) are the trivial and signature characters of $\mathfrak S_n$,
respectively.

An $\Irr(G)$-partition is a family $\boldsymbol\lambda=(\lambda_\gamma
)_{\gamma\in\Irr(G)}$ indexed by $\Irr(G)$ of partitions. The size of
an $\Irr(G)$-partition $\boldsymbol\lambda$ is the number
$\|\boldsymbol\lambda\|=\sum_{\gamma\in\Irr(G)}|\lambda_\gamma|$. We
define the dual of $\boldsymbol\lambda$ as the $\Irr(G)$-partition
$\boldsymbol\lambda^*=\bigl(\gamma\mapsto\lambda_{\gamma^*}\bigr)$.

Given an $\Irr(G)$-partition $\boldsymbol\lambda$ of size $n$, one
constructs a complex representation of $G\wr\mathfrak S_n$ as follows.
One enumerates the irreducible characters $\gamma_1$, $\gamma_2$, \dots,
$\gamma_r$ of $G$ and picks up $\mathbb CG$-modules $V_1$, $V_2$, \dots,
$V_r$ that afford them. Let us set $c_i=|\lambda_i|$. Since $\mathfrak
S_{c_i}$ is a quotient of $G\wr\mathfrak S_{c_i}$, we may view the
Specht module $S_{\lambda_{ \gamma_i}}$ as a representation of
$G\wr\mathfrak S_{c_i}$ and we may then multiply it by $\eta_{c_i}(V_i)$.
The outer product
$$\bigl(S_{\lambda_{\gamma_1}}\otimes\eta_{c_1}(V_1)\bigr)\otimes
\bigl(S_{\lambda_{\gamma_2}}\otimes\eta_{c_2}(V_2)\bigr)\otimes\cdots
\otimes\bigl(S_{\lambda_{\gamma_r}}\otimes\eta_{c_r}(V_r)\bigr)$$
is then a representation of $(G\wr\mathfrak S_{c_1})\times(G\wr\mathfrak
S_{c_2})\times\cdots\times(G\wr\mathfrak S_{c_r})\cong\bigl(G\wr\mathfrak
S_{(c_1,c_2,\ldots,c_r)}\bigr)$, which we can induce to $G\wr\mathfrak
S_n$. The result of this induction does not depend up to isomorphism on
the choice of the enumeration $\gamma_1$, $\gamma_2$, \dots, $\gamma_r$.
Its character depends therefore only of $\boldsymbol\lambda$; we denote
it by $\boldsymbol\chi^{\boldsymbol\lambda}$. The map
$\boldsymbol\lambda\mapsto\boldsymbol\chi^{\boldsymbol\lambda}$
affords a bijection from the set of $\Irr(G)$-partitions of size $n$ onto
the set $\Irr\bigl(G\wr\mathfrak S_n\bigr)$. The complex-conjugate of
the character $\boldsymbol\chi^{\boldsymbol\lambda}$ is the character
$\boldsymbol\chi^{\boldsymbol\lambda^*}$.

Each representation ring $R(G\wr\mathfrak S_n)$ is a ring endowed with
its fundamental linear and bilinear forms $\varphi$ and $\beta$, this
latter being an associative, symmetric and perfect pairing. Considering
all $n$ at the same time yields however extra structures. We consider
therefore the direct sum $\Rep(G)=\bigoplus_{n\geq0}R(G\wr\mathfrak S_n)$.

We define the induction product $\psi*\psi'$ of two characters $\psi$ of
$G\wr\mathfrak S_n$ and $\psi'$ of $G\wr\mathfrak S_{n'}$ as the induction
$$\psi*\psi'=\Ind_{(G\wr\mathfrak S_n)\times(G\wr\mathfrak S_{n'})}^{G\wr
\mathfrak S_{n+n'}}\bigl(\psi\otimes\psi'\bigr),$$
where $(G\wr\mathfrak S_n)\times(G\wr\mathfrak S_{n'})$ is viewed as
the subgroup $G\wr\mathfrak S_{(n,n')}$ of $G\wr\mathfrak S_{n+n'}$.
The bilinear extension of the external product to $\Rep(G)\times\Rep(G)$
endows the space $\Rep(G)$ with the structure of a graded associative
and commutative algebra.

Likewise, the restriction coproduct $\Delta(\psi)$ of a character
$\psi$ of $G\wr\mathfrak S_n$ is defined to be the sum over
$n'\in\{0,1,\ldots,n\}$ of the restrictions
$$\Res_{(G\wr\mathfrak S_{n'})\times(G\wr\mathfrak S_{n-n'})}^{G\wr
\mathfrak S_n}\bigl(\psi\bigr).$$
This notation implicitely identifies characters of the group
$(G\wr\mathfrak S_{n'})\times(G\wr\mathfrak S_{n-n'})$ with elements of
$\Rep_{n'}(G)\otimes\Rep_{n-n'}(G)$, so that $\Delta(\psi)\in\Rep(G)
\otimes\Rep(G)$. The linear extension of $\Delta$ to the whole space
$\Rep(G)$ endows the latter with the structure of a graded coassociative
and cocommutative coalgebra. Mackey's subgroup theorem implies that
$(\Rep(G),*,\Delta)$ is a graded commutative cocommutative bialgebra.

In order to make the situation more alike to the structures seen in
Sections~\ref{se:Duality} and \ref{se:InternProd}, we extend the
product and the fundamental linear and bilinear forms $\varphi_n$
and $\beta_n$ defined on each $R(G\wr\mathfrak S_n)$ to operations
defined on the whole space $\Rep(G)$ by setting
$$fg=\sum_{n\geq0}f_ng_n,\quad\varphi_{\mathrm{tot}}(f)=\sum_{n\geq0}
\varphi_n(f_n)\quad\text{and}\quad\beta_{\mathrm{tot}}(f,g)=
\varphi_{\mathrm{tot}}(fg)=\sum_{n\geq0}\beta_n(f_n,g_n)$$
for any $f=\sum_{n\geq0}f_n$ and $g=\sum_{n\geq0}g_n$, where $f_n$
and $g_n$ in $R(G\wr\mathfrak S_n)$. Then $\beta_{\mathrm{tot}}$ is
a perfect symmetric pairing on $\Rep(G)$, with respect to which the
induction product $*$ and the restriction coproduct $\Delta$ are
adjoint to each other by Frobenius reciprocity. Moreover, Mackey's
tensor product theorem (more precisely, the particular case stated
in Corollary~(10.20) of~\cite{Curtis-Reiner81}) implies the following
splitting formula: for any $f$, $g_1$, $g_2$, \dots, $g_l$ in $\Rep(G)$,
there holds
\begin{equation}
f(g_1*g_2*\cdots*g_l)=\sum_{(f)}(f_{(1)}g_1)*(f_{(2)}g_2)*\dots*
(f_{(l)}g_l).
\label{eq:SplitFormRep}
\end{equation}

We denote by $\Lambda$ the ring of symmetric functions. This is indeed
a graded bialgebra (see I, 5, Ex.~25 in \cite{Macdonald95}). As is
well-known, the complete symmetric functions $h_n$ are algebraically
independent generators of the commutative $\mathbb Z$-algebra $\Lambda$.
On the other hand, the Schur functions $s_\lambda$, where $\lambda$ is
a partition, is a basis of the $\mathbb Z$-module $\Lambda$. Let
$\Lambda(\Irr(G))$ be the tensor product of a family
$(\Lambda(\gamma))_{\gamma\in\Irr(G)}$ of copies of $\Lambda$. For any
$\gamma\in\Irr(G)$, we denote by $P(\gamma)$ the element in the tensor
factor $\Lambda(\gamma)$ that corresponds to the symmetric function
$P\in\Lambda$. Given an $\Irr(G)$-partition $\boldsymbol\lambda=
(\lambda_\gamma)_{\gamma\in\Irr(G)}$, we set
$$\mathbf s_{\boldsymbol\lambda}=\prod_{\gamma\in\Irr(G)}
s_{\lambda_\gamma}(\gamma).$$
Then the elements $h_n(\gamma)$ are algebraically independent generators
of the commutative $\mathbb Z$-algebra $\Lambda(\Irr(G))$, where $n\geq1$
and $\gamma\in\Irr(G)$, and the elements $\mathbf s_{\boldsymbol\lambda}$
form a basis of the $\mathbb Z$-module $\Lambda(\Irr(G))$, where
$\boldsymbol\lambda$ is an $\Irr(G)$-partition. Finally we endow the
graded bialgebra $\Lambda(\Irr(G))$ with a perfect symmetric pairing
$\langle?,?\rangle$ defined on the basis of Schur functions by
$$\langle\mathbf s_{\boldsymbol\lambda},\mathbf s_{\boldsymbol\lambda'}
\rangle=\begin{cases}1&\text{if $\boldsymbol\lambda'=\boldsymbol
\lambda^*$,}\\0&\text{otherwise.}\end{cases}$$

Our interest in $\Lambda(\Irr(G))$ is that it gives a model that allows to
calculate in $\Rep(G)$. More precisely, there is an isomorphism of graded
bialgebras $\ch:\Rep(G)\stackrel\simeq\longrightarrow\Lambda(\Irr(G))$,
called the Frobenius characteristic, such that
$$\ch\bigl(\eta_n(\gamma)\bigr)=h_n(\gamma)\quad\text{and}\quad
\ch\bigl(\boldsymbol\chi^{\boldsymbol\lambda}\bigr)=\mathbf
s_{\boldsymbol\lambda}$$
for any $n\geq1$, any $\gamma\in\Irr(G)$ and any $\Irr(G)$-partition
$\boldsymbol\lambda$. Moreover $\ch$ is compatible in the obvious sense
with the perfect symmetric pairings $\beta_{\mathrm{tot}}$ on $\Rep(G)$
and $\langle?,?\rangle$ on $\Lambda(\Irr(G))$.

What precedes implies that
$$\Rep(G)\cong\bigotimes_{\gamma\in\Irr(G)}\mathbb Z\bigl[\eta_1
(\gamma),\eta_2(\gamma),\ldots\bigr].$$
The following proposition, which will be used in Section
\ref{ss:PartCaseZ/2Z}, explains how to find the expression of
$\eta_n(\zeta)$ as a polynomial in the $\eta_n(\gamma)$ when the
effective character $\zeta$ is not irreducible. In order to state it,
we introduce a last notation: viewing the signature character $\sgn$ of
$\mathfrak S_n$ as a character of $G\wr\mathfrak S_n$ through the
quotient map $(G\wr\mathfrak S_n)\to\mathfrak S_n$, we denote the
product $\sgn\cdot\eta_n(\zeta)$ in the ring $R(G\wr\mathfrak S_n)$
by $\varepsilon_n(\zeta)$.

\begin{proposition}
\label{pr:DescrEta}
There exists a morphism of groups $H:R(G)\to\bigl(\Rep(G)[[u]]
\bigr)^\times$ such that
\begin{equation}
H(\zeta)=\sum_{n\geq0}\eta_n(\zeta)u^n\quad\text{and}\quad
H(-\zeta)=\sum_{n\geq0}(-1)^n\varepsilon_n(\zeta)u^n
\label{eq:DescrEta}
\end{equation}
for all effective characters $\zeta$.
\end{proposition}
\begin{proof}
We extend the Frobenius characteristic $\ch$ to an isomorphism of rings
from $\Rep(G)[[u]]$ onto $\Lambda(\Irr(G))[[u]]$. Since $R(G)$ is a free
$\mathbb Z$-module with basis $\Irr(G)$, there exists an homomorphism
of abelian group $H:R(G)\to\bigl(\Rep(G)[[u]]\bigr)^\times$ such that
for each $\gamma\in\Irr(G)$,
$$H(\gamma)=\sum_{n\geq0}\eta_n(\gamma)u^n.$$

Now let $\zeta=\sum_{\gamma\in\Irr(G)}a_\gamma\gamma$ be an effective
character of $G$. A slight modification of the calculation made in
Appendix~B of Chap.~I, (8.3) in \cite{Macdonald95} yields
$$\ch\left(\sum_{n\geq0}\eta_n(\zeta)u^n\right)=\prod_{\gamma\in\Irr(G)}
\left(\sum_{n\geq0}h_n(\gamma)u^n\right)^{a_\gamma}$$
and
\begin{align}
\ch\left(\sum_{n\geq0}(-1)^n\varepsilon_n(\zeta)u^n\right)
&=\prod_{\gamma\in\Irr(G)}\left(\sum_{n\geq0}(-1)^ne_n(\gamma)u^n
\right)^{a_\gamma}\notag\\[3pt]
&=\prod_{\gamma\in\Irr(G)}\left(\sum_{n\geq0}h_n(\gamma)u^n
\right)^{-a_\gamma},\label{eq:PfDescrEta}
\end{align}
where $e_n\in\Lambda$ is the symmetric elementary function of degree
$n$. It follows that
\begin{align*}
\ch\left(\sum_{n\geq0}\eta_n(\zeta)u^n\right)&=\prod_{\gamma\in\Irr(G)}
\left(\sum_{n\geq0}\ch\bigl(\eta_n(\gamma)\bigr)u^n\right)^{a_\gamma}\\[3pt]
&=\ch\left(\prod_{\gamma\in\Irr(G)}\left(\sum_{n\geq0}\eta_n(\gamma)
u^n\right)^{a_\gamma}\right)\\[3pt]
&=\ch\bigl(H(\zeta)\bigr),
\end{align*}
and likewise $\ch\left(\sum_{n\geq0}(-1)^n\varepsilon(\zeta)u^n\right)
=\ch\bigl(H(-\zeta)\bigr)$. We conclude that (\ref{eq:DescrEta}) holds,
as required.
\end{proof}

\subsection{The Solomon homomorphism}
\label{ss:SoloHomo}
The representation theory presented in Section~\ref{ss:CharacWreathProd}
allows the use of the model $\Lambda(\Irr(G))$ to compute the
character tables of all the groups $G\wr\mathfrak S_n$ and to study the
inductions and the restrictions with respect to the Young subgroup.
However $\Lambda(\Irr(G))$ does not make the computation of the ring
structure of $R(G\wr\mathfrak S_n)$ particularly easy. In this section,
we construct explicitly a surjective ring homomorphism from a subring of
$\mathbb Z\bigl[G\wr\mathfrak S_n\bigr]$ onto $R(G\wr\mathfrak S_n)$.
However, we must restrict ourselves to the case where $G$ is abelian.
As usual, it is convenient to do this simultaneously for all $n$.

In this section and in the following one, $\mathbb K$ is the ring
$\mathbb Z$. Some variants are indeed possible, but this choice
simplifies slightly the notation. The letter $G$ denotes a finite
abelian group. The dual group of $G$, denoted by $G^\wedge$ or by
$\Gamma$, is the set $\Irr(G)$ endowed with the ordinary product of
characters. Although $G$ and $\Gamma$ are isomorphic as abstract
groups, we do not identify them. On the other hand, we observe that
the group ring $\mathbb Z\Gamma$ coincides with the representation
ring $R(G)$; indeed even the pairings and the trace forms which turn
these rings into symmetric $\mathbb Z$-algebras agree.

We construct the graded bialgebra $\mathscr F(\mathbb Z\Gamma)$ with
its external product $*$ and its coproduct $\Delta$; it is further
endowed with the internal product $\cdot$, the linear form
$\tau_{\mathrm{tot}}$ and the pairing $\varpi_{\mathrm{tot}}$ (see
Section~\ref{ss:CaseGpAlg}). On the other hand, we have the graded
bialgebra $\Rep(G)$ with the induction product~$*$ and the coproduct
$\Delta$, with also the fundamental linear and bilinear forms
$\varphi_{\mathrm{tot}}$ and $\beta_{\mathrm{tot}}$; moreover the
graded components $R(G\wr\mathfrak S_n)$ of $\Rep(G)$ are symmetric
algebras. Our aim now is to show that $\Rep(G)$ is a subquotient of
$\mathscr F(\mathbb Z\Gamma)$.

Since $\mathbb Z\Gamma$ is a cocommutative bialgebra, the
Mantaci-Reutenauer subbialgebra $\MR(\mathbb Z\Gamma)$ of $\mathscr
F(\mathbb Z\Gamma)$ is defined. This is a graded subbialgebra,
whose homogeneous component of degree $n$, say, will be denoted by
$\MR_n(\mathbb Z\Gamma)$. By Corollary~\ref{co:SubalgIntSigmaB}, each
$\MR_n(\mathbb Z\Gamma)$ is a subalgebra of $\mathscr F_n(\mathbb
Z\Gamma)$ for the internal product. Moreover, it follows from
Proposition~\ref{pr:SigmaIsFreeAlg} that with respect to the external
product, the associative algebra $\MR(\mathbb Z\Gamma)$ is freely
generated by the elements $y_{n,\gamma}$, where $n\geq1$ and
$\gamma\in\Gamma$. Thus there is a unique morphism of algebras
$\theta_G:\MR(\mathbb Z\Gamma)\to\Rep(G)$ that maps $y_{n,\gamma}$ to
$\eta_n(\gamma)$. We call this map $\theta_G$ the Solomon homomorphism.

\begin{theorem}
\label{th:SoloHomo}
\begin{enumerate}
\item\label{it:ThSHa}The Solomon homomorphism $\theta_G$ is a surjective
homomorphism of graded bialgebras with respect to the products $*$
and the coproducts $\Delta$ on $\MR(\mathbb Z\Gamma)$ and $\Rep(G)$;
its kernel is the ideal generated by the elements $(y_{m,\gamma}*
y_{n,\delta}-y_{n,\delta}*y_{m,\gamma})$, where $m\geq1$, $n\geq1$,
and $(\gamma,\delta)\in\Gamma^2$.
\item\label{it:ThSHb}For every degree $n$, the restriction of the
Solomon map $\theta_G:\MR_n(\mathbb Z\Gamma)\to R(G\wr\mathfrak S_n)$
is a surjective homomorphism of rings; its kernel is the Jacobson
radical of the ring $\MR(\mathbb Z\Gamma)$.
\item\label{it:ThSHc}The Solomon homomorphism
is compatible with the linear and bilinear forms $\tau_{\mathrm{tot}}$
and $\varpi_{\mathrm{tot}}$ on $\MR(\mathbb Z\Gamma)$ and
$\varphi_{\mathrm{tot}}$ and $\beta_{\mathrm{tot}}$ on
$R(G\wr\mathfrak S_n)$, in the sense that
\begin{equation}
\tau_{\mathrm{tot}}=\varphi_{\mathrm{tot}}\circ\theta_G\quad\text{and}
\quad\varpi_{\mathrm{tot}}=\beta_{\mathrm{tot}}\bigl(\theta_G(?),
\theta_G(?)\bigr).\label{eq:ThSHc}
\end{equation}
The kernel of $\theta_G$ is equal to the kernel $\MR(\mathbb Z\Gamma)
\cap\MR(\mathbb Z\Gamma)^\circ$ of the pairing $\varpi_{\mathrm{tot}}
\bigr|_{\MR(\mathbb Z\Gamma)\times\MR(\mathbb Z\Gamma)}$, where the
polar $\MR(\mathbb Z\Gamma)^\circ$ is defined in the ambient space
$\mathscr F(\mathbb Z\Gamma)$ with respect to the perfect pairing
$\varpi_{\mathrm{tot}}$.
\end{enumerate}
\end{theorem}
\begin{proof}
\begin{enumerate}
\item The algebra $\MR(\mathbb Z\Gamma)$ is the free associative
$\mathbb Z$-algebra generated by the elements $y_{n,\gamma}$, where
$n\geq1$ and $\gamma\in\Gamma$, whilst $\Rep(G)$ is the free
associative commutative $\mathbb Z$-algebra generated by the elements
$\eta_n(\gamma)$. It follows that $\theta_G$ is surjective and that
its kernel is the ideal generated by the commutators $(y_{m,\gamma}*
y_{n,\delta}-y_{n,\delta}*y_{m,\gamma})$. Moreover $\theta_G$ is
graded, for $y_{n,\gamma}$ and $h_n(\gamma)$ have both degree $n$.

It is easy to see that
$$\Res_{(G\wr\mathfrak S_{n'})\times(G\wr\mathfrak S_{n-n'})}^{G\wr
\mathfrak S_n}\bigl(\eta_n(\gamma)\bigr)=\eta_{n'}(\gamma)\otimes
\eta_{n-n'}(\gamma)$$
for any $\gamma\in\Gamma$ and any integers $n$ and $n'$ with $0\leq
n'\leq n$, and therefore
\begin{equation}
\Delta\bigl(\eta_n(\gamma)\bigr)=\sum_{n'=0}^n\eta_{n'}(\gamma)
\otimes\eta_{n-n'}(\gamma)
\label{eq:CoprodEta}
\end{equation}
in $\Rep(G)$. It follows then by comparison with Equation
(\ref{eq:CoprodY}) that the set
$$\{x\in\MR(\mathbb Z\Gamma)\mid\Delta\circ\theta_G(x)=(\theta_G
\otimes\theta_G)\circ\Delta(x)\}$$
contains the elements $y_{n,\gamma}$. Since this set is a subalgebra, it
is the whole $\MR(\mathbb Z\Gamma)$. The compatibility of $\theta_G$ with
the counit is trivial, and we conclude that $\theta_G$ is a morphism of
coalgebras. Assertion~\ref{it:ThSHa} is proved.
\item We first prove that $\theta_G$ maps the internal product of
$\MR_n(\mathbb Z\Gamma)$ to the ordinary product of characters in
$R(G\wr\mathfrak S_n)$. This fact may be shown by a direct computation
using Mantaci and Reutenauer's rule (Corollary~\ref{co:MantReutRule})
and Mackey's tensor product theorem; it may also be obtained by
the following reasoning, that is actually grounded on the same
combinatorial foundations.

A straightforward calculation, based on Equations~(\ref{eq:SplitForm})
and (\ref{eq:SplitFormRep}) and on the fact that $\theta_G$ is a morphism
of bialgebras for the operations $*$ and $\Delta$, shows that
$$E=\{z\in\MR(\mathbb Z\Gamma)\mid\forall y\in\MR(\mathbb Z\Gamma),\
\theta_G(y\cdot z)=\theta_G(y)\theta_G(z)\}$$
is a subalgebra of $\MR(\mathbb Z\Gamma)$ for the external product $*$.
On the other hand, every generator $y_{n,\delta}$ of $\MR(\mathbb Z\Gamma)$
belongs to $E$. Indeed any element in $\MR_n(\mathbb Z\Gamma)$ is a
linear combination of elements of the form
$$y=y_{c_1,\gamma_1}*y_{c_2,\gamma_2}*\cdots*y_{c_k,\gamma_k},$$
where $\mathbf c=(c_1,c_2,\ldots,c_k)$ is a composition of $n$ and
$\gamma_1$, $\gamma_2$, \dots, $\gamma_k$ are elements of $\Gamma$, and for
such a $y$, Formulas~(\ref{eq:SplitFormRep}) and (\ref{eq:CoprodEta}) imply
\begin{align*}
\theta_G(y\cdot y_{n,\delta})
&=\theta_G\bigl(x_{\mathbf c}\cdot\bigl[(\gamma_1^{\otimes c_1}\otimes
\gamma_2^{\otimes c_2}\otimes\cdots\otimes\gamma_k^{\otimes c_k})
\#e_n\bigr]\cdot(\delta^{\otimes n}\#e_n)\bigr)\\[3pt]
&=\theta_G\bigl(x_{\mathbf c}\cdot\bigl[((\gamma_1\delta)^{\otimes c_1}
\otimes(\gamma_2\delta)^{\otimes c_2}\otimes\cdots\otimes
(\gamma_k\delta)^{\otimes c_k})\#e_n\bigr]\bigr)\\[3pt]
&=\theta_G\bigl(y_{c_1,\gamma_1\delta}*y_{c_2,\gamma_2\delta}*\cdots*
y_{c_k,\gamma_k\delta}\bigr)\\[3pt]
&=\eta_{c_1}(\gamma_1\delta)*\eta_{c_2}(\gamma_2\delta)*\cdots*
\eta_{c_k}(\gamma_k\delta)\\[3pt]
&=\bigl(\eta_{c_1}(\gamma_1)\eta_{c_1}(\delta)\bigr)*
\bigl(\eta_{c_2}(\gamma_2)\eta_{c_2}(\delta)\bigr)*\cdots*
\bigl(\eta_{c_k}(\gamma_k)\eta_{c_k}(\delta)\bigr)\\[3pt]
&=\sum_{(\eta_n(\delta))}
\bigl(\eta_{c_1}(\gamma_1)(\eta_n(\delta))_{(1)}\bigr)*
\bigl(\eta_{c_2}(\gamma_2)(\eta_n(\delta))_{(2)}\bigr)*\cdots*
\bigl(\eta_{c_k}(\gamma_k)(\eta_n(\delta))_{(k)}\bigr)\\[3pt]
&=\bigl(\eta_{c_1}(\gamma_1)*\eta_{c_2}(\gamma_2)*\cdots*
\eta_{c_k}(\gamma_k)\bigr)\;\eta_n(\delta)\\[3pt]
&=\theta_G(y)\;\theta_G(y_{n,\delta}).
\end{align*}
Therefore $E=\MR(\mathbb Z\Gamma)$. Observing moreover that $\theta_G$
maps the unit element of $\MR_n(\mathbb Z\Gamma)$, namely $(1^{\otimes
n}\#e_n)=y_{n,1}$, to the unit of $R(G\wr\mathfrak S_n)$, namely the
trivial character $\eta_n(1)$ of $G\wr\mathfrak S_n$, we conclude that
the degree $n$ part of $\theta_G$ is an homomorphism of rings from
$\MR_n(\mathbb Z\Gamma)$ to $R(G\wr\mathfrak S_n)$.

Assertion~\ref{it:ThSHa} implies that this homomorphism is surjective,
which entails that the Jacobson radical of $\MR_n(\mathbb Z\Gamma)$ is
contained in the preimage by $\theta_G\bigm|_{\MR_n(\mathbb Z\Gamma)}$
of the Jacobson radical of $R(G\wr\mathfrak S_n)$. By Proposition
\ref{pr:RepRingWoRad}), this translates readily into the inclusion
$\rad\MR_n(\mathbb Z\Gamma)\subseteq\ker\theta_G\bigm|_{\MR_n(\mathbb
Z\Gamma)}$.

To prove the reverse inclusion, we will use the result stated in
Assertion~\ref{it:ThSHc}. (Though its validity has not yet been
established, no vicious circle arises in the reasoning.) So let us
suppose that some element $x\in\MR_n(\mathbb Z\Gamma)$ belongs to the
kernel of $\theta_G$. This element $x$ acts by left multiplication on
the algebra $\mathscr F_n(\mathbb Z\Gamma)$. Since this latter is a
free $\mathbb Z$-module, this action can be represented by a matrix
with entries in $\mathbb Z$. For any positive integer $k$, the $k$-th
power of this matrix represents the action of the left multiplication
by $x^k$ and therefore its trace is
$$\rk\mathscr F_n(\mathbb Z\Gamma)\;\tau_{\mathrm{tot}}(x^k)$$
by the interpretation of $\tau_{\mathrm{tot}}$ given at the end of
Section~\ref{ss:CaseGpAlg}. However our assumption that $x\in\ker\theta_G$
and Assertion~\ref{it:ThSHc} yield
$$\tau_{\mathrm{tot}}(x^k)=\varpi_{\mathrm{tot}}(x,x^{k-1})=
\beta_{\mathrm{tot}}\bigl(\theta_G(x),\theta_G\bigl(x^{k-1}
\bigr)\bigr)=0$$
for all $k\geq1$. It follows that our matrix is nilpotent, and therfore
that $x$ itself is nilpotent. This argument shows that all elements of
the ideal $\ker\theta_G\bigm|{\MR_n(\mathbb Z\Gamma)}$ of
$\MR_n(\mathbb Z\Gamma)$ are nilpotent. This kernel is thus contained
in the radical of $\MR_n(\mathbb Z\Gamma)$, which completes the proof
of Assertion~\ref{it:ThSHb}.

\item Elements of the form
$$x=y_{c_1,\gamma_1}*y_{c_2,\gamma_2}*\cdots*y_{c_k,\gamma_k}$$
span the $\mathbb Z$-module $\MR(\mathbb Z\Gamma)$. Putting
$n=c_1+c_2+\cdots+c_k$ and $\mathbf c=(c_1,c_2,\ldots,c_k)$, we
observe that
$$\tau_{\mathrm{tot}}(x)=\tau_{\mathrm{tot}}\bigl(x_{\mathbf c}\cdot
\bigl[(\gamma_1^{\otimes c_1}\otimes\gamma_2^{\otimes c_2}\otimes
\cdots\otimes\gamma_k^{\otimes c_k})\#e_n\bigr]\bigr)$$
is $1$ if all the elements $\gamma_i$ are equal to $1$ and is $0$
otherwise. On the other hand, Frobenius reciprocity implies that
\begin{align*}
\varphi_{\mathrm{tot}}\circ\theta_G(x)
&=\dim\Hom_{G\wr\mathfrak S_n}\bigl(\eta_{c_1}(\gamma_1)*\eta_{c_2}
(\gamma_2)*\cdots*\eta_{c_k}(\gamma_k),1_{G\wr\mathfrak S_n}\bigr)\\[2pt]
&=\dim\Hom_{G\wr\mathfrak S_n}\bigl(\Ind_{G\wr\mathfrak S_{\mathbf
c}}^{G\wr\mathfrak S_n}\bigl(\eta_{c_1}(\gamma_1)\otimes\eta_{c_2}
(\gamma_2)\otimes\cdots\otimes\eta_{c_k}(\gamma_k)\bigr),
1_{G\wr\mathfrak S_n}\bigr)\\[2pt]
&=\dim\Hom_{G\wr\mathfrak S_{\mathbf c}}\bigl(\eta_{c_1}(\gamma_1)
\otimes\eta_{c_2}(\gamma_2)\otimes\cdots\otimes\eta_{c_k}(\gamma_k),
1_{G\wr\mathfrak S_{\mathbf c}}\bigr).
\end{align*}
The character $\eta_{c_1}(\gamma_1)\otimes\eta_{c_2}(\gamma_2)\otimes
\cdots\otimes\eta_{c_k}(\gamma_k)$ of $G\wr\mathfrak S_{\mathbf c}$ is
one-dimensional. Therefore $\varphi_{\mathrm{tot}}\circ\theta_G(x)$ is
$1$ if this character is trivial, that is, if all the elements
$\gamma_i$ are equal to $1$, and is $0$ otherwise. The equality
$\tau_{\mathrm{tot}}(x)=\varphi_{\mathrm{tot}}\circ\theta_G(x)$ being
valid for each $x$ in a spanning set for $\MR(\mathbb Z\Gamma)$, we
conclude that $\tau_{\mathrm{tot}}=\varphi_{\mathrm{tot}}\circ\theta_G$.
In turn, this implies that
$$\varpi_{\mathrm{tot}}(x,y)=\tau_{\mathrm{tot}}(x\cdot y)=
\varphi_{\mathrm{tot}}\circ\theta_G(x\cdot y)=
\varphi_{\mathrm{tot}}(\theta_G(x)\theta_G(y))=
\beta_{\mathrm{tot}}(\theta_G(x),\theta_G(y))$$
for any $(x,y)\in\MR(\mathbb Z\Gamma)^2$.

An immediate consequence of this last equality is that $\ker\theta_G$
is contained in the kernel $\MR(\mathbb Z\Gamma)\cap\MR(\mathbb
Z\Gamma)^\circ$ of the symmetric pairing $\varpi_{\mathrm{tot}}
\bigr|_{\MR(\mathbb Z\Gamma)\times\MR(\mathbb Z\Gamma)}$. The reverse
inclusion holds also because $\theta_G$ is surjective and
$\beta_{\mathrm{tot}}$ is a perfect pairing on $\Rep(G)$.
Assertion~\ref{it:ThSHc} is proved.
\end{enumerate}
\end{proof}

Assertion~\ref{it:ThSHb} of Theorem~\ref{th:SoloHomo} says that
the representation ring $R(G\wr\mathfrak S_n)$ can be obtained as
a quotient of the subring $\MR_n(\mathbb Z\Gamma)$ of $\mathscr
F_n(\mathbb Z\Gamma)\cong\mathbb Z\bigl[\Gamma\wr\mathfrak S_n\bigr]$.
Since $\Gamma$ and $G$ are isomorphic, this entails that the
representation ring of the group $G\wr\mathfrak S_n$ can be realized
as a quotient of a subring of its group algebra. In other words, there
exists a Solomon descent theory for the wreath product $G\wr\mathfrak S_n$.
However it is not canonical, for it depends on the choice of an
isomorphism between $G$ and its dual.

The notation used above suggests the existence of some kind of
functoriality. In order to state a precise statement, we define a
category $\mathscr V$. Objects of $\mathscr V$ are $\mathbb N$-graded
abelian groups $A=\bigoplus_{n\geq0}A_n$; each graded piece $A_n$ is
further endowed with the structure of a ring, and the whole space $A$
is endowed with the structure of a graded $\mathbb Z$-bialgebra through
another, graded product, a unit, a coproduct and a counit. Morphisms in
$\mathscr V$ are maps that respect the $\mathbb N$-graduation, all
products with their units and the coproduct with its counit. In the
statement below, we denote the dual of a finite abelian group $G$ by
$G^\wedge$; the assignment $G\rightsquigarrow G^\wedge$ is a
contravariant endofunctor of the category of finite abelian groups.

\begin{proposition}
\label{pr:SolomonIsNatural}
The assignments $G\rightsquigarrow\mathscr F(\mathbb Z[G^\wedge])$,
$G\rightsquigarrow\MR(\mathbb Z[G^\wedge])$ and $G\rightsquigarrow
\Rep(G)$ are contravariant functors from the category of finite
abelian groups to the category $\mathscr V$. The assignment
$G\rightsquigarrow\theta_G$ is a natural transformation from
$\MR(\mathbb Z[?^\wedge])$ to $\Rep(?)$.
\end{proposition}

We leave the proof as a (rather tedious) exercise. The naturality of
$G\rightsquigarrow\theta_G$ means that for each morphism $f:G\to G'$
between two finite abelian groups, the diagram
$$\xymatrix@C+=70pt{\mathscr F(\mathbb Z[G^\wedge])&\mathscr
F(\mathbb Z[(G')^\wedge])\ar[l]_{\mathscr F(\mathbb Z[f^\wedge])}\\
\MR(\mathbb Z[G^\wedge])\ar@{^{(}->}[]+U+/u2pt/;[u]\ar@{->>}
[d]_{\theta_G}&\MR(\mathbb Z[(G')^\wedge])\ar@{^{(}->}[]+U+/u2pt/;[u]
\ar@{->>}[d]^{\theta_{G'}}\ar[l]_{\MR(\mathbb Z[f^\wedge])}\\
\Rep(G)&\Rep(G')\ar[l]_{\Rep(f)}}$$
is commutative.

To conclude this section, let us observe that Formula~(\ref{eq:ThSHc})
implies the commutativity of the diagram
\begin{equation}
\raisebox{44pt}{\xymatrix{&\mathscr F(\mathbb Z\Gamma)\ar[r]^{\simeq}
_{\varpi_{\mathrm{tot}}{}^\flat}&\mathscr F(\mathbb Z\Gamma)^\vee
\ar@{->>}[dr]&\\\MR(\mathbb Z\Gamma)\ar@{^{(}->}[ur]\ar@{->>}
[dr]_{\theta_G}&&&\MR(\mathbb Z\Gamma)^\vee.\\
&\Rep(G)\ar[r]^{\simeq}_{\beta_{\mathrm{tot}}{}^\flat}&\Rep(G)^\vee
\ar@{^{(}->}[ur]_{(\theta_G)^\vee}&}}\label{eq:DiagrThetaVee}
\end{equation}
(The surjectivity of the map $\mathscr F(\mathbb Z\Gamma)^\vee\to\MR(
\mathbb Z\Gamma)^\vee$ comes from Proposition~\ref{pr:MRIsDirSum}.)
By Theorem~\ref{th:SoloHomo}, $\theta_G$ is surjective with kernel
$\MR(\mathbb Z\Gamma)\cap\MR(\mathbb Z\Gamma)^\circ$, which implies
that $\theta_G$ defines an isomorphism of graded bialgebras
$$\overline{\theta_G}:\MR(\mathbb Z\Gamma)/(\MR(\mathbb Z\Gamma)\cap
\MR(\mathbb Z\Gamma)^\circ)\stackrel\simeq\longrightarrow\Rep(G).$$
We can therefore add an horizontal line in the middle of the diagram
(\ref{eq:DiagrThetaVee}) and get
$$\xymatrix@C-8.7pt{&\mathscr F(\mathbb Z\Gamma)\ar[r]^{\simeq}
_{\varpi_{\mathrm{tot}}{}^\flat}&\mathscr F(\mathbb Z\Gamma)^\vee
\ar@{->>}[dr]&\\\!\MR(\mathbb Z\Gamma)\ar@{^{(}->}[ur]\ar@{->>}[r]
\ar@{->>}[dr]_{\theta_G}&\MR(\mathbb Z\Gamma)/(\MR(\mathbb Z\Gamma)
\cap\MR(\mathbb Z\Gamma)^\circ)\ar[d]_*[left]{\scriptstyle\simeq}
^{\overline{\theta_G}}&\bigl(\MR(\mathbb Z\Gamma)/(\MR(\mathbb Z\Gamma)
\cap\MR(\mathbb Z\Gamma)^\circ)\bigr){}^\vee\ar@{^{(}->}[]+R+/r2pt/;[r]
&\MR(\mathbb Z\Gamma)^\vee.\\&\Rep(G)\ar[r]^\simeq_{\beta_{\mathrm{tot}}
{}^\flat}&\Rep(G)^\vee\ar[u]_*[right]{\scriptstyle\simeq}^{\overline{
\theta_G}{}^\vee}\ar@{^{(}->}[ur]_{(\theta_G)^\vee}&}$$
This is of course an occurence of the diagram~(\ref{eq:AutodualDiagr})
with $V=\mathbb Z\Gamma$ and $S=\MR(\mathbb Z\Gamma)$. As a bonus,
we see that the pairing induced by $\varpi_{\mathrm{tot}}$ on
$S/(S\cap S^\circ)$ is perfect in the present situation.

\begin{other}{Remark}
\label{rk:HivNovThib}
In this remark, we consider the case $G=\mathbb Z/r\mathbb Z$. Hiver,
Novelli and Thibon~\cite{Hivert-Novelli-Thibon04} have found an
embodiment of the lower half of the diagram~(\ref{eq:DiagrThetaVee})
in terms of the representation theory of a suitable limit at $q=0$ of
the Ariki-Koike algebra $\mathscr H_{n,r}(q)$. More precisely, these
authors propose to identify as $\mathbb Z$-modules the degree $n$
components $\MR_n(\mathbb Z\Gamma)$ and $\MR_n(\mathbb Z\Gamma)^\vee$
with the Grothendieck groups $K_0(\mathscr H_{n,r}(0))$ and
$G_0(\mathscr H_{n,r}(0))$, respectively. They claim that in this
identification, the map $\theta_G^\vee\circ\theta_G^{}$ coincides
with the Cartan homomorphism $c:K_0(\mathscr H_{n,r}(0))\to
G_0(\mathscr H_{n,r}(0))$, which describes the Jordan-Hölder
multiplicities of the simple modules in a projective module. They also
assert that the maps $\theta_G$ and $\theta_G^\vee$ can be interpreted
as arrows in a Cartan-Brauer $cde$ triangle
$$\xymatrix{K_0(\mathscr H_{n,r}(0))\ar[rr]^c\ar[dr]_e&&G_0(\mathscr
H_{n,r}(0))\\&K_0(\mathscr H_{n,r}(q))\ar[ru]_d&},$$
the bottom vertex being the Grothendieck group of the semisimple
category of finitely generated $\mathscr H_{n,r}(q)$-modules, where
$q$ is generic. Yet the fact, apparent in our constructions, that the
Cartan map $c$ is a morphism of $K_0(\mathscr H_{n,r}(0))$-bimodules
is missing in this picture.
\end{other}

\subsection{Symmetry property of the Solomon homomorphism}
\label{ss:SymmProp}
Given any finite group $H$, the data of a complex-valued function on $H$
is the same thing as the data of a $\mathbb Z$-linear map from $\mathbb
ZH$ into $\mathbb C$; we can therefore evaluate an element of $R(H)$ on
an element of $\mathbb ZH$. Applying this remark to the case of the group
$G\wr\mathfrak S_n$, we can evaluate an element $\theta_G(y)$ with
$y\in\MR_n(\mathbb Z\Gamma)$ on an element of $\mathbb Z\bigl[G\wr
\mathfrak S_n\bigr]=\mathscr F_n(\mathbb ZG)$, and in particular on an
element $x$ of $\MR_n(\mathbb ZG)$. Now $G$ and $\Gamma$ play symmetric
roles, so that the problem of comparing $\theta_G(y)(x)$ and
$\theta_\Gamma(x)(y)$ arises.

\begin{theorem}
\label{th:SymmSoloHomo}
For any $n\geq1$, $x\in\MR_n(\mathbb ZG)$ and $y\in\MR_n(\mathbb Z\Gamma)$,
there holds $\theta_G(y)(x)=\theta_\Gamma(x)(y)$.
\end{theorem}
\begin{proof}
In order to better put in evidence the symmetry between $G$ and $\Gamma$,
we denote the evaluation of a character $\gamma\in\Gamma$ at a point
$g\in G$ by a bracket $\langle\gamma,g\rangle$; the same notation can
then also be used to denote the evaluation of $g$, viewed as a character
of $\Gamma$, at the point $\gamma$.

We check the property asserted by the theorem for $x=y_{c_1,g_1}*
y_{c_2,g_2}*\cdots*y_{c_k,g_k}$ and $y=y_{d_1,\gamma_1}*y_{d_2,\gamma_2}*
\cdots*y_{d_k,\gamma_k}$, where $\mathbf c=(c_1,c_2,\ldots,c_k)$ and
$\mathbf d=(d_1,d_2,\ldots,d_l)$ are two compositions of $n$,
$(g_1,g_2,\ldots,g_k)\in G^k$ and $(\gamma_1,\gamma_2,\ldots,\gamma_l)
\in\Gamma^l$. Let us set
\begin{align*}
\mathbf{\tilde g}&=\tilde g_1\otimes\tilde g_2\otimes\cdots\otimes
\tilde g_n=(g_1^{\otimes c_1})\otimes(g_2^{\otimes c_2})\otimes\cdots
\otimes(g_k^{\otimes c_k}),\\[5pt]
\boldsymbol{\tilde\gamma}&=\tilde\gamma_1\otimes\tilde\gamma_2\otimes
\cdots\otimes\tilde\gamma_n=(\gamma_1^{\otimes d_1})\otimes(\gamma_2^{
\otimes d_2})\otimes\cdots\otimes(\gamma_l^{\otimes d_l}).
\end{align*}

We first compute $\theta_G(y)(x)$. Let $\rho\in X_{\mathbf c}$. Noting
that the composed map
$$\mathfrak S_n\hookrightarrow G\wr\mathfrak S_n\twoheadrightarrow(G
\wr\mathfrak S_n)/(G\wr\mathfrak S_{\mathbf d})$$
induces a bijection from $\mathfrak S_n/\mathfrak S_{\mathbf d}$ onto
$(G\wr\mathfrak S_n)/(G\wr\mathfrak S_{\mathbf d})$ and setting
$u_j=d_1+d_2+\cdots+d_j$, we compute
\begin{align*}
\theta_G(y)(\rho\cdot(\mathbf{\tilde g}\#e_n))
&=\Ind_{G\wr\mathfrak S_{\mathbf d}}^{G\wr\mathfrak S_n}
\bigl(\eta_{d_1}(\gamma_1)\otimes\eta_{d_2}(\gamma_2)\otimes\cdots
\otimes\eta_{d_l}(\gamma_l)\bigr)\bigl(\rho\cdot(\mathbf{\tilde
g}\#e_n)\bigr)\\[6pt]
&=\sum_{\substack{\pi\in\mathfrak S_n/\mathfrak S_{\mathbf d}\\[2pt]
\pi^{-1}\rho\pi\in\mathfrak S_{\mathbf d}}}\bigl(\eta_{d_1}(\gamma_1)
\otimes\eta_{d_2}(\gamma_2)\otimes\cdots\otimes\eta_{d_l}(\gamma_l)
\bigr)\bigl(\pi^{-1}\rho\cdot(\mathbf{\tilde g}\#e_n)\cdot\pi\bigr)\\[2pt]
&=\sum_{\substack{\pi\in\mathfrak S_n/\mathfrak S_{\mathbf d}\\[2pt]
\pi^{-1}\rho\pi\in\mathfrak S_{\mathbf d}}}\prod_{j=1}^l\
\langle\gamma_j,\tilde g_{\pi(u_{j-1}+1)}\tilde g_{\pi(u_{j-1}+2)}
\cdots\tilde g_{\pi(u_j)}\rangle\\[6pt]
&=\sum_{\substack{\pi\in X_{\mathbf d}\\[2pt]\rho\pi\mathfrak
S_{\mathbf d}=\pi\mathfrak S_{\mathbf d}}}\langle\tilde\gamma_1,
\tilde g_{\pi(1)}\rangle\langle\tilde\gamma_2,\tilde g_{\pi(2)}\rangle
\cdots\langle\tilde\gamma_n,\tilde g_{\pi(n)}\rangle.
\end{align*}
Taking the sum for all $\rho\in X_{\mathbf c}$, we find
\begin{equation}
\theta_G(y)(x)=\sum_{\rho\in X_{\mathbf c}}\theta_G(y)(\rho\cdot
(\mathbf{\tilde g}\#e_n))=\sum_{\substack{\rho\in X_{\mathbf c},\
\pi\in X_{\mathbf d}\\[2pt]\rho\pi\mathfrak S_{\mathbf d}=\pi\mathfrak
S_{\mathbf d}}}\langle\tilde\gamma_1,\tilde g_{\pi(1)}\rangle\langle
\tilde\gamma_2,\tilde g_{\pi(2)}\rangle\cdots\langle\tilde\gamma_n,
\tilde g_{\pi(n)}\rangle.
\label{eq:PfThSSH}
\end{equation}

In Section~\ref{ss:DblCosetSymGp}, we have parametrized double cosets
$C\in\mathfrak S_{\mathbf c}\backslash\mathfrak S_n/\mathfrak
S_{\mathbf d}$ by matrices $M\in\mathscr M_{\mathbf c,\mathbf d}$: to the
matrix $M=(m_{ij})$ corresponds the double coset $C(M)$. We aim now
at splitting the sum in the right-hand side of (\ref{eq:PfThSSH})
according to the double coset $C$ containing $\pi$. For that, we set
$$\mathcal F(\mathbf c,\mathbf d,M)=\bigl\{(\rho,\pi)\in
X_{\mathbf c}\times(X_{\mathbf d}\cap C(M))\bigm|\rho\pi
\mathfrak S_{\mathbf d}=\pi\mathfrak S_{\mathbf d}\bigr\}$$
and we observe that
$$\langle\tilde\gamma_1,\tilde g_{\pi(1)}\rangle\langle\tilde\gamma_2,
\tilde g_{\pi(2)}\rangle\cdots\langle\tilde\gamma_n,\tilde g_{\pi(n)}
\rangle=\prod_{i=1}^k\prod_{j=1}^l\;\langle g_i,\gamma_j\rangle^{m_{ij}}$$
if $\pi\in C(M)$. Then (\ref{eq:PfThSSH}) reads
$$\theta_G(y)(x)=\sum_{M\in\mathscr M_{\mathbf c,\mathbf d}}\Biggl[\;
\bigl|\mathcal F(\mathbf c,\mathbf d,M)\bigr|\;\prod_{i=1}^k\prod_{j=1}^l
\;\langle g_i,\gamma_j\rangle^{m_{ij}}\Biggr],$$
and by symmetry,
\begin{align*}
\theta_\Gamma(x)(y)
&=\sum_{N\in\mathscr M_{\mathbf d,\mathbf c}}\Biggl[\;\bigl|\mathcal
F(\mathbf d,\mathbf c,N)\bigr|\;\prod_{j=1}^l\prod_{i=1}^k
\;\langle\gamma_j,g_i\rangle^{n_{ji}}\Biggr]\\[3pt]
&=\sum_{M\in\mathscr M_{\mathbf c,\mathbf d}}\Biggl[\;\bigl|\mathcal
F(\mathbf d,\mathbf c,M^T)\bigr|\;\prod_{i=1}^k\prod_{j=1}^l\;\langle
g_i,\gamma_j\rangle^{m_{ij}}\Biggr],
\end{align*}
where $M^T$ denote the transpose of the matrix $M$. Observing now that
the double coset $C(M^T)\in\mathfrak S_{\mathbf d}\backslash\mathfrak
S_n/\mathfrak S_{\mathbf c}$ is equal to $C(M)^{-1}=\{\pi^{-1}\mid\pi\in
C(M)\}$, we deduce from Theorem~1.2 and Corollary~2.2 of
\cite{Blessenohl-Hohlweg-Schocker04}, applied to the group $\mathfrak
S_n$, that
$$\bigl|\mathcal F(\mathbf c,\mathbf d,M)\bigr|=\bigl|\mathcal
F(\mathbf d,\mathbf c,M^T)\bigr|$$
for each matrix $M\in\mathscr M_{\mathbf c,\mathbf d}$. The theorem
follows.
\end{proof}

This kind of question was first investigated by Jöllenbeck and Reutenauer
in \cite{Jollenbeck-Reutenauer01}; their result corresponds to the
(already non-trivial) case where $G$ is the group with one element. A
similar symmetry result holds also for the original Solomon descent
algebra and the original Solomon homomorphism of an arbitrary finite
Coxeter group (see~\cite{Blessenohl-Hohlweg-Schocker04}); the critical
point in the proof above is a theorem from this latter work.

\subsection{The particular case $G=\mathbb Z/2\mathbb Z$}
\label{ss:PartCaseZ/2Z}
In this section, we apply our results to the case where $G=\{\pm1\}$ is
the group with two elements. The peculiarity of this case is \vspace{2pt}
that $W_n=G\wr\mathfrak S_n$ is then the Coxeter group of type B$_n$.
Thus Solomon's constructions \cite{Solomon76} can be applied to it: there
is a certain subring $\tilde\Sigma$ of the group ring $\mathbb ZW_n$ and
a certain homomorphism of rings $\tilde\theta$ from $\tilde\Sigma$ to the
representation ring $R(W_n)$. This map $\tilde\theta$ is not surjective,
but Bonnafé and Hohlweg \cite{Bonnafe-Hohlweg04} manage to correct the
situation. They notice that the subring $\tilde\Sigma$ of $\mathbb
ZW_n\cong \mathscr F_n(\mathbb ZG)$ is contained in the Mantaci-Reutenauer
algebra $\MR_n(\mathbb ZG)$ and show how to extend $\tilde\theta$ to
$\MR_n(\mathbb ZG)$. The resulting map, still denoted by $\tilde\theta$,
is a surjective homomorphism of rings from $\MR_n(\mathbb ZG)$ onto $R(W_n)$.
The situation now looks like our Theorem~\ref{th:SoloHomo}~\ref{it:ThSHb},
which says that the homomorphism $\theta_G$ is a surjective ring
homomorphism from $\MR_n(\mathbb Z\Gamma)$ onto
$R(G\wr\mathfrak S_n)=R(W_n)$. Indeed we may identify $G$ and $\Gamma$
in the present case $G=\{\pm1\}$, because there is a unique isomorphism
between $G$ and $\Gamma$. Then both $\tilde\theta$ and $\theta_G$ are
surjective ring homomorphisms from $\MR_n(\mathbb Z\Gamma)$ onto $R(W_n)$;
our aim in this section is to explain the relationship between them.

We begin by setting the notation, following~\cite{Bonnafe-Hohlweg04}.
Let $n$ be a non-negative integer. We set $G=\{\pm1\}$ and $W_n=G\wr
\mathfrak S_n$. The group $W_n$ contains $\mathfrak S_n$ as a subgroup;
it is generated by the transpositions $s_i\in\mathfrak S_n$ (see the proof
of Corollary~\ref{co:CombiSymmGp}) and the element $(-1,1,1,\ldots,1)\in
G^n$. Endowed with this system of generators, $W_n$ becomes a Coxeter
system.

We agree to denote the subgroup $\mathfrak S_n$ of $W_n$ by the somewhat
strange convention $W_{-n}$. Likewise, we denote the trivial subgroup
with one element of $G^c$ by $G^{-c}$, for any positive integer $c$. We
define a signed composition of $n$ as a finite sequence
$C=(c_1,c_2,\ldots,c_k)$ of non-zero integers such that
$|c_1|+|c_2|+\cdots+|c_k|=n$; then the sequence
$C^+=(|c_1|,|c_2|,\ldots,|c_k|)$ is a composition of $n$.
Given such a sequence $C$, we observe that the Young subgroup $\mathfrak
S_{C^+}$ of $\mathfrak S_n$, acting on $G^n$, leaves stable the subgroup
$G^C=G^{c_1}\times G^{c_2}\times\cdots\times G^{c_k}$. We can thus make
the semidirect product $W_C=G^C\rtimes\mathfrak S_{C^+}$; this is a
subgroup of $G^n\rtimes\mathfrak S_{C^+}=G\wr\mathfrak S_{C^+}$, hence
of $W_n$. For instance, $W_{(-n)}=\mathfrak S_n=W_{-n}$.

Let $C$ be a signed composition of $n$. By Proposition~2.8 of
\cite{Bonnafe-Hohlweg04}, each left coset $wW_C$ of $W_n$ modulo $W_C$
contains a unique element of minimal length, called the distinguished
representative of this coset. Following \cite{Bonnafe-Hohlweg04}, we denote
the set of all these distinguished representatives by $X_C$ and we define
the element $\tilde x_C=\sum_{w\in X_C}w$ in the group ring $\mathbb ZW_n$.

The dual group $\Gamma$ of $G$ has also two elements, namely the trivial
character $t$ and the sign character~$s$. Since $\Gamma$ is canonically
isomorphic to $G$, the group ring $\mathbb Z\bigl[\Gamma\wr\mathfrak
S_n\bigr]=\mathscr F_n(\mathbb Z\Gamma)$ is canonically isomorphic to
$\mathbb ZW_n$. The elements $\tilde x_C$ can therefore be viewed as
elements in $\mathscr F_n(\mathbb Z\Gamma)$. To complete the notation,
we set
$$z_n=(\underbrace{s,s,\ldots,s}_{\text{$n$ times}})\cdot[n(n-1)\cdots 1]$$
for any $n\geq1$, where $[n(n-1)\cdots 1]$ is the longest permutation in
$\mathfrak S_n$, and we agree that $\tilde x_{(0)}$, $y_{0,s}$, $y_{0,t}$ and
$z_0$ are all equal to the unit of $\mathscr F_0(\mathbb Z\Gamma)=\mathbb Z$.

\begin{proposition}
\label{pr:EltsBonnHohl}
\begin{enumerate}
\item\label{it:PrEBHa}We have the following relations:
\begin{gather}
\tilde x_{(n)}=y_{n,t},\label{eq:PrEBHa}\\
\tilde x_{(-n)}=\sum_{i=0}^nz_i*y_{n-i,t},\label{eq:PrEBHb}\\[3pt]
\sum_{i=0}^n(-1)^i\,z_i*y_{n-i,s}=\sum_{i=0}^n(-1)^i\,y_{n-i,s}*z_i
=\begin{cases}1&\text{if $n=0$,}\\0&\text{if $n>0$,}\end{cases}
\label{eq:PrEBHc}\\[3pt]
\tilde x_C=\tilde x_{(c_1)}*\tilde x_{(c_2)}*\cdots*\tilde x_{(c_k)}.
\label{eq:PrEBHd}
\end{gather}
for any non-negative integer $n$ and any signed composition
$C=(c_1,c_2,\ldots,c_k)$ of $n$.
\item\label{it:PrEBHb}The elements $\tilde x_C$ form a basis of
$\MR_n(\mathbb Z\Gamma)$, where $C$ is a signed composition of $n$.
\end{enumerate}
\end{proposition}
\begin{proof}
\begin{enumerate}
\item Formula~(\ref{eq:PrEBHa}) holds because both members are equal
to the unit of the ring $\mathbb ZW_n$, by definition.

By Example~2.23 in \cite{Bonnafe-Hohlweg04}, we know that $X_{(-n)}$
is the set of all elements $w\in W_n$ of the form
$$w=\sigma\cdot(\underbrace{-1,-1,\ldots,-1}_{i\text{ times}},
\underbrace{1,1,\ldots,1}_{n-i\text{ times}}),$$
where $\sigma\in\mathfrak S_n$ is decreasing on the interval $[1,i]$
and increasing on the interval $[i+1,n]$. This entails Formula
(\ref{eq:PrEBHb}), since in the identification of $G$ with $\Gamma$,
the elements $1$ and $-1$ correspond to $t$ and $s$, respectively.

Let $n$ be a positive integer. The set of all compositions of $n$ is a
ranked poset when endowed with the refinement order $\preccurlyeq$; here
the rank function is the map which associates to a composition $\mathbf d$
its number of parts $l(\mathbf d)$. The equality
$$x_{\mathbf c}=\sum_{\substack{\mathbf d\models n\\[2pt]\mathbf d
\preccurlyeq\mathbf c}}\ \sum_{\substack{\sigma\in\mathfrak S_n
\\[2pt]D(\sigma)=\mathbf d}}\ \sigma,$$
valid for each composition $\mathbf c$ of $n$, entails by Möbius inversion
\begin{equation}
\sum_{\substack{\sigma\in\mathfrak S_n\\[2pt]D(\sigma)=\mathbf c}}\
\sigma=\sum_{\substack{\mathbf d\models n\\[2pt]\mathbf d
\preccurlyeq\mathbf c}}\ (-1)^{l(\mathbf c)-l(\mathbf d)}\ x_{\mathbf d}.
\label{eq:PfPrEBH}
\end{equation}
Taking $\mathbf c=(1,1,\ldots,1)$ ($n$ times) in Formula~(\ref{eq:PfPrEBH})
and multiplying by $(s,s,\ldots,s)$, we obtain
$$z_n=\sum_{\substack{\mathbf d\models n\\[2pt]\mathbf d=(d_1,d_2,\ldots,
d_l)}}(-1)^{n-l}\ y_{d_1,s}*y_{d_2,s}*\cdots*y_{d_l,s}.$$
From there, one deduces easily Formula~(\ref{eq:PrEBHc}).

Finally Formula~(\ref{eq:PrEBHd}) is Example~5.3 in \cite{Bonnafe-Hohlweg04}.
\item Formula~(\ref{eq:PrEBHc}) implies that each element $z_n$ belongs
to $\MR(\mathbb Z\Gamma)$. Using Formulas~(\ref{eq:PrEBHa}),
(\ref{eq:PrEBHb}) and (\ref{eq:PrEBHd}), we then deduce
that each element $\tilde x_C$ belongs to $\MR(\mathbb Z\Gamma)$, where
$C$ is a signed composition. In other words, the submodule $\MR'$ of
$\mathscr F(\mathbb Z\Gamma)$ spanned over $\mathbb Z$ by the elements
$\tilde x_C$ is contained in $\MR(\mathbb Z\Gamma)$. Formula
(\ref{eq:PrEBHd}) shows furthermore that $\MR'$ is a subalgebra for the
external product~$*$. Observing then that $\MR'$ contains all the
elements $y_{n,t}$ and $\tilde x_{(-n)}$, an easy induction based on
Formulas~(\ref{eq:PrEBHb}) and (\ref{eq:PrEBHc}) shows that each $z_n$
and each $y_{n,s}$ is in $\MR'$. This implies that $\MR'$ contains
$\MR(\mathbb Z\Gamma)$ because the latter is generated as an algebra by
the elements $y_{n,t}$ and $y_{n,s}$. It follows that the $\mathbb Z$-module
$\MR(\mathbb Z\Gamma)=\MR'$ is spanned by the elements $\tilde x_C$.

Now Proposition~\ref{pr:SigmaIsFreeAlg} (or more precisely, its
consequence stated at the end of Section~\ref{ss:SubbialgSigmaW}) implies
that $\MR_n(\mathbb Z\Gamma)$ is a free $\mathbb Z$-module whose rank $r$
is equal to the number of words $y_{c_1,\gamma_1}*y_{c_2,\gamma_2}*\cdots
*y_{c_k,\gamma_k}$, where $(c_1,c_2,\ldots,c_k)$ is a composition of $n$
and each $\gamma_i\in\{t,s\}$. Then any generating family of
$\MR_n(\mathbb Z\Gamma)$ with $r$ elements is a basis thereof. We conclude
that the family of elements $\tilde x_C$, where $C$ is a signed
composition of $n$, is a basis of $\MR_n(\mathbb Z\Gamma)$.
\end{enumerate}
\end{proof}

Bonnafé and Hohlweg call the submodule spanned by the elements
$\tilde x_C$ the `generalized descent algebra' and observe that it
coincides with the Mantaci-Reutenauer algebra $\MR_n(\mathbb Z\Gamma)$
(see \S3.1 in \cite{Bonnafe-Hohlweg04}). Assertion~\ref{it:PrEBHb} of
Proposition~\ref{pr:EltsBonnHohl} is roughly equivalent to this
observation, and indeed our proof follows closely the analysis in
\cite{Bonnafe-Hohlweg04}.

The associative algebra $\MR(\mathbb Z\Gamma)$ is freely generated
by the elements $y_{n,t}$ and $y_{n,s}$, where $n\geq1$. On the
other side, we have defined in Section~\ref{ss:CharacWreathProd} the
characters $\eta_n(t)$ and $\varepsilon_n(s)$. Thus there exists a
unique morphism of algebras $\tilde\theta:\MR(\mathbb Z\Gamma)\to\Rep(G)$
that maps $y_{n,t}$ and $y_{n,s}$ to $\eta_n(t)$ and $\varepsilon_n(s)$,
respectively.

\begin{proposition}
\label{pr:SoloBonnHohl}
\begin{enumerate}
\item\label{it:PrSBHa}The map $\tilde\theta$ enjoys all the properties
stated in Theorem~\ref{th:SoloHomo} for the map~$\theta_G$.
\item\label{it:PrSBHb}For any signed composition $C$ of a positive integer
$n$, $\tilde\theta(\tilde x_C)$ is the character of $W_n$ induced from the
trivial character of $W_C$.
\end{enumerate}
\end{proposition}
\begin{proof}
\begin{enumerate}
\item The graded bialgebra $\Lambda$ of symmetric functions has a
canonical involution $\omega$, which exchanges the complete symmetric
function $h_n$ with the elementary symmetric function $e_n$ of the same
degree (see I, (2.7) in \cite{Macdonald95}). Now $\Lambda(\Irr(G))$ is
the tensor product $\Lambda(t)\otimes\Lambda(s)$ of two copies of
$\Lambda$, so that $\id_\Lambda(t)\otimes\omega(s)$ is an involutive
automorphism of $\Lambda(\Irr(G))$. Equation~(\ref{eq:PfDescrEta}) shows
that the Frobenius characteristic $\ch$ maps $\varepsilon_n(s)$ to the
element $e_n(s)$ of $\Lambda(\Irr(G))$. Therefore the homomorphism
$\ch\circ\,\tilde\theta$ maps the two elements $y_{n,t}$ and $y_{n,s}$ to
$h_n(t)$ and $e_n(s)$, respectively, while $\ch\circ\,\theta_G$ maps
these elements to $h_n(t)$ and $h_n(s)$. Thus the diagram
$$\xymatrix{&\Rep(G)\ar[r]_{\simeq}^\ch&\Lambda(\Irr(G))\ar[dd]^{
\id_\Lambda(t)\otimes\omega(s)}\\\MR(\mathbb Z\Gamma)\ar[ur]^{\theta_G}
\ar[dr]_{\tilde\theta}&&\\&\Rep(G)\ar[r]_{\simeq}^\ch&\Lambda(\Irr(G))}$$
is commutative. Since $\ch$ and $\id_\Lambda(t)\otimes\omega(s)$ are
isomorphisms of graded bialgebras, $\tilde\theta$ inherits from
$\theta_G$ the properties stated in Assertion~\ref{it:ThSHa} of
Theorem~\ref{th:SoloHomo}. The proof of Assertions~\ref{it:ThSHb}
and~\ref{it:ThSHc} of Theorem~\ref{th:SoloHomo} presented in
Section~\ref{ss:SoloHomo} can be repeated with evident adjustments to
the case of $\tilde\theta$; the main difference lies in the proof of
the multiplicativity of $\tilde\theta$ with respect to the internal
product of $\MR_n(\mathbb Z\Gamma)$ and the ordinary product of $R(W_n)$,
where one must use the equalities $\eta_n(t)\varepsilon_n(s)
=\varepsilon_n(s)$ and $\varepsilon_n(s)\varepsilon_n(s)=\eta_n(t)$.
\item Let $u$ be an indeterminate. Applying $\tilde\theta$ to
Formulas~(\ref{eq:PrEBHb}) and (\ref{eq:PrEBHc}) and summing over $n$,
we find
\begin{gather*}
\left(\sum_{n\geq0}\tilde\theta(z_n)u^n\right)*
\left(\sum_{n\geq0}(-1)^n\,\tilde\theta(y_{n,s})u^n\right)=
\left(\sum_{n\geq0}(-1)^n\,\tilde\theta(y_{n,s})u^n\right)*
\left(\sum_{n\geq0}\tilde\theta(z_n)u^n\right)=1,\\
\sum_{n\geq0}\tilde\theta(\tilde x_{(-n)})u^n=
\left(\sum_{n\geq0}\tilde\theta(z_n)u^n\right)*
\left(\sum_{n\geq0}\tilde\theta(y_{n,t})u^n\right).
\end{gather*}
In Proposition~\ref{pr:DescrEta}, we have constructed an homomorphism $H$
from the additive group $R(G)$ into $\bigl(\Rep(G)[[u]]\bigr)^\times$
such that
\begin{align*}
H(t)&=\sum_{n\geq0}\eta_n(t)u^n=\sum_{n\geq0}\tilde\theta(y_{n,t})
u^n,\\[3pt]H(-s)&=\sum_{n\geq0}(-1)^n\varepsilon_n(s)u^n=
\sum_{n\geq0}(-1)^n\tilde\theta(y_{n,s})u^n.
\end{align*}
Then
$$\sum_{n\geq0}\tilde\theta(z_n)u^n=\left(\sum_{n\geq0}(-1)^n\tilde
\theta(y_{n,s})u^n\right)^{-1}=H(-s)^{-1}=H(s),$$
which implies in turn
$$\sum_{n\geq0}\tilde\theta(\tilde x_{(-n)})u^n=H(s)*H(t)=H(t+s)
=\sum_{n\geq0}\eta_n(t+s)u^n.$$
It follows that $\tilde\theta(\tilde x_{(-n)})=\eta_n(t+s)$. Now
the character $\eta_n(t+s)$ of $G\wr\mathfrak S_n$ is induced from the
trivial representation of $\mathfrak S_n$, because $t+s$ is the regular
character of $G$. Therefore $\tilde\theta(\tilde x_{(-n)})$ is the
character of $W_n$ induced from the trivial character of $W_{-n}$.

On the other side, $\tilde\theta(\tilde x_{(n)})=\tilde\theta(y_{n,t})=
\eta_n(t)$ is the trivial character of $W_n$. Using the transitivity
of induction, we thus find that for any signed composition
$C=(c_1,c_2,\ldots,c_k)$ of~$n$,
\begin{align*}
\tilde\theta(\tilde x_C)
&=\tilde\theta\bigl(\tilde x_{(c_1)}\bigr)*\tilde\theta\bigl(\tilde
x_{(c_2)}\bigr)*\cdots*\tilde\theta\bigl(\tilde x_{(c_k)}\bigr)\\[3pt]
&=\Ind_{G\wr\mathfrak S_{C^+}}^{G\wr\mathfrak S_n}\Bigl(
\Ind_{W_{c_1}}^{W_{|c_1|}}1\otimes\Ind_{W_{c_2}}^{W_{|c_2|}}1\otimes
\cdots\otimes\Ind_{W_{c_k}}^{W_{|c_k|}}1\Bigr)\\[3pt]
&=\Ind_{W_C}^{W_n}1,
\end{align*}
taking into account the identifications
$$G\wr\mathfrak S_n=W_n,\quad G\wr\mathfrak S_{C^+}\cong W_{|c_1|}
\times W_{|c_2|}\times\cdots\times W_{|c_k|}\quad\text{and}\quad
W_{c_1}\times W_{c_2}\times\cdots\times W_{c_k}\cong W_C.$$
This concludes the proof.
\end{enumerate}
\end{proof}

Assertion~\ref{it:PrSBHb} of this proposition says that our homomorphism
$\tilde\theta$ is equal to the homomorphism defined by Bonnafé and Hohlweg
\S3.1 in \cite{Bonnafe-Hohlweg04}. It follows then from the results of
these authors that $\tilde\theta$ extends Solomon's original homomorphism.

On the contrary, $\theta_G$ does not extend Solomon's original
homomorphism. Indeed we observe that the parabolic subgroups of the
Coxeter system $W_n$ are the subgroups $W_C$, where the signed
composition $C=(c_1,c_2,\ldots,c_k)$ has all its parts negative with
the possible exception of $c_1$. Therefore the original Solomon algebra
of $W_n$ is the submodule $\tilde\Sigma$ of the group ring $\mathbb ZW_n$
spanned by the elements $\tilde x_C$ for such signed compositions.
Taking $n=2$ and using Relations~(\ref{eq:PrEBHa})--(\ref{eq:PrEBHd}), we
check that the element $y_{2,s}=\tilde x_{(-1,-1)}-\tilde x_{(1,-1)}+
\tilde x_{(2)}-\tilde x_{(-2)}$ belongs to $\tilde\Sigma$. Its image
under $\ch\circ\,\theta_G$, namely $\ch(\eta_2(s))=h_2(s)$, is different
from its image under $\ch\circ\,\tilde\theta$, namely
$\omega(h_2(s))=e_2(s)$; it follows that $\theta_G\bigr|_{\tilde\Sigma}$
does not coincide with Solomon's original homomorphism
$\tilde\theta\bigr|_{\tilde\Sigma}$.

On the other side, our map $\theta_G$ shares with Solomon's
original homomorphism a property that Bonnafé and Hohlweg's extension
$\tilde\theta$ does not have, namely the symmetry property of
Theorem~\ref{th:SymmSoloHomo}. Again a counterexample can be
found already for $n=2$: one can indeed check that the value of the
character $\tilde\theta(\tilde x_{(-2)})$ on the element $\tilde
x_{(1,1)}$ is $6$, while the value of $\tilde\theta(\tilde x_{(1,1)})$
on $\tilde x_{(-2)}$ is $4$.

\section{Coloured combinatorial Hopf algebras}
\label{se:ColCombHopfAlg}
In the previous sections, we have presented our main constructions
and the applications which have motivated them. In the case where $G$
is the group with one element, the diagram~(\ref{eq:DiagrThetaVee})
is the usual diagram which relates the different kinds of symmetric
functions: ordinary, noncommutative, quasisymmetric, free quasisymmetric.
This diagram can be enriched with other bialgebras: the plactic and the
coplactic bialgebras \cite{Poirier-Reutenauer95}, the Loday-Ronco
bialgebra \cite{Loday-Ronco98}, the peak algebra \cite{Stembridge97},
etc.

Here we define analogues of some of the plactic and the coplactic
bialgebras and we insert them in~(\ref{eq:DiagrThetaVee}). The analogue
of the coplactic bialgebra presents two interests: first, the Solomon
homomorphism $\theta_G$ can be extended to it in a natural way; second
it is related to a construction already present in the litterature,
which we call the Robinson-Schensted-Okada correspondence.

\subsection{Categorical framework}
\label{ss:CatFmwk}
We start by setting up quickly a clean framework adapted to our goal.
We define a category $\mathscr E$ as follows. The objects of $\mathscr E$
are pairs $(B,?^*)$, where $B$ is a finite set and $?^*:b\mapsto b^*$
is an involutive map from $B$ to $B$. We generally omit the involution
in the notation, writing just $B$. Given two objects $B$ and $C$ of
$\mathscr E$, an homomorphism from $B$ to $C$ is the data of a
bijection $\varphi$ from a subset $B'$ of $B$ onto a subset $C'$ of $C$.

Any finite group $\Gamma$ can be considered as an object of $\mathscr E$,
where the involution $?^*$ is the map $\gamma\mapsto\gamma^{-1}$. The set
$\Irr(H)$ of irreducible characters of a finite group $H$ can also be
considered as an object of $\mathscr E$, with the complex conjugation
as involution $?^*$.

An object $(B,?^*)$ of $\mathscr E$ is viewed as the basis of the free
$\mathbb K$-module $\mathbb KB$. We define a pairing $\varpi$ on
$\mathbb KB$ by setting $\varpi(b,b')$ equal to $1$ if $b'=b^*$ and
equal to $0$ otherwise; this pairing is perfect and symmetric, for
$b\mapsto b^*$ is involutive. A morphism $f:B\to C$ induces a linear
map $\mathbb Kf:\mathbb KB\to\mathbb KC$ as follows: if $f$ is defined
by the bijection $\varphi:B'\to C'$, then $\mathbb Kf$ maps an element
$b$ of the basis $B$ to $\varphi(b)$ if $b\in B'$ and to $0$ otherwise.

Let $B$ be an object of $\mathscr E$. We can trace the constructions
of Section~\ref{se:GenConstr} at the level of bases. In more details,
the group $\mathfrak S_n$ acts by permutation on $B^n$. We denote the
cartesian product of $B^n$ with $\mathfrak S_n$ by $B\wr\mathfrak S_n$
and endow it with the following two-sided action of
$\mathfrak S_n$:
\begin{align*}
\pi\cdot(b_1,b_2,\ldots,b_n;\sigma)&=(b_{\pi^{-1}(1)},b_{\pi^{-1}(2)},
\ldots,b_{\pi^{-1}(n)});\pi\sigma)\\[3pt]
(b_1,b_2,\ldots,b_n;\sigma)\cdot\pi&=(b_1,b_2,\ldots,b_n;\sigma\pi).
\end{align*}
We can view $B\wr\mathfrak S_n$ as a basis of the free
$\mathbb K$-module $\mathscr F_n(\mathbb KB)$ by identifying the
element $(b_1,b_2,\ldots,b_n;\sigma)$ of $B\wr\mathfrak S_n$ with the
element $[(b_1\otimes b_2\otimes\cdots\otimes b_n)\#\sigma\bigr]$
of $\mathscr F_n(\mathbb KB)$.

We can then continue the construction and obtain from $B$ the free
quasisymmetric graded bialgebra $\mathscr F(\mathbb KB)$ and the
Novelli-Thibon algebra $\NT(\mathbb KB)$. Now the construction of the
Mantaci-Reutenauer bialgebra requires additionnally the data of a
coalgebra structure. But given a finite set $B$, one can always
define a structure of a coalgebra on $\mathbb KB$ by requiring that the
elements of $B$ are group-like; in other words, one agrees that the
coproduct $\delta$ and the counit $\varepsilon$ are defined by
$$\delta(b)=b\otimes b\quad\text{and}\quad\varepsilon(b)=1$$
for any $b\in B$. Endowing $\mathbb KB$ with this structure, we can
construct the Mantaci-Reutenauer bialgebra $\MR(\mathbb KB)$; to
translate into the notation the fact that this bialgebra depends on
the choice of the basis $B$ of $\mathbb KB$, we denote it by
$\mathscr D(B)$. By Proposition~\ref{pr:SigmaIsFreeAlg},
the associative algebra $\mathscr D(B)$ is freely generated by the
elements $y_{n,b}$ with $n\geq1$ and $b\in B$. The assignments
$B\rightsquigarrow\mathscr F(\mathbb KB)$ and $B\rightsquigarrow
\mathscr D(B)$ are covariant functors from the category $\mathscr E$
to the category of graded bialgebras.

Finally, given an object $B$ of $\mathscr E$, the perfect symmetric
pairing $\varpi$ on $\mathbb KB$ can be extended to a perfect
symmetric pairing $\varpi_{\mathrm{tot}}$ on $\mathscr F(\mathbb
KB)$ (see Section~\ref{ss:DualFunctF}). The basis $B\wr\mathfrak S_n$
is dual to itself with respect to $\varpi_{\mathrm{tot}}$; more
precisely, the basis element dual to $\alpha=(b_1,b_2,\ldots,b_n;\sigma)$
is $\alpha^*=\bigl[\sigma^{-1}\cdot(b_1^*,b_2^*,\ldots,b_n^*;e_n)\bigr]$.

\subsection{Coloured descent compositions}
\label{ss:ColDescComp}
Let $B$ be an object of $\mathscr E$. Since $\mathscr D(B)$ is a direct
summand of $\mathscr F(\mathbb KB)$ by Proposition~\ref{pr:MRIsDirSum},
the dual bialgebra $\mathscr D(B)^\vee$ is canonically isomorphic to
the quotient $\mathscr F(\mathbb KB)/\mathscr D(B)^\circ$. Our aim is
to study the subbialgebra $\mathscr D(B)$ and the quotient bialgebra
$\mathscr F(\mathbb KB)/\mathscr D(B)^\circ$ of $\mathscr F(\mathbb KB)$
on the level of basis in a combinatorial way.

We begin with definitions. A $B$-composition is a finite sequence
$\mathbf c=((c_1,b_1),(c_2,b_2),\ldots,\linebreak(c_k,b_k))$ of elements
of $\mathbb Z_{>0}\times B$. The size of $\mathbf c$ is the integer
$\|\mathbf c\|=c_1+c_2+\cdots+c_k$. The dual of $\mathbf c$ is
the $B$-composition $\mathbf c^*=((c_1,b_1^*),(c_2,b_2^*),\ldots,
(c_k,b_k^*))$. Given two $B$-compositions $\mathbf c=((c_1,b_1),
(c_2,b_2),\ldots,(c_k,b_k))$ and $\mathbf d=((d_1,b'_1),(d_2,b'_2),
\ldots,(d_l,b'_l))$ of the same size $n$, we say that $\mathbf c$ is
a refinement of $\mathbf d$ and we write $\mathbf c\succcurlyeq\mathbf d$
if there holds
\begin{align*}
(c_1,c_2,\ldots,c_k)&\succcurlyeq(d_1,d_2,\ldots,d_l),\\
(\underbrace{b_1,\ldots,b_1}_{c_1\text{ times}},
\underbrace{b_2,\ldots,b_2}_{c_2\text{ times}},\ldots,
\underbrace{b_k,\ldots,b_k}_{c_k\text{ times}})&=
(\underbrace{b'_1,\ldots,b'_1}_{d_1\text{ times}},
\underbrace{b'_2,\ldots,b'_2}_{d_2\text{ times}},\ldots,
\underbrace{b'_k,\ldots,b'_k}_{d_k\text{ times}}).
\end{align*}
The relation $\preccurlyeq$ is a partial order on the set of
$B$-compositions of $n$.

We associate to each element $\alpha\in B\wr\mathfrak S_n$ two
$B$-compositions $D(\alpha)$ and $R(\alpha)$ of size $n$. The `descent
composition' $D(\alpha)$ is constructed by the following procedure,
due to Mantaci and Reutenauer \cite{Mantaci-Reutenauer95}. We first
write $\alpha$ as $\bigl(b_{\sigma^{-1}(1)},b_{\sigma^{-1}(2)},\ldots,
b_{\sigma^{-1}(n)};\sigma\bigr)=\bigl[\sigma\cdot(b_1,b_2,\ldots,b_n;e_n)
\bigr]$, as before. Then one decomposes the interval $[1,n]$ into the
largest subintervals on which the map $i\mapsto b_i$ is constant, and
after that, one decomposes each such subinterval into the largest
subsubintervals on which the map $i\mapsto\sigma(i)$ is increasing.
Each subsubinterval yields a pair formed by its length and the value
taken by the map $i\mapsto b_i$. Then $D(\alpha)$ is the ordered list
of all these pairs. We define the `receding composition' of $\alpha$ by
the equality $R(\alpha)=D(\alpha^*)^*$. An example illustrates these
definitions. We take $a$ and $b$ in $B$, $n=7$ and $\alpha=(a,a,b,a,b,
b,a;1426735)=\bigl[1426735\cdot(a,a,a,b,a,b,b;e_7)\bigr]$; then
$D(\alpha)=((2,a),(1,a),(1,b),(1,a),(2,b))$ and
$R(\alpha)=((2,a),(1,b),(1,a),(1,b),(1,b),(1,a))$.

For any $i\in\{1,2,\ldots,n-1\}$, we denote by $s_i$ the transposition
in $\mathfrak S_n$ that exchanges $i$ and $i+1$. We say that two
elements $\alpha$ and $\alpha'$ of $B\wr\mathfrak S_n$ are related by
an Atkinson relation and we write $\alpha\underset A\sim\alpha'$ if
there exists an index $i$ such that:
\begin{itemize}
\item $\alpha'=\alpha\cdot s_i$;
\item writing $\alpha=(b_1,b_2,\ldots,b_n;\sigma)$, the map $j\mapsto
b_j$ is not constant on the interval $[\sigma(i),\sigma(i+1)]$ or
the inequality $|\sigma(i+1)-\sigma(i)|>1$ holds.
\end{itemize}
(In the case where $\sigma(i+1)<\sigma(i)$, the notation $[\sigma(i),
\sigma(i+1)]$ means the interval $[\sigma(i+1),\sigma(i)]$.)
The Atkinson relation is clearly symmetric.

The following proposition explains the relation between these combinatorial
definitions and the maps $\mathscr D(B)\hookrightarrow\mathscr
F(\mathbb KB)$ and $\mathscr F(\mathbb KB)\twoheadrightarrow
\mathscr F(\mathbb KB)/\mathscr D(B)^\circ$.
\begin{proposition}
\label{pr:CombMRBialg}
Let $B$ be an object of $\mathscr E$ and let $n$ be a non-negative integer.
\begin{enumerate}
\item\label{it:PrCMRBa}The submodule $\mathscr D(B)\cap\mathscr
F_n(\mathbb KB)$ of $\mathscr F_n(\mathbb KB)$ is spanned over
$\mathbb K$ by the elements
$$\sum_{\substack{\alpha\in B\wr\mathfrak S_n\\[2pt]D(\alpha)=\mathbf
c}}\alpha,$$
where $\mathbf c$ is a $B$-composition.
\item\label{it:PrCMRBb}Two elements $\alpha$ and $\alpha'$ in
$B\wr\mathfrak S_n$ have the same receding composition $R(\alpha)=
R(\alpha')$ if and only if there exists a sequence of elements
$\alpha_1$, $\alpha_2$, \dots, $\alpha_k$ such that
$$\alpha=\alpha_1\underset A\sim\alpha_2\underset A\sim\cdots
\underset A\sim\alpha_k=\alpha'.$$
\item\label{it:PrCMRBc}The module $\mathscr D(B)^\circ\cap
\mathscr F_n(\mathbb KB)$ is spanned over $\mathbb K$ by the set
$$\{\alpha-\alpha'\mid\text{$\alpha$ and $\alpha'$ in $B\wr\mathfrak
S_n$ with $\alpha\underset A\sim\alpha'$}\}.$$
\end{enumerate}
\end{proposition}
\begin{proof}
\begin{enumerate}
\item We observe that for any $B$-composition $\mathbf c=((c_1,b_1),
(c_2,b_2),\ldots,(c_k,b_k))$ of $n$,
\begin{align}
y_{c_1,b_1}*y_{c_2,b_2}*\cdots*y_{c_k,b_k}
&=x_{(c_1,c_2,\ldots,c_k)}\cdot(b_1^{\otimes c_1}\otimes
b_2^{\otimes c_1}\otimes\cdots\otimes b_k^{\otimes c_k}\#e_n)\notag\\[5pt]
&=\sum_{\substack{\sigma\in\mathfrak S_n\\[2pt]D(\sigma)\preccurlyeq
(c_1,c_2,\ldots,c_k)}}\sigma(b_1^{\otimes c_1}\otimes
b_2^{\otimes c_1}\otimes\cdots\otimes b_k^{\otimes c_k}\#e_n)\notag\\[5pt]
&=\sum_{\substack{\alpha\in B\wr\mathfrak S_n\\[2pt]D(\alpha)
\preccurlyeq\mathbf c}}\alpha.\label{eq:PfPrCMRB}
\end{align}
Assertion~\ref{it:PrCMRBa} follows easily from (\ref{eq:PfPrCMRB})
and from the fact that $\mathscr D(B)\cap\mathscr F_n(\mathbb KB)$
is spanned over $\mathbb K$ by such products $y_{c_1,b_1}*y_{c_2,b_2}
*\cdots*y_{c_k,b_k}$.
\item For any two elements $\alpha$ and $\alpha'$ in $B\wr\mathfrak S_n$,
the relation $\alpha^*\underset A\sim\alpha^{\prime*}$ is equivalent to
the existence of an index $i\in\{1,2,\ldots,n-1\}$ such that:
\begin{itemize}
\item$\alpha'=s_i\cdot\alpha$;
\item writing $\alpha=\bigl[\sigma\cdot(b_1,b_2,\ldots,b_n;e_n)\bigr]$,
the map $j\mapsto b_j$ is not constant on the interval
$[\sigma^{-1}(i),\sigma^{-1}(i+1)]$ or the inequality
$|\sigma^{-1}(i+1)-\sigma^{-1}(i)|>1$ holds.
\end{itemize}
An easy verification shows then that $D(\alpha)=D(\alpha')$ as soon as
$\alpha^*\underset A\sim\alpha^{\prime*}$, and therefore as soon as
$\alpha^*$ and $\alpha^{\prime*}$ are related by a sequence of
Atkinson relations.

Let now $\mathbf c=((c_1,b_1),(c_2,b_2),\ldots,(c_k,b_k))$ be a
$B$-composition of $n$ and set
$$(\tilde b_1,\tilde b_2,\ldots,\tilde b_n)=(\underbrace{b_1,b_1,\ldots,
b_1}_{\text{$c_1$ times}},\underbrace{b_2,b_2,\ldots,b_2}_{\text{$c_2$
times}},\ldots,\underbrace{b_k,b_k,\ldots,b_k}_{\text{$c_k$ times}}).$$
Each element $\alpha$ in $B\wr\mathfrak S_n$ such that $D(\alpha)=
\mathbf c$ can be written $\alpha=\bigl[\sigma\cdot(\tilde b_1,\tilde
b_2,\ldots,\tilde b_n;e_n)\bigr]$. We now apply successive Atkinson
relations to $\alpha^*$ to reduce as much as possible the number of
inversions of $\sigma$, obtaining eventually an element $\alpha_0^*$.
By the previous paragraph, the descent composition is preserved at
each step of the process, so that $D(\alpha_0)=\mathbf c$.

We now observe that $\alpha_0$ depends only on $\mathbf c$ and not on
the element $\alpha$ from which we started or on the choices made during
the reduction process. Indeed let us write $\alpha_0=\bigl[\sigma_0\cdot
(\tilde b_1,\tilde b_2,\ldots,\tilde b_n;e_n)\bigr]$. The equality
$D(\alpha_0)=\mathbf c$ holds, and there is no permutation
$\sigma'\in\mathfrak S_n$ with smaller number of inversions than
$\sigma_0$ such that $\alpha^*\underset A\sim\alpha^{\prime*}$, where
$\alpha'=\bigl[\sigma'\cdot(\tilde b_1,\tilde b_2,\ldots,\tilde
b_n;e_n)\bigr]$. Setting $t_j=c_1+c_2+\cdots+c_j$ for each $j\in\{1,2,
\ldots,k-1\}$, these constraints imply in turn the equivalence of
the three following assertions for each $i\in[1,n-1]$:
\begin{itemize}
\item there exists an index $j\in\{1,2,\ldots,k-1\}$ such that
$i=t_j$ and $b_j=b_{j+1}$;
\item$\sigma(i)>\sigma(i+1)$;
\item$\sigma(i)=\sigma(i+1)+1$.
\end{itemize}
The uniqueness of $\sigma_0$, hence of $\alpha_0$, can be easily
derived from this.

Summarizing, we have seen that for any two elements $\alpha$ and
$\alpha'$ in $B\wr\mathfrak S_n$:
\begin{itemize}
\item If $\alpha^*$ and $\alpha^{\prime*}$ are related by a sequence
of Atkinson relations, then $D(\alpha)=D(\alpha')$.
\item If $D(\alpha)=D(\alpha')$, then starting from $\alpha^*$ as well
as $\alpha^{\prime*}$, one may reach the same element $\alpha_0^*$ by
applying a sequence of Atkinson relations.
\end{itemize}
Therefore $\alpha^*$ and $\alpha^{\prime*}$ are related by a sequence
of Atkinson relations if and only if $D(\alpha)=D(\alpha')$. This fact
is equivalent to Assertion~\ref{it:PrCMRBb}.
\item Let $x=\sum_{\alpha\in B\wr\mathfrak S_n}a_\alpha\,\alpha$ be
an element of $\mathscr F_n(\mathbb KB)$, where $a_\alpha\in\mathbb K$.
Then for any $B$-composition $\mathbf c$ of size $n$,
$$\varpi_{\mathrm{tot}}\left(x,\ \sum_{\substack{\alpha\in B\wr\mathfrak
S_n\\[2pt]D(\alpha)=\mathbf c}}\alpha\right)=\sum_{\substack{\alpha\in
B\wr\mathfrak S_n\\[2pt]D(\alpha)=\mathbf c}}a_{\alpha^*}=\sum_{\substack
{\alpha\in B\wr\mathfrak S_n\\[2pt]R(\alpha)=\mathbf c^*}}a_{\alpha}.$$
The element $x$ is orthogonal to $\mathscr D(B)$ if and only if this
quantity vanishes for all $\mathbf c$. Assertion~\ref{it:PrCMRBc} is then
a direct consequence of Assertion~\ref{it:PrCMRBb}.
\end{enumerate}
\end{proof}

The result stated in Proposition~\ref{pr:CombMRBialg}~\ref{it:PrCMRBb}
above was first obtained by Atkinson (see \cite{Atkinson86}, Corollary
on~p.~352) for the case where $B$ has only one element.

We already mentioned in Section~\ref{ss:CatFmwk} that the
assignments $B\rightsquigarrow\mathscr F(\mathbb KB)$ and
$B\rightsquigarrow\mathscr D(B)$ are covariant functors from the
category $\mathscr E$ to the category of graded bialgebras. By
Proposi\-tion~\ref{pr:CombMRBialg}~\ref{it:PrCMRBc}, the biideal
$\mathscr D(B)^\circ$ of $\mathscr F(\mathbb KB)$ is functorial in $B$,
which implies that $B\rightsquigarrow\mathscr F(\mathbb KB)/\mathscr
D(B)^\circ$ is a covariant functor from $\mathscr E$ to the category of
graded bialgebras. One may also observe that $B\rightsquigarrow\mathscr
D(B)^\vee$ is a contravariant functor between the same categories, and
that the two graded bialgebras $\mathscr F(\mathbb KB)/\mathscr
D(B)^\circ$ and $\mathscr D(B)^\vee$ are isomorphic.

\subsection{Tableaux and the Robinson-Schensted-Knuth correspondence}
\label{ss:TabRSKCorr}
In this section, we recall some classical stuff to fix the notations
needed to present the Robinson-Schensted-Okada correspondence.

Let $\mathscr A$ be a totally ordered set (an alphabet). An
$\mathscr A$-weight is a finite multiset of $\mathscr A$, that is, a map
$\mu:\mathscr A\to\mathbb N$ with finite support. Thus for instance a
$\mathbb Z_{>0}$-weight is an infinite sequence $\mu=(\mu_1,\mu_2,\ldots)$
of non-negative integers, all of whose terms but a finite number vanish.
The size of a weight $\mu$ is the sum of its values; we denote it by
$|\mu|$. The weight of a word $w=a_1a_2\cdots a_n$ with letters in
$\mathscr A$ is the $\mathscr A$-weight $\mu$ such that any letter
$a\in\mathscr A$ occurs $\mu(a)$ times in $w$; we denote it by $\wt(w)$.

A semistandard tableau $T$ with entries in $\mathscr A$ is a Young
diagram whose boxes are labelled by letters in $\mathscr A$ in such a
way that the rows are weakly increasing from left to right and the
columns are strictly increasing from top to bottom. The shape of $T$ is
the partition $\lambda=(\lambda_1,\lambda_2,\ldots)$ such that $T$ has
$\lambda_1$ boxes in the first row, $\lambda_2$ boxes in the second row,
and so on; we denote it by $\sh(T)$. The weight of $T$ is the $\mathscr
A$-weight $\mu$ such that any letter $a\in\mathscr A$ occurs $\mu(a)$
times as the label in a box of $T$; we denote it by $\wt(T)$. A tableau
$T$ filled with positive integers is said to be standard if its weight is
$$(\underbrace{1,1,\ldots,1}_{\mysmash{|\sh(T)|\text{ times}}},
0,0,\ldots).$$

To a word $w=a_1a_2\cdots a_n$ with letters in $\mathscr A$, the
Robinson-Schensted correspondence associates a pair $(P,Q)$ of tableaux
with the same shape, such that $\wt(P)=\wt(w)$ and $Q$ is standard.
The insertion tableau $P$ is constructed inductively using the well-known
`bump' procedure; the label in a box of the record tableau $Q$ indicate
the number of the step at which this box appears during the making of $P$.

One says that two words $w=a_1a_2\cdots a_n$ and $w'=a'_1a'_2\cdots a'_n$
with letters in $\mathscr A$ and of the same length are related by a
Knuth relation and one writes $w\underset K\sim w'$ if one can find two
decompositions $w=xuy$ and $w'=xu'y$ of $w$ and $w'$ as the concatenation
of subwords in such a way that one of the two following conditions holds:
\begin{description}
\item[(a)]There exist three letters $a\leq b<c$ in $\mathscr A$ such
that $\{u,u'\}=\{acb,cab\}$.
\item[(b)]There exist three letters $a<b\leq c$ in $\mathscr A$ such
that $\{u,u'\}=\{bac,bca\}$.
\end{description}
The following results can be found in \cite{Knuth70}.

\begin{proposition}
\label{pr:RSKCorresp}
\begin{enumerate}
\item\label{it:PrRSKCa}Let $(P,Q)$ be the image of the word
$w=a_1a_2\cdots a_n$ under the Robinson-Schensted correspondence. Then
$a_i>a_{i+1}$ if and only if the box of $Q$ that contains the label
$i+1$ appears south or south-west to the box that contains the label $i$.
\item\label{it:PrRSKCb}Two words $w$ and $w'$ with letters in $\mathscr A$
have the same insertion tableau $P$ under the Robinson-Schensted
correspondence if and only if there exists a sequence of words $w_1$,
$w_2$, \dots, $w_k$ such that $$w=w_1\underset K\sim w_2\underset
K\sim\cdots\underset K\sim w_k=w'.$$
\end{enumerate}
\end{proposition}

Knuth has extended the scope of the Robinson-Schensted correspondence
to a slightly more general situation, which we recall now. Let $\mathscr
B$ be a second alphabet. Given an $\mathscr A$-weight $\mu$ and a
$\mathscr B$-weight $\nu$, we denote by $\mathscr M_{\mu,\nu}$ the set
of matrices $M=(m_{ab})_{(a,b)\in\mathscr A\times\mathscr B}$ with
non-negative integral entries and with row-sum $\mu$ and column-sum
$\nu$, that is,
$$\mu_a=\sum_{b\in\mathscr B}m_{ab}\quad\text{for all $a$ and}\quad
\nu_b=\sum_{a\in\mathscr A}m_{ab}\quad\text{for all $b$.}$$
(This condition tacitely implies that all but a finite number of entries
of $M$ vanish and that $\mu$ and $\nu$ have the same size. The notation
$\mathscr M_{\mathbf c,\mathbf d}$ used in Section~\ref{ss:DblCosetSymGp}
is a particular case of this one.)

We order the product $\mathscr B\times\mathscr A$ lexicographically.
An element $M\in\mathscr M_{\mu,\nu}$ can be seen as a finite multiset
of $\mathscr B\times\mathscr A$, whose elements can be listed in
increasing order: $((b_1,a_1),(b_2,a_2),\ldots,\linebreak(b_n,a_n))$.
In this way, $M$ determines two words $w_{\mathscr A}=a_1a_2\cdots
a_n$ and $w_{\mathscr B}=b_1b_2\cdots b_n$, with the obvious property
that $\mu=\wt(w_{\mathscr A})$ and $\nu=\wt(w_{\mathscr B})$. The
Robinson-Schensted correspondence applied to $w_{\mathscr A}$ yields
a pair of tableaux $(P,\tilde Q)$. Substituting in each box of
$\tilde Q$ the label $j$ by the letter $b_j$, we obtain a tableau $Q$.
With these notations, Knuth has shown in \cite{Knuth70} that the map
$T\mapsto(P,Q)$ is a bijection from $\mathscr M_{\mu,\nu}$ onto
\begin{equation*}
\left\{(P,Q)\Biggm|\begin{aligned}\text{$P$ and }&\text{$Q$ tableaux
with $\sh(P)=\sh(Q)$,}\\&\text{$\wt(P)=\mu$ and $\wt(Q)=\nu$}
\end{aligned}\right\}.
\end{equation*}
Furthermore the transposition of $M$ corresponds to the exchange of
$P$ and $Q$. It is usual to call this map the RSK correspondence.

\subsection{The Robinson-Schensted-Okada correspondence}
\label{ss:RSOCorr}
Let $B$ be an object of $\mathscr E$. We define a $B$-partition as
a family $\boldsymbol\lambda=(\lambda_b)_{b\in B}$ of partitions.
The size of $\boldsymbol\lambda$ is the integer $\|\boldsymbol
\lambda\|=\sum_{b\in B}|\lambda_b|$. The dual of $\boldsymbol\lambda$
is the $B$-partition $\boldsymbol\lambda^*=\bigl(b\mapsto\lambda_{b^*}
\bigr)$.

Let now $\mathscr A$ be an alphabet. We define a $B$-tableau with
entries in $\mathscr A$ as a family $\mathbf T=(T_b)_{b\in B}$ of
tableaux whose boxes are filled by elements of $\mathscr A$. The
shape of $\mathbf T$ is the $B$-partition $\sh(\mathbf T)=
(\sh(T_b))_{b\in B}$. A $B$-tableau $\mathbf T$ with entries in
$\mathbb Z_{>0}$ is said to be standard if
$$\sum_{b\in B}\wt(T_b)=(\underbrace{1,1,\ldots,1}_{\mysmash{\|
\sh(\mathbf T)\|\text{ times}}},0,0,\ldots).$$
In other words, all the labels $1$, $2$, \dots, $n$ are used once and
only once to fill the boxes of the tableaux $T_b$, where
$n=\|\sh(\mathbf T)\|$ is the total number of boxes in $\mathbf T$.

Now let $w=x_1x_2\cdots x_n$ be a word whose letters $x_i=(a_i,b_i)$
\vspace{2pt}belong to $\mathscr A\times B$. For each $b\in B$, we form a
matrix $M^{(b)}=\Bigl(m^{(b)}_{aj}\Bigr)_{(a,j)\in\mathscr A\times[1,n]}$
by \vspace{2pt}setting $m^{(b)}_{aj}$ equal to $1$ if $(a_j,b_j)=(a,b)$
and equal to $0$ otherwise. From the matrix $M^{(b)}$, the RSK
correspondence produces a pair of tableaux $(P_b,Q_b)$ with the same
shape. The family $\mathbf P=(P_b)_{b\in B}$ is a $B$-tableau with
entries in $\mathscr A$ such that $\sum_{b\in B}\wt(P_b)$ is the weight
of the word $a_1a_2\cdots a_n$; the family $\mathbf Q=(Q_b)_{b\in B}$
is a standard $B$-tableau; the tableaux $\mathbf P$ and $\mathbf Q$
have the same shape. We say that $\mathbf P$ and $\mathbf Q$ are the
insertion and record tableaux of $w$, respectively, and we call the map
$w\mapsto(\mathbf P,\mathbf Q)$ the RSO correspondence (for
Robinson-Schensted-Okada).

One can adapt the Knuth relations to the RSO correspondence in the
following way. We say that two words $w=x_1x_2\cdots x_n$ and
$w'=x'_1x'_2\cdots x'_n$ of the same length with letters $x_i=(a_i,b_i)$
and $x'_i=(a'_i,b'_i)$ are related by a Knuth relation and we write
$w\underset K\sim w'$ if there exists an index $i$ such that one of
the following two conditions holds:
\begin{description}
\item[(c)]$b_i\neq b_{i+1}$, $x_i=x'_{i+1}$, $x_{i+1}=x'_i$, and
$x_j=x'_j$ for all $j\not\in\{i,i+1\}$.
\item[(d)]$b_i=b_{i+1}=b_{i+2}=b'_i=b'_{i+1}=b'_{i+2}$, the two words
$u=a_ia_{i+1}a_{i+2}$ and $u'=a'_ia'_{i+1}a'_{i+2}$ are as in Condition
(a) or (b), and $x_j=x'_j$ for all $j\not\in\{i,i+1,i+2\}$.
\end{description}
Then we have the following analogue of Knuth's theorem.

\begin{proposition}
\label{pr:KnuthRSORel}
Two words $w$ and $w'$ with letters in $\mathscr A\times B$ have the same
insertion tableau $\mathbf P$ under the RSO correspondence if and only if
there exists a sequence of words $w_1$, $w_2$, \dots, $w_k$ such that
$$w=w_1\underset K\sim w_2\underset K\sim\cdots\underset K\sim w_k=w'.$$
\end{proposition}
\begin{proof}
Let $w=x_1x_2\cdots x_n$ be a word with the letters $x_i=(a_i,b_i)$, and
let $\mathbf P=(P_b)_{b\in B}$ be the insertion tableau of $w$. For each
$B\in B$, we form the word $w^{(b)}=a_{j_1}a_{j_2}\cdots a_{j_k}$, where
$(j_1,j_2,\ldots,j_k)$ is the list in increasing order of all indices $j$
for which $b_j=b$. By construction, $P_b$ is the insertion tableau in the
RSK image of the matrix $M^{(b)}$, so $P_b$ is the insertion tableau of
the word $w^{(b)}$. We fix an enumeration $b_1$, $b_2$, \dots, $b_l$ of
the elements of $B$ and we form the word $\overline w=w^{(b_1)}w^{(b_2)}
\cdots w^{(b_l)}$ by concatenation. Obviously $w$ and $\overline w$ are
related by a sequence of Knuth relations of type~(c).

Let now $w'$ be a word with the same length as $w$. We produce the words
$w^{\prime(b)}$ and $\overline{w'}=w^{\prime(b_1)}w^{\prime(b_2)}\cdots
w^{\prime(b_l)}$ in the same way as we formed $w^{(b)}$ and $\overline w$
from $w$. The words $w$ and $w'$ are related by a sequence of Knuth
relations of type~(c) or~(d) if and only if the words $\overline w$
and $\overline{w'}$ are related by a sequence of Knuth relations of
type~(d). By definition, this happens if and only if for each $b\in B$,
the words $w^{(b)}$ and $w^{\prime(b)}$ are related by a sequence of
Knuth relations as in Section~\ref{ss:TabRSKCorr}. On the other hand,
$w$ and $w'$ have the same insertion tableau $\mathbf P$ if and only if
for each $b\in B$, the words $w^{(b)}$ and $w^{\prime(b)}$ have the same
insertion tableau. The desired result now follows directly from
Proposition~\ref{pr:RSKCorresp}~\ref{it:PrRSKCb}.
\end{proof}

We now explain why we have added Okada's name after those of Robinson
and Schensted. Any element $\alpha\in B\wr\mathfrak S_n$ can be written
uniquely in the form $\alpha=\bigl[\sigma\cdot(b_1,b_2,\ldots,b_n;e_n)
\bigr]$, where $\sigma\in\mathfrak S_n$ and $(b_1,b_2,\ldots,b_n)\in
B^n$. It thus determines the word
$$w(\alpha)=(\sigma(1),b_1)(\sigma(2),b_2)\cdots(\sigma(n),b_n)$$
with letters in $[1,n]\times B$. We denote the RSO correspondent of
$w(\alpha)$ by $(\mathbf P(\alpha),\mathbf Q(\alpha))$. The element
$\alpha$ can be recovered from the data of $w(\alpha)$; it is therefore
characterized by $(\mathbf P(\alpha),\mathbf Q(\alpha))$. Finally, we
define the dual of a $B$-tableau $\mathbf T=(T_b)_{b\in B}$ as the
$B$-tableau $\mathbf T^*=\bigl(b\mapsto T_{b^*}\bigr)$, where
$b\mapsto b^*$ is the involution on $B$. The following result is in
substance a theorem of Okada~\cite{Okada90}.

\begin{proposition}
\label{pr:OkadaCorresp}
The map $\alpha\mapsto(\mathbf P(\alpha),\mathbf Q(\alpha))$ is a
bijection from $B\wr\mathfrak S_n$ onto the set of pairs of standard
$B$-tableaux with the same shape. For any element $\alpha$ of
$B\wr\mathfrak S_n$, there holds $\mathbf Q(\alpha^*)=\mathbf P(\alpha)^*$.
\end{proposition}

As an example, we consider the same situation as in Section
\ref{ss:ColDescComp}, that is, we take $a$, $b$ in~$B$, $n=7$ and
$\alpha=\bigl[1426735\cdot(a,a,a,b,a,b,b;e_7)\bigr]$. Then the
matrices $M^{(a)}$ and $M^{(b)}$ are
$$M^{(a)}=\begin{pmatrix}1&0&0&0&0&0&0\\0&0&1&0&0&0&0\\0&0&0&0&0&0&0\\
0&1&0&0&0&0&0\\0&0&0&0&0&0&0\\0&0&0&0&0&0&0\\0&0&0&0&1&0&0\end{pmatrix}
\qquad\text{and}\qquad M^{(b)}=\begin{pmatrix}0&0&0&0&0&0&0\\
0&0&0&0&0&0&0\\0&0&0&0&0&1&0\\0&0&0&0&0&0&0\\0&0&0&0&0&0&1\\
0&0&0&1&0&0&0\\0&0&0&0&0&0&0\end{pmatrix},$$
and we find
$$P_a=\young(127,4)\ ,\quad Q_a=\young(125,3)\ ,\quad P_b=\young(35,6)\
,\quad Q_b=\young(47,6)\ .$$

\noindent We will write $\alpha\underset K\sim\alpha'$ whenever the
words $w(\alpha)$ and $w(\alpha')$ are related by a Knuth relation.
Writing $\alpha=(b_1,b_2,\ldots,b_n;\sigma)$, one checks easily that
$\alpha\underset K\sim\alpha'$ if and only if $\alpha'=\alpha\cdot s_i$
for an index $i$ such that at least one of the following three
conditions is satisfied:
\begin{itemize}
\item$b_{\sigma(i)}\neq b_{\sigma(i+1)}$;
\item$\sigma(i-1)\in[\sigma(i),\sigma(i+1)]$ and
$b_{\sigma(i-1)}=b_{\sigma(i)}=b_{\sigma(i+1)}$;
\item$\sigma(i+2)\in[\sigma(i),\sigma(i+1)]$ and
$b_{\sigma(i+2)}=b_{\sigma(i)}=b_{\sigma(i+1)}$.
\end{itemize}
(Here again the notation $[\sigma(i),\sigma(i+1)]$ means the interval
$[\sigma(i+1),\sigma(i)]$ if ever $\sigma(i+1)<\sigma(i)$.) It follows
then from Proposition~\ref{pr:KnuthRSORel} that the insertion tableaux
$\mathbf P(\alpha)$ and $\mathbf P(\alpha')$ of two elements $\alpha$
and $\alpha'$ of $B\wr\mathfrak S_n$ are equal if and only if there
exists a sequence $\alpha_1$, $\alpha_2$, \dots, $\alpha_k$ such that
$$\alpha=\alpha_1\underset K\sim\alpha_2\underset K\sim\cdots
\underset K\sim\alpha_k=\alpha'.$$

\subsection{The plactic and the coplactic bialgebras}
\label{ss:PlacCoplacBialg}
In this section, we fix an object $B$ of the category $\mathscr E$
and we use the Robinson-Schensted-Okada correspondence to define a
subbialgebra and a quotient bialgebra of $\mathscr F(\mathbb KB)$,
called respectively the coplactic and the plactic bialgebra.

Given a standard $B$-tableau $\mathbf T=(T_b)_{b\in B}$, we define
an element $t_{\mathbf T}$ of $\mathscr F(\mathbb KB)$ by setting
$$t_{\mathbf T}=\sum_{\substack{\alpha\in B\wr\mathfrak S_n\\[2pt]
\mathbf Q(\alpha)=\mathbf T}}\alpha,$$
where $n=\|\sh(\mathbf T)\|$ is the total number of boxes in
$\mathbf T$. Clearly, the elements $t_{\mathbf T}$ are linearly
independent and the $\mathbb K$-submodule $\mathscr Q(B)$ that they
span is a direct summand of $\mathscr F(\mathbb KB)$. This submodule
depends on $B$ and not only on $\mathbb KB$. The following result is
the analogue of Proposition~\ref{pr:CombMRBialg}~\ref{it:PrCMRBc}.

\begin{proposition}
\label{pr:CoplacIsSubbialg}
Let $B$ be an object of $\mathscr E$ and let $n$ be a non-negative integer.
\begin{enumerate}
\item\label{it:PrCISa}The module $\mathscr Q(B)^\circ\cap\mathscr
F_n(\mathbb KB)$ is spanned over $\mathbb K$ by the set
$$\{\alpha-\alpha'\mid\text{$\alpha$ and $\alpha'$ in $B\wr\mathfrak
S_n$ with $\alpha\underset K\sim\alpha'$}\}.$$
\item\label{it:PrCISb}The submodules $\mathscr Q(B)$ and $\mathscr Q
(B)^\circ$ are respectively a graded subbialgebra and a graded
biideal of the graded bialgebra $(\mathscr F(\mathbb KB),*,\Delta)$.
\item\label{it:PrCISc}The submodule $\mathscr Q(B)\cap\mathscr Q(B)^\circ$
is spanned over $\mathbb K$ by the set
$$\{t_{\mathbf T}-t_{\mathbf T'}\mid\text{$\mathbf T$ and $\mathbf T'$
standard $B$-tableaux with $\sh(\mathbf T)=\sh(\mathbf T')$}\}.$$
\item\label{it:PrCISd}The submodule $\mathscr Q(B)+\mathscr Q(B)^\circ$
is a direct summand of $\mathscr F(\mathbb KB)$.
\end{enumerate}
\end{proposition}
\begin{proof}
Let $x=\sum_{\alpha\in B\wr\mathfrak S_n}a_\alpha\alpha$ be an element
of $\mathscr F_n(\mathbb KB)$, where each $a_\alpha\in\mathbb K$. Then
for any standard $B$-tableau $\mathbf T$,
$$\varpi_{\mathrm{tot}}(x,t_\mathbf T)=\sum_{\substack{\alpha\in
B\wr\mathfrak S_n\\[2pt]\mathbf Q(\alpha^*)=\mathbf T}}a_\alpha=
\sum_{\substack{\alpha\in B\wr\mathfrak S_n\\[2pt]\mathbf
P(\alpha)=\mathbf T^*}}a_\alpha.$$
The element $x$ is orthogonal to $\mathscr Q(B)$ if and only if
these quantities vanish for all $\mathbf T$. Assertion~\ref{it:PrCISa}
now follows from Proposition~\ref{pr:KnuthRSORel}, or more precisely,
from its consequence explained at the end of Section~\ref{ss:RSOCorr}.

Now let $\alpha=(b_1,b_2,\ldots,b_n;\sigma)$ and $\alpha'=(b'_1,b'_2,
\ldots,b'_n;\sigma')$ be two elements in $B\wr\mathfrak S_n$ that are
related by a Knuth relation. Then for each element $\alpha''=(b''_1,
b''_2,\ldots,b''_{n'};\sigma'')$ in $B\wr\mathfrak S_{n'}$ and each
permutation $\rho\in X_{n,n'}$, the two elements
$$\rho\cdot(b_1,b_2,\ldots,b_n,b''_1,b''_2,\ldots,b''_{n'};
\sigma\times\sigma'')\quad\text{and}\quad
\rho\cdot(b'_1,b'_2,\ldots,b'_n,b''_1,b''_2,\ldots,b''_{n'};
\sigma'\times\sigma'')$$
of $B\wr\mathfrak S_{n+n'}$ are related by a Knuth relation, because
$\rho$ is increasing on the interval $[1,n]$. Therefore
$(\alpha-\alpha')*\alpha''$, which is equal to the sum
$$\sum_{\rho\in X_{n,n'}}\bigl[\rho\cdot(b_1,b_2,\ldots,b_n,b''_1,
b''_2,\ldots,b''_{n'};\sigma\times\sigma'')-\rho\cdot(b'_1,b'_2,
\ldots,b'_n,b''_1,b''_2,\ldots,b''_{n'};\sigma'\times\sigma'')\bigr],$$
belongs to $\mathscr Q(B)^\circ$. Since $\mathscr Q(B)^\circ$
is spanned by such differences $\alpha-\alpha'$, we conclude that
$\mathscr Q(B)^\circ$ is a left ideal of $\mathscr F(\mathbb KB)$.
A similar reasoning shows that $\mathscr Q(B)^\circ$ is a right ideal.

Consider again an element $\alpha=(b_1,b_2,\ldots,b_n;\sigma)$ in
$B\wr\mathfrak S_n$. Given an integer $n'\in[0,n]$, we denote the
standardizations of the words $\sigma^{-1}(1)\;\sigma^{-1}(2)\;\cdots\;
\sigma^{-1}(n')$ and $\sigma^{-1}(n'+1)\;\sigma^{-1}(n'+2)\;\cdots\;
\sigma^{-1}(n)$ by $\pi_{n'}\in\mathfrak S_{n'}$ and $\pi'_{n-n'}\in
\mathfrak S_{n-n'}$, respectively. A straightforward but tedious
verification shows that whenever $\alpha$ undergoes a Knuth relation,
either both of
$$(b_1,b_2,\ldots,b_{n'};\pi_{n'})\quad\text{and}\quad
(b_{n'+1},b_{n'+2},\ldots,b_n;\pi'_{n-n'})$$
are left unchanged, or one of them remains the same and the other
undergoes a Knuth relation. This fact implies that the class
modulo $\mathscr Q(B)^\circ\otimes\mathscr F(\mathbb KB)+\mathscr
F(\mathbb KB)\otimes\mathscr Q(B)^\circ$ of
$$\Delta(\alpha)=\sum_{n'=0}^n(b_1,b_2,\ldots,b_{n'};\pi_{n'})
\otimes(b_{n'+1},b_{n'+2},\ldots,b_n;\pi'_{n-n'})$$
does not change when $\alpha$ undergoes a Knuth relation. We conclude
that
\begin{equation}
\Delta\bigl(\mathscr Q(B)^\circ\bigr)\subseteq\mathscr Q(B)^\circ
\otimes\mathscr F(\mathbb KB)+\mathscr F(\mathbb KB)\otimes
\mathscr Q(B)^\circ.
\label{eq:PolCoplacIsBiideal}
\end{equation}
Observing then that all homogeneous elements of $\mathscr Q(B)^\circ$
have positive degree, we see that the counit of $\mathscr F(\mathbb KB)$
vanishes on $\mathscr Q(B)^\circ$. Jointly with Equation
(\ref{eq:PolCoplacIsBiideal}), this means that $\mathscr Q(B)^\circ$
is a coideal of $\mathscr F(\mathbb KB)$.

We have therefore proved that $\mathscr Q(B)^\circ$ is a graded
biideal of $\mathscr F(\mathbb KB)$. Since $\mathscr Q(B)$ is a direct
summand of $\mathscr F(\mathbb KB)$, this is equivalent to the fact
that $\mathscr Q(B)$ is a subbialgebra of $\mathscr F(\mathbb KB)$,
which concludes the proof of Assertion~\ref{it:PrCISb}.

Proposition~\ref{pr:OkadaCorresp} implies that for each positive
integer $n$ and each pair $(\mathbf T,\mathbf T')$ of standard
$B$-tableaux with $n$ boxes,
\begin{align*}
\varpi_{\mathrm{tot}}(t_{\mathbf T},t_{\mathbf T'})
&=\bigl|\{\alpha\in B\wr\mathfrak S_n\mid\mathbf Q(\alpha)=\mathbf T,\
\mathbf Q(\alpha^*)=\mathbf T'\}\bigr|\\[3pt]
&=\bigl|\{\alpha\in B\wr\mathfrak S_n\mid\mathbf Q(\alpha)=\mathbf T,\
\mathbf P(\alpha)=\mathbf T^{\prime*}\}\bigr|\\[3pt]
&=\begin{cases}1&\text{if $\mathbf T$ and $\mathbf T'$ have the same
shape,}\\0&\text{otherwise.}\end{cases}
\end{align*}
Assertion~\ref{it:PrCISc} follows easily from this fact.

Given an finite index set $I$, we denote by $\Mat_I(\mathbb K)$ the set
of matrices with lines and columns indexed by $I$ and with entries in
$\mathbb K$. The subspace $\Mat_I^r(\mathbb K)$ of matrices
$(m_{ij})_{(i,j)\in I^2}$ such that all row sums $\sum_{j\in I}m_{ij}$
are equal is a direct summand of $\Mat_I(\mathbb K)$.\vspace{2pt}

For each $B$-partition $\boldsymbol\lambda=(\lambda_b)_{b\in B}$,
let $\mathscr T_{\boldsymbol\lambda}$ be the set of all standard
$B$-tableaux of shape $\boldsymbol\lambda$. Let $n$ be a positive
integer, and let $\Part_B(n)$ be the set of $B$-partitions of size $n$.
Using the RSO correspondence, we define a linear bijection between
\vspace{2pt}$\mathscr F_n (\mathbb KB)$ and $\prod_{\boldsymbol\lambda
\in\Part_B(n)}\Mat_{\mathscr T_{\boldsymbol\lambda}}(\mathbb K)$\linebreak
as follows: an element $\sum_{\alpha\in B\wr\mathfrak S_n}a_\alpha\alpha$
of $\mathscr F_n(\mathbb KB)$ \vspace{2pt}corresponds to a family of
matrices\linebreak$(M_{\boldsymbol\lambda})_{\boldsymbol\lambda\in
\Part_B(n)}$ if and only if for each $\alpha\in B\wr\mathfrak S_n$, the
coefficient $a_\alpha$ is equal to the entry in $M_{\boldsymbol\lambda}$
with row index $\mathbf P(\alpha)$ and column index $\mathbf Q(\alpha)$,
where $\boldsymbol\lambda=\sh(\mathbf P(\alpha))$. One checks without
much difficulty that $\bigl(\mathscr Q(B)+\mathscr Q(B)^\circ\bigr)
\cap\mathscr F_n(\mathbb KB)$ is mapped by this bijection to the
product $\prod_{\boldsymbol\lambda\in\Part_B(n)}\Mat_{\mathscr
T_{\boldsymbol\lambda}}^r(\mathbb K)$. Assertion \ref{it:PrCISd} follows.
\end{proof}

The subbialgebra $\mathscr Q(B)$ is called the coplactic bialgebra. We
denote the quotient\linebreak$\mathscr F(\mathbb KB)/\mathscr Q(B)^\circ$
by $\mathscr P(B)$ and we name it the plactic bialgebra. Both
assignments $B\rightsquigarrow\mathscr Q(B)$ and $B\rightsquigarrow
\mathscr P(B)$ are covariant functors from $\mathscr E$ to the category
of graded bialgebras, and moreover the graded bialgebras $\mathscr P(B)$
and $\mathscr Q(B)^\vee$ are isomorphic for each $B$.

\begin{other}{Remark}
It turns out that $\mathscr Q(B)$ is neither a left nor a right
internal $\mathscr D(B)$-submodule of $\mathscr F(\mathbb KB)$.
Indeed there is the following counterexample in degree $4$.
The set $B$ does not play any role here; we take it reduced to one
element and abbreviate $\mathscr D(B)$ and $\mathscr Q(B)$ to
$\mathscr D$ and $\mathscr Q$, respectively. We consider the
standard tableau $T=\;${\scriptsize$\young(13,24)$} and the elements
$$t_T=3142+2143\quad\text{and}\quad y=x_{(1,2,1)}-x_{(3,1)}-x_{(1,3)}+
x_{(4)}=3142+2143+4132+4231+3241.$$
Then $t_T$ belongs to $\mathscr Q$ and $y$ belongs to $\mathscr D$.
A direct computation yields
$$y\cdot t_T=4321+4231+1324+1234+3421+3412+4312+1423+2413+2314.$$
We observe that the permutation $3421$ appears with a positive
coefficient in $y\cdot t_T$, which is not the case of the permutation
$1432$, although they have the same record tableau. Therefore
$y\cdot t_T$ does not belong to $\mathscr Q$. One checks similarly
that $t_T\cdot y$ does not either belong to $\mathscr Q$.
\end{other}

\subsection{An homomorphism onto a bialgebra of coloured symmetric functions}
\label{ss:HomToColSymmFcts}
Our aim now is to extend the work of Poirier and Reutenauer
\cite{Poirier-Reutenauer95} to the present framework. We compare the
coplactic bialgebra $\mathscr Q(B)$ and the plactic bialgebra
$\mathscr P(B)$ with the Mantaci-Reutenauer algebra $\mathscr D(B)$
and its dual $\mathscr D(B)^\vee\cong\mathscr F(\mathbb KB)/\mathscr
D(B)^\circ$, and we insert them in a commutative diagram similar to
(\ref{eq:DiagrThetaVee}).

We need some preparation, and to begin with, we define the descent
composition $D(\mathbf T)$ of a standard $B$-tableau $\mathbf
T=(T_b)_{b\in B}$ in the following way. Let $n$ be the total number
of boxes in $\mathbf T$ and let $\beta:\{1,2,\ldots,n\}\to B$ the map
which sends a label $i$ to the element $b$ such that $i$ appears in a
box of $T_b$. We decompose the interval $[1,n]$ into the largest
subintervals on which the map $\beta$ takes a constant value; in turn
we decompose each subinterval into the largest possible subsubintervals,
so that for any two numbers $i$ and $j$ located in the same subsubinterval,
$i<j$ if and only if $i$ is located west or south-west to $j$. For each
subsubinterval, we form the pair consisting of its length $c$ and the
value $b$ taken on it by the map $\beta$. The ordered list of all these
pairs is the $B$-composition $D(\mathbf T)$. For instance, with $B=\{a,b\}$,
the descent composition of the $B$-tableau $\mathbf T$ given by
$$T_a=\young(125,3)\quad\text{and}\quad T_b=\young(47,6)$$
is $D(\mathbf T)=((2,a),(1,a),(1,b),(1,a),(2,b))$.

\begin{lemma}
\label{le:DescCompoStdTab}
\begin{enumerate}
\item\label{it:LeDCSTa}The descent composition of an element $\alpha
\in B\wr\mathfrak S_n$ coincides with the descent composition of its
record tableau $\mathbf Q(\alpha)$
\item\label{it:LeDCSTb}Let $\boldsymbol\lambda=(\lambda_b)_{b\in B}$ be
a $B$-partition and $\mathbf c=((c_1,b_1),(c_2,b_2),\ldots,(c_k,b_k))$
be a $B$-composi\-tion, both of the same size. For each $b\in B$, we
define a $\mathbb Z_{>0}$-weight $\mu^{(b)}=(\mu_1^{(b)},\mu_2^{(b)},
\ldots,\mu_k^{(b)},\linebreak0,0,\ldots)$ by setting $\mu_j^{(b)}=c_j$
if $b_j=b$ and $\mu_j^{(b)}=0$ otherwise. Then the two sets
\begin{equation*}
\left\{\mathbf T\Biggm|\begin{aligned}\text{$\mathbf T$ standard }&
\text{$B$-tableau with}\\\text{$\sh(\mathbf T)=\boldsymbol\lambda$ }&
\text{and $D(\mathbf T)\preccurlyeq\mathbf c$}\end{aligned}\right\}
\ \text{and}\ \,
\left\{\mathbf U\Biggm|\begin{aligned}\text{$\mathbf U=(U_b)_{b\in B}$
\thinspace$B$-tableau }&\text{with entries in $\mathbb Z_{>0}$}\\
\text{such that $\sh(\mathbf U)=\boldsymbol\lambda$ and}&
\text{ $\ \forall b,\ \wt(U_b)=\mu^{(b)}$}\end{aligned}\right\}
\end{equation*}
are equipotent.
\end{enumerate}
\end{lemma}
\begin{proof}
Assertion~\ref{it:LeDCSTa} is a direct consequence of Proposition
\ref{pr:RSKCorresp}~\ref{it:PrRSKCa}. Let us prove Assertion
\ref{it:LeDCSTb}. We set $t_i=c_1+c_2+\cdots+c_i$; we denote the first
set by $X$ and the second set by $Y$. Our aim is to construct mutually
inverse bijections from $X$ onto $Y$ and from $Y$ onto $X$.

Let first $\mathbf T=(T_b)_{b\in B}$ be an element of $X$. For each $b$,
we construct a tableau $U_b$ by substituting in each box of $T_b$ the
label $j$ it contains by the index $i$ such that $j\in[t_{i-1}+1,t_i]$.
Since $D(\mathbf T)\preccurlyeq\mathbf c$, each index $i$ appears
$c_i$ times in $U_{b_i}$ and does not appear in the other tableaux
$U_b$. Therefore each tableau $U_b$ has $\mu^{(b)}$ for weight. It
follows that the $B$-tableau $\mathbf U=(U_b)_{b\in B}$, which has
visibly the same shape as $\mathbf T$, namely $\boldsymbol\lambda$,
belongs to $Y$.

In the other direction, let $\mathbf U=(U_b)_{b\in B}$ be an element of
$Y$. By definition, any label $i\in\{1,2,\ldots,k\}$ appears $c_i$ times
in $U_{b_i}$. We replace these entries $i$ in the boxes of $U_{b_i}$
by the numbers $t_{i-1}+1$, $t_{i-1}+2$, \dots, $t_i$, proceeding in
increasing order whilst going south-west to north-east. These
substitutions transform the $B$-tableau $\mathbf U$ in a standard
$B$-tableau $\mathbf T$. By construction, $\mathbf T$ has the same shape
as $\mathbf U$, namely $\boldsymbol\lambda$, and satisfies $D(\mathbf T)
\preccurlyeq\mathbf c$; it thus belongs to $X$.

Routine verifications show that these correspondences are inverse
bijections, which entails Assertion~\ref{it:LeDCSTb}.
\end{proof}

\begin{corollary}
The inclusion $\mathscr D(B)\subseteq\mathscr Q(B)$ holds.
\end{corollary}
\begin{proof}
By Proposition~\ref{pr:CombMRBialg}~\ref{it:PrCMRBa}, the module
$\mathscr D(B)$ is spanned by elements of the form
$$\sum_{\substack{\alpha\in B\wr\mathfrak S_n\\[2pt]D(\alpha)
=\mathbf c}}\alpha,$$
where $n$ is a positive integer and $\mathbf c$ is a $B$-composition
of $n$.
By Lemma~\ref{le:DescCompoStdTab}~\ref{it:LeDCSTa}, such a sum may be
rewritten as
\begin{equation}
\sum_{\substack{\alpha\in B\wr\mathfrak S_n\\[2pt]D(\alpha)
=\mathbf c}}\alpha=
\sum_{\substack{\mathbf T\text{ standard $B$-tableau}\\[2pt]
D(\mathbf T)=\mathbf c}}\ \sum_{\substack{\alpha\in B\wr\mathfrak
S_n\\[2pt]\mathbf Q(\alpha)=\mathbf T}}\alpha
=\sum_{\substack{\mathbf T\text{ standard $B$-tableau}\\[2pt]
D(\mathbf T)=\mathbf c}}t_{\mathbf T}.
\label{eq:DescClsAsCoplacCls}
\end{equation}
It belongs therefore to $\mathscr Q(B)$. The corollary follows.
\end{proof}

Changing slightly the notation used in Section~\ref{ss:CharacWreathProd},
we use now the symbol $\Lambda$ to denote the algebra of symmetric
functions with coefficients in $\mathbb K$. It is indeed a bialgebra
(see I, 5, Ex.~25 in \cite{Macdonald95}). We keep the notation $h_n$
and $s_\lambda$ to denote the complete symmetric functions and the
Schur functions, where $n$ is a positive integer and $\lambda$ is a
partition. We consider a family $(\Lambda(b))_{b\in B}$ of copies of
$\Lambda$: given $b\in B$, we denote by $P(b)$ the image in $\Lambda(b)$
of an element $P\in\Lambda$. We carry out the tensor product
$\Lambda(B)=\bigotimes_{b\in B}\Lambda(b)$. Given a $B$-partition
$\boldsymbol\lambda=(\lambda_b)_{b\in B}$, we set
$\mathbf s_{\boldsymbol\lambda}=\prod_{b\in B}s_{\lambda_b}(b)$;
these elements $\mathbf s_{\boldsymbol\lambda}$ form a basis of the
$\mathbb K$-module $\Lambda(B)$. The pairing $\langle?,?\rangle$ on
$\Lambda(B)$ defined on this basis by
$$\langle\mathbf s_{\boldsymbol\lambda},\mathbf s_{\boldsymbol\lambda'}
\rangle=\begin{cases}1&\text{if $\boldsymbol\lambda'=\boldsymbol
\lambda^*$,}\\0&\text{otherwise,}\end{cases}$$
is then a perfect and symmetric pairing.

Let $\Theta_B:\mathscr Q(B)\to\Lambda(B)$ be the $\mathbb K$-linear map
such that $\Theta_B(t_\mathbf T)=\mathbf s_{\sh(\mathbf T)}$, for each
standard $B$-tableau $\mathbf T$. The following lemma will help us to
understand the behaviour of $\Theta_B$ on the subspace $\mathscr D(B)$
of $\mathscr Q(B)$.

\begin{lemma}
\label{le:ThetaCoplac}
For any $B$-composition $\mathbf c=((c_1,b_1),(c_2,b_2),\ldots,
(c_k,b_k))$, there holds
$$\Theta_B(y_{c_1,b_1}*y_{c_2,b_2}*\cdots*y_{c_k,b_k})=h_{c_1}(b_1)
h_{c_2}(b_2)\cdots h_{c_k}(b_k).$$
\end{lemma}
\begin{proof}
It is known (see I, (6.4) in \cite{Macdonald95} for a proof) that in the
ring $\Lambda$ of symmetric functions,
\begin{equation}
h_{\mu_1}h_{\mu_2}\cdots=\sum_{\lambda\text{ partition}}\left|\left\{
U\Biggm|\begin{aligned}\text{$U$}&\text{ tableau with entries in
$\mathbb Z_{>0}$}\\\text{such }&\text{that $\sh(U)=\lambda$ and
$\wt(U)=\mu$}\end{aligned}\right\}\right|\;s_\lambda
\label{eq:TransHS}
\end{equation}
for any $\mathbb Z_{>0}$-weight $\mu=(\mu_1,\mu_2,\ldots)$.

We fix a $B$-composition $\mathbf c=((c_1,b_1),(c_2,b_2),\ldots,(c_k,b_k))$
as in the statement of the lemma, of size say $n$, and we construct
a family $(\mu^{(b)})_{b\in B}$ of $\mathbb Z_{>0}$-weights as in
Lemma~\ref{le:DescCompoStdTab}~\ref{it:LeDCSTb}. Regrouping the factors
in the product $h_{c_1}(b_1)h_{c_2}(b_2)\cdots h_{c_k}(b_k)$ that
correspond to the different indices $b$ and applying
Formula~(\ref{eq:TransHS}), we find
\begin{multline*}
h_{c_1}(b_1)h_{c_2}(b_2)\cdots h_{c_k}(b_k)\\=\sum_{\boldsymbol\lambda
\text{ $B$-partition}}\left|\left\{\mathbf U\Biggm|\begin{aligned}
\text{$\mathbf U=(U_b)_{b\in B}$ \thinspace$B$-tableau }&\text{with
entries in $\mathbb Z_{>0}$}\\\text{such that $\sh(\mathbf U)
=\boldsymbol\lambda$ and}&\text{ $\ \forall b,\ \wt(U_b)=\mu^{(b)}$}
\end{aligned}\right\}\right|\;\mathbf s_{\boldsymbol\lambda}.
\end{multline*}

On the other hand, Equation~\ref{eq:PfPrCMRB} and
Lemma~\ref{le:DescCompoStdTab}~\ref{it:LeDCSTa} imply that
$$y_{c_1,b_1}*y_{c_2,b_2}*\cdots*y_{c_k,b_k}=\sum_{\substack{\alpha\in
B\wr\mathfrak S_n\\[2pt]D(\alpha)\preccurlyeq\mathbf c}}\alpha=
\sum_{\substack{\mathbf T\text{ standard $B$-tableau}\\[2pt]
D(\mathbf T)\preccurlyeq\mathbf c}}t_{\mathbf T},$$
so that
$$\Theta_B(y_{c_1,b_1}*y_{c_2,b_2}*\cdots*y_{c_k,b_k})=\sum_{\boldsymbol
\lambda\text{ $B$-partition}}\left|\left\{\mathbf T\Biggm|\begin{aligned}
\text{$\mathbf T$ standard }&\text{$B$-tableau with}\\\text{$\sh(\mathbf
T)=\boldsymbol\lambda$ }&\text{and $D(\mathbf T)\preccurlyeq\mathbf c$}
\end{aligned}\right\}\right|\;\mathbf s_{\boldsymbol\lambda}.$$
The desired result follows now from Lemma~\ref{le:DescCompoStdTab}
\ref{it:LeDCSTb}.
\end{proof}

We can now state and prove the main properties of $\Theta_B$.
\begin{theorem}
\label{th:ThetaCoplac}
The map $\Theta_B:\mathscr Q(B)\to\Lambda(B)$ is a surjective morphism
of graded bialgebras, with kernel $\mathscr Q(B)\cap\mathscr Q(B)^\circ$.
It is compatible with the pairings $\varpi_{\mathrm{tot}}$ on $\mathscr
Q(B)$ and $\langle?,?\rangle$ on $\Lambda(B)$, in the sense that
$$\varpi_{\mathrm{tot}}=\langle\Theta_B(?),\Theta_B(?)\rangle.$$
The restriction of $\Theta_B$ to $\mathscr D(B)$ is the unique algebra
homomorphism that maps $y_{n,b}$ to $h_n(b)$, where $n$ is a positive
integer and $b\in B$; this restriction also is surjective, with kernel
$\mathscr D(B)\cap\mathscr D(B)^\circ$.
\end{theorem}
\begin{proof}
Lemma~\ref{le:ThetaCoplac} implies that the restriction of $\Theta_B$ to
$\mathscr D(B)$ is a morphism of graded algebras, for the associative
algebra $\mathscr D(B)$ is generated by the elements $y_{n,b}$, and that
this restriction is surjective, for the algebra $\Lambda(B)$ is generated
by the elements $h_n(b)$. The coproducts of $\mathscr D(B)$ and
$\Lambda(B)$ being characterized by the equations
$$\Delta(y_{n,b})=\sum_{n'=0}^ny_{n',b}\otimes y_{n-n',b}\quad\text{and}
\quad\Delta(h_n(b))=\sum_{n'=0}^nh_{n'}(b)\otimes h_{n-n'}(b)$$
(with the convention that $y_{0,b}$ and $h_0(b)$ are the unit of
the algebras $\mathscr D(B)$ and $\Lambda(B)$, respectively), we also
see that $\Theta_B\bigr|_{\mathscr D(B)}$ preserves the coproducts. To
sum up, $\Theta_B$ is a surjective morphism of graded bialgebras.

We have seen in the proof of Proposition~\ref{pr:CoplacIsSubbialg}
\ref{it:PrCISc} that for each pair $(\mathbf T,\mathbf T')$ of standard
$B$-tableaux with the same number of boxes, there holds
$$\varpi_{\mathrm{tot}}(t_{\mathbf T},t_{\mathbf T'})
=\begin{cases}1&\text{if $\sh(\mathbf T)=\sh(\mathbf T')$,}\\0&
\text{otherwise,}\end{cases}$$
which implies that $\varpi_{\mathrm{tot}}(t_{\mathbf T},t_{\mathbf T'})
=\langle\mathbf s_{\sh(\mathbf T)},\mathbf s_{\sh(\mathbf T')}\rangle
=\langle\Theta_B(t_{\mathbf T}),\Theta_B(t_{\mathbf T'})\rangle$.
Therefore $\Theta_B$ is compatible with the pairings $\varpi_{\mathrm{tot}}$
and $\langle?,?\rangle$. In turn, this assertion, the fact that
$\langle?,?\rangle$ is a perfect pairing on $\Lambda(B)$ and the
surjectivity of $\Theta_B$ imply that the kernel of $\Theta_B$ is
equal to $\mathscr Q(B)\cap\mathscr Q(B)^\circ$. The surjectivity of
the restriction $\Theta_B\bigr|_{\mathscr D(B)}$ implies likewise that
the kernel of $\Theta_B\bigr|_{\mathscr D(B)}$ is equal to
$\mathscr D(B)\cap\mathscr D(B)^\circ$.

We then arrive at the following commutative diagram of graded bialgebras
$$\xymatrix@!0@C=148pt@R=42pt{&\mathscr D(B)\ar@{->>}
[dl]_{\Theta_B\bigr|_{\mathscr D(B)}}\ar@{^{(}->}[r]\ar@{->>}[d]&
\mathscr Q(B)\ar@{->>}[d]\\\Lambda(B)&\mathscr D(B)/(\mathscr D(B)
\cap\mathscr D(B)^\circ)\ar[r]^\simeq\ar[l]_(.6)\simeq&\mathscr Q(B)
/(\mathscr Q(B)\cap\mathscr Q(B)^\circ).}$$
An easy chase in this diagram shows that there exists a unique
homomorphism of $\mathbb K$-modules from $\mathscr Q(B)$ to
$\Lambda(B)$ which factorizes through $\mathscr Q(B)/(\mathscr Q(B)\cap
\mathscr Q(B)^\circ)$ and which extends $\Theta_B\bigr|_{\mathscr D(B)}$,
and that this homomorphism is a morphism of graded bialgebras. This
isomorphism is of course $\Theta_B$, which concludes the proof of the
theorem.
\end{proof}

\begin{corollary}
There holds
$$\mathscr D(B)\subseteq\mathscr Q(B)\subseteq\mathscr D(B)+\mathscr
Q(B)^\circ\quad\text{and}\quad\mathscr D(B)\cap\mathscr Q(B)^\circ=
\mathscr D(B)\cap\mathscr D(B)^\circ.$$
\end{corollary}
\begin{proof}
The inclusion $\mathscr D(B)\subseteq\mathscr Q(B)$ gives rise to an
injective map
$$\mathscr D(B)/\ker\bigl(\Theta_B\bigr|_{\mathscr D(B)}\bigr)
\hookrightarrow\mathscr Q(B)/\ker\Theta_B.$$
This latter is surjective, for the restriction $\Theta_B
\bigr|_{\mathscr D(B)}$ has the same image as $\Theta_B$. Using
Theorem~\ref{th:ThetaCoplac}, we arrive at the isomorphism $\mathscr
D(B)/(\mathscr D(B)\cap\mathscr D(B)^\circ)\stackrel\simeq
\longrightarrow\mathscr Q(B)/(\mathscr Q(B)\cap\mathscr Q(B)^\circ)$.
The corollary follows from this by standard arguments.
\end{proof}

We now have a big commutative diagram of graded bialgebras
\begin{equation}
\raisebox{43pt}{\xymatrix@!0@C=60pt@R=40pt{&&\mathscr F(\mathbb KB)
\ar[rr]^{\simeq}_{\varpi_{\mathrm{tot}}{}^\flat}\ar@{->>}[dr]&&\mathscr
F(\mathbb KB)^\vee\ar@{->>}[dr]&&\\&\mathscr Q(B)\ar@{^{(}->}[ur]
\ar@{->>}[dr]_{\Theta_B}&&\mathscr P(B)\ar[rr]^{\simeq}&&\mathscr
Q(B)^\vee\ar@{->>}[dr]&\\\mathscr D(B)\ar@{->>}[rr]\ar@{^{(}->}[ur]
&&\;\Lambda(B)\;\ar[rr]^{\simeq}_{\langle?,?\rangle^\flat}\ar@{^{(}->}[ur]
&&\;\Lambda(B)^\vee\;\ar@{^{(}->}[rr]\ar@{^{(}->}[ur]_{\Theta_B{}^\vee}
&&\mathscr D(B)^\vee.}}
\label{eq:DiagThetaCoplac}
\end{equation}

Given a standard $B$-tableau $\mathbf T$, let us denote by
$u_{\mathbf T}$ the class modulo $\mathscr Q(B)^\circ$ of an
$\alpha\in B\wr\mathfrak S_n$ such that $\mathbf P(\alpha)=\mathbf T$
(this class does not depends on the choice of $\alpha$). Using the
pairings, one checks rather easily that for any $B$-partition
$\boldsymbol\lambda$, the map from $\Lambda(B)$ to $\mathscr P(B)$ in
the diagram (\ref{eq:DiagThetaCoplac}) sends an element
$\mathbf s_{\boldsymbol\lambda}$ to
$$\sum_{\substack{\mathbf T\text{ standard $B$-tableau}\\[2pt]
\sh(\mathbf T)=\boldsymbol\lambda}}u_{\mathbf T}.$$

Finally, one may observe that the sequences~(\ref{eq:DualSubquot})
of homomorphisms, applied to the case $M=\mathscr F(\mathbb KB)$,
$S=\mathscr Q(B)$ and $T=\varpi_{\mathrm{tot}}{}^\flat(\mathscr Q(B))$,
show the existence of a symmetric pairing on $\mathscr Q(B)/(\mathscr
Q(B)\cap\mathscr Q(B)^\circ)$, which is perfect thanks to
Proposition~\ref{pr:CoplacIsSubbialg}~\ref{it:PrCISd}. This pairing is
of course equal to $\langle?,?\rangle$ under the isomorphism
$\mathscr Q(B)/(\mathscr Q(B)\cap\mathscr Q(B)^\circ)\cong\Lambda(B)$
defined by $\Theta_B$.

\subsection{Consequences for the Solomon descent theory}
\label{ss:CsqSoloTheo}
We now use the construction presented in the previous section to
complement the results of Section~\ref{ss:SoloHomo}. We consider a
finite abelian group $G$, we call $\Gamma=\Irr(G)$ its dual, and
we view $\Gamma$ as an object of $\mathscr E$ as explained in
Section~\ref{ss:CatFmwk}. The Frobenius characteristic $\ch$ is an
isomorphism of graded bialgebras from $\Rep(G)$ onto $\Lambda(\Gamma)$,
and there holds $\Theta_\Gamma\bigr|_{\mathscr D(\Gamma)}=\ch\circ\,\theta_G$,
because both members are homomorphisms of algebras which map $y_{n,\gamma}$
to $h_n(\gamma)$, where $n$ is a positive integer and $\gamma\in\Gamma$.
Therefore the diagrams~(\ref{eq:DiagrThetaVee}) and
(\ref{eq:DiagThetaCoplac}) agree and can be fused together.

We have recalled in Section~\ref{ss:CharacWreathProd} the construction
of the irreducible characters $\boldsymbol\chi^{\boldsymbol\lambda}$
of the wreath product $G\wr\mathfrak S_n$, indexed by the
$\Gamma$-partitions $\boldsymbol\lambda$ of $n$. From the diagram
$$\xymatrix{\MR(\mathbb K\Gamma)\ar@{->>}[d]_{\theta_G}\ar@{=}[r]&
\mathscr D(\Gamma)\,\ar@{^{(}->}[r]&\mathscr Q(\Gamma)\ar@{->>}[d]_{
\Theta_\Gamma}\\\Rep(G)\ar[rr]^\simeq_{\ch}&&\Lambda(\Gamma),}$$
we see that the homomorphism of $\mathbb K$-modules $\tilde\theta_G$
from $\mathscr Q(\Gamma)$ to $\Rep(G)$, defined by $\tilde\theta_G
(t_{\mathbf T})=\boldsymbol\chi^{\sh(\mathbf T)}$ for any standard
$\Gamma$-tableau $\mathbf T$, is a graded morphism of bialgebras
which extends $\theta_G$ and which is compatible with the pairings
$\varpi_{\mathrm{tot}}$ on $\mathscr Q(\Gamma)$ and $\beta_{\mathrm{tot}}$
on $\Rep(G)$.

Let $n$ be a positive integer \vspace{2pt}and $\mathbf c=((c_1,\gamma_1),
(c_2,\gamma_2),\ldots,(c_k,\gamma_k))$ be a $\Gamma$-composition of
size $n$. Set $\mathbf c^+=(c_1,c_2,\ldots,c_k)$ and
$\boldsymbol{\tilde\gamma}=\gamma_1^{\otimes c_1}\otimes\gamma_2^{\otimes
c_2}\otimes\cdots\otimes\gamma_k^{\otimes c_k}$. The character
$$\theta_G\bigl(y_{c_1,\gamma_1}*y_{c_2,\gamma_2}*\cdots*y_{c_k,\gamma_k}
\bigr)=\Ind_{G\wr\mathfrak S_{\mathbf c^+}}^{G\wr\mathfrak S_n}
\bigl(\eta_{c_1}(\gamma_1)\otimes\eta_{c_2}(\gamma_2)\otimes\cdots
\otimes\eta_{c_k}(\gamma_k)\bigr)$$
of $G\wr\mathfrak S_n$ is induced from a linear character of
$G\wr\mathfrak S_{\mathbf c^+}$. It can therefore be realized by a
representation on the $\mathbb C$-vector space with basis
$(G\wr\mathfrak S_n)/(G\wr\mathfrak S_{\mathbf c^+})$. As we saw during
the proof of Theorem~\ref{th:SymmSoloHomo}, this set is in natural
bijection with $X_{\mathbf c^+}\cong\{\rho\cdot(\boldsymbol{\tilde
\gamma}\#e_n)\mid\rho\in X_{\mathbf c^+}\}$. After translation in the
notation of Section~\ref{ss:ColDescComp}, this result means that
$$\theta_G\left(\sum_{\substack{\alpha\in\Gamma\wr\mathfrak S_n\\[2pt]
D(\alpha)\preccurlyeq\mathbf c}}\alpha\right)$$
is the character of a representation of $G\wr\mathfrak S_n$ on the
$\mathbb C$-vector space with basis
$$\{\alpha\in\Gamma\wr\mathfrak S_n\mid D(\alpha)\preccurlyeq\mathbf c\}.$$
Using a quotient construction, we may substitute equalities to the
inequalities $D(\alpha)\preccurlyeq\mathbf c$ in both formulas above. A
representation of $G\wr\mathfrak S_n$ whose character is 
$$\theta_G\left(\sum_{\substack{\alpha\in\Gamma\wr\mathfrak S_n\\[2pt]
D(\alpha)=\mathbf c}}\alpha\right)$$
is called a descent representation. Descent representations are
in particular studied in \cite{Adin-Brenti-Roichman04} by Adin, Brenti
and Roichman in the case $G=\{\pm1\}$ and in the forthcoming
paper~\cite{Bagno-Biagioli??} by Bagno and Biagioli in the case
$G=\mathbb Z/r\mathbb Z$. Our methods provide alternative proofs for some of
their results; for instance, Formula (\ref{eq:DescClsAsCoplacCls}) and the
equality $\tilde\theta_G{}(t_{\mathbf T})=\boldsymbol\chi^{\sh(\mathbf T)}$
imply the decomposition into irreducible characters
$$\theta_G\left(\sum_{\substack{\alpha\in\Gamma\wr\mathfrak S_n\\[2pt]
D(\alpha)=\mathbf c}}\alpha\right)=
\sum_{\substack{\text{$\boldsymbol\lambda$ $\Gamma$-partition}\\[2pt]
\|\boldsymbol\lambda\|=\|\mathbf c\|}}\left|\left\{\mathbf T\Biggm|
\begin{aligned}\text{$\mathbf T$ standard }&\text{$\Gamma$-tableau with}\\
\text{$\sh(\mathbf T)=\boldsymbol\lambda$ }&\text{and $D(\mathbf T)=
\mathbf c$}\end{aligned}\right\}\right|\;\boldsymbol{\chi^{\lambda}},$$
which generalizes Theorems~4.1 and~5.9 in \cite{Adin-Brenti-Roichman04}.


\section{Coloured quasisymmetric functions}
\label{se:ColQSymFun}
Motivated by problems of enumeration of permutations having a given
descent type, Gessel discovered in 1984 a link between Solomon's
descent algebra for the symmetric group and symmetric functions. More
precisely, he introduces in \cite{Gessel84} an algebra $\QSym$ of
`quasisymmetric functions,' which are polynomials in a countable and
totally ordered set of variables enjoying a certain symmetry property.
The algebra $\QSym$ is graded by the degree of polynomials (that is, the
homogeneous components of a quasisymmetric function are quasisymmetric),
which we write $\QSym=\bigoplus_{n\geq0}\QSym_n$. Gessel endows each
graded component $\QSym_n$ with the structure of a coalgebra and observes
that the dual algebra $\QSym_n^*$ can be identified with Solomon's
descent algebra $\Sigma_{\mathfrak S_n}$ for the symmetric group. Gessel
observes further that $\QSym_n$ contains the set $\Lambda_n$ of
homogeneous symmetric polynomials of degree $n$. Now $\Lambda_n$ is
isomorphic to its dual thanks to the usual inner product on symmetric
polynomials, and it is also isomorphic to the character ring $R(\mathfrak
S_n)$ of the symmetric group $\mathfrak S_n$ thanks to the characteristic
map. The inclusion $\Lambda_n\hookrightarrow\QSym_n$ gives then by
duality a surjection $\Sigma_{\mathfrak S_n}\cong\QSym_n^\vee
\twoheadrightarrow\Lambda_n^\vee\cong R(\mathfrak S_n)$, which Gessel
identifies with the Solomon map $\theta_{\mathfrak S_n}$.
$$\xymatrix{&\mathbb Z\mathfrak S_n&&&&\\\Sigma_{\mathfrak S_n}
\ar@{^{(}->}[ur]\ar@{->>}[dr]_{\theta_{\mathfrak S_n}}\ar[rr]^\simeq
&&\QSym_n^\vee\ar@{->>}[dr]&&&\QSym_n.\\&R(\mathfrak S_n)\ar[rr]^\simeq
&&\Lambda_n^\vee\ar@{=}[r]^\sim&\Lambda_n\ar@{^{(}->}[ur]&}$$

This picture was completed in 1995 by two independent groups of people.
On the one hand, Malvenuto and Reutenauer \cite{Malvenuto-Reutenauer95}
endow the space $\mathscr F=\bigoplus_{n\geq0}\mathbb Z\mathfrak{S_n}$
with the structure of a graded bialgebra by defining the external
product and the coproduct. They show that $\Sigma=\bigoplus_{n\geq0}
\Sigma_{\mathfrak S_n}$ is a graded subbialgebra of $\mathscr F$.
They endow Gessel's algebra $\QSym$ with a second coproduct, different
from Gessel's one, and they observe that this operation turns the
algebra $\QSym$ into a graded bialgebra, which they identify with the
graded dual of $\Sigma$.

On the other hand, Gelfand, Krob, Lascoux, Leclerc, Retakh and Thibon
\cite{GKLLRT95} introduce a graded module $\Sym=\bigoplus_{n\geq0}\Sym_n$
of `non-commutative symmetric functions.' They endow each graded
component $\Sym_n$ with an associative product with unit, which they
call the internal product, and find an explicit isomorphism between the
resulting algebra $\Sym_n$ and Solomon's descent algebra
$\Sigma_{\mathfrak S_n}$. Defining an external product and a coproduct,
they also endow $\Sym$ with the structure of a graded bialgebra, in such
a way that $\Sym$ can be identified as a graded bialgebra to $\Sigma$
and to the graded dual of $\QSym$. The pairing between $\Sym$ and $\QSym$
is made explicit through the use of bases; it reminds of the inner product
on the bialgebra $\Lambda$ of symmetric functions. Finally $\Lambda$
can be recovered as the quotient of $\Sym$ obtained by making commutative
the variables.

We want to generalize these works to the multidimensional case. To
this aim, we fix a finite set $B$ endowed with a linear order. As
in Section~\ref{ss:CatFmwk}, we denote the free $\mathbb K$-module
with basis $B$ by $\mathbb KB$ and define a Mantaci-Reutenauer
subbialgebra $\mathscr D(B)$ in $\mathscr F(\mathbb KB)$. In
Section~\ref{ss:WordRealFG}, we present a realization of the algebra
$\mathscr F(\mathbb KB)$ in terms of `coloured' free quasisymmetric
functions. The dependence of our realization on the linear order on
$B$ may seem cumbersome, but is a necessary step so that the quotient
map $\mathscr F(\mathbb KB)\to \mathscr F(\mathbb KB)/\mathscr
D(B)^\circ$ corresponds to make the variables commutative. In
Section~\ref{ss:PoirQuasisymFonc}, we show that our construction
yields some of Poirier's quasisymmetric functions.

We fix for the whole Section~\ref{se:ColQSymFun} an infinite alphabet
$\mathscr A$ and we endow the product $\mathscr A\times B$ with the
lexicographical order.

\subsection{The word realization of $\mathscr F(k\Gamma)$}
\label{ss:WordRealFG}
Let $w=x_1x_2\cdots x_n$ be a word with letters $x_i=(a_i,b_i)$ in
$\mathscr A\times B$. Denoting by $\sigma\in\mathfrak S_n$ the
standardization of $w$, we may form the element
$$\std_B(w)=\sigma\cdot(b_1,b_2,\ldots,b_n;e_n)=(b_{\sigma^{-1}(1)},
b_{\sigma^{-1}(2)},\ldots,b_{\sigma^{-1}(n)};\sigma)$$
of $B\wr\mathfrak S_n$; we call it the $B$-standardization of $w$.
(This element $\std_B(w)$ is called `standard signed permutation' of
$w$ by Poirier; see \cite{Poirier98}, p.~322.) As an example, let
$B=\{\,\Bar\ ,\Bar{\Bar\ }\,\}$ with $\Bar\ <\Bar{\Bar\ }$, and let
$\mathscr A=\{u,v,w,\ldots\}$ with the usual alphabetical order. We
denote the letters $(u,\Bar\ \,)$, $(u,\Bar{\Bar\ }\,)$, etc.\ by
$\Bar u$, $\Bar{\Bar u}$, etc. Then the standardization of the word
$w=\Bar u\Bar v\Bar u\Bar{\Bar v}\Bar w\Bar{\Bar u}\Bar v$ is
$\std_B(w)=\bigl[(1426735)\cdot(\,\Bar\ ,\Bar\ ,\Bar\ ,\Bar{\Bar\
},\Bar\ ,\Bar{\Bar\ },\Bar\ ;e_7)\bigr]$.

We denote the algebra of non-commutative formal power series on
the set $\mathscr A\times B$ with coefficients in $\mathbb K$ by
$\mathbb K\llangle\mathscr A\times B\rrangle$; thus elements of
$\mathbb K\llangle\mathscr A\times B\rrangle$ are (possibly infinite)
linear combinations of words on the alphabet $\mathscr A\times B$.
We denote the algebra of commutative formal power series on the
set $\mathscr A\times B$ with coefficients in $\mathbb K$ by
$\mathbb K[[\mathscr A\times B]]$; elements of this algebra may be
viewed as (possibly infinite) linear combinations of $\mathscr A\times
B$-weights. There is an obvious morphism of $\mathbb K$-algebras from
$\mathbb K\llangle\mathscr A\times B\rrangle$ onto
$\mathbb K[[\mathscr A\times B]]$, which maps each word $w$ on the
alphabet $\mathscr A\times B$ to its weight.

We denote by $\Phi:\mathscr F(\mathbb KB)\to\mathbb K\llangle\mathscr
A\times B\rrangle$ the map which sends an element $\alpha\in
(B\wr\mathfrak S_n)$ to the sum of all words $w$ such that $\alpha$ is
the $B$-standardization of $w$:
$$\Phi(\alpha)=\sum_{\substack{w\in\langle\mathscr A\times B\rangle\\[2pt]
\std_B(w)=\alpha}}w.$$

\begin{theorem}
\label{th:WordReal}
\begin{enumerate}
\item\label{it:ThWRa}The map $\Phi$ is an injective morphism of algebras
from $\mathscr F(\mathbb KB)$ to $\mathbb K\llangle\mathscr A\times
B\rrangle$.
\item\label{it:ThWRb}Let $I$ be the kernel of the canonical morphism
from $\mathbb K\llangle\mathscr A\times B\rrangle$ onto
$\mathbb K[[\mathscr A\times B]]$. Then $\Phi^{-1}(I)=\mathscr D(B)^\circ$.
\end{enumerate}
\end{theorem}
\begin{proof}
\begin{enumerate}
\item Let $n$ and $n'$ be two positive integers and let $w$ and $w'$ be
two words on the alphabet $\mathscr A\times B$ of length $n$ and $n'$,
respectively. If we denote by $\sigma\in\mathfrak S_n$, $\sigma'\in
\mathfrak S_{n'}$ and $\pi\in\mathfrak S_{n+n'}$ the standardizations
of the words $w$, $w'$ and $ww'$, respectively, then $\sigma$ is the
standardization of the word $\pi(1)\pi(2)\cdots\pi(n)$ and $\sigma'$
is the standardization of the word $\pi(n+1)\pi(n+2)\cdots\pi(n+n')$;
in other words, there exists $\rho\in X_{(n,n')}$ such that $\pi=\rho
(\sigma\times\sigma')$.

Now let $\alpha\in B\wr\mathfrak S_n$ and $\alpha'\in B\wr\mathfrak
S_{n'}$. We write $\alpha=\sigma\cdot(b_1,b_2,\ldots,b_n;e_n)$,
$\alpha'=\sigma'\cdot(b'_1,b'_2,\ldots,b'_{n'};e_{n'})$, $w=x_1x_2\cdots
x_n$ and $w'=x'_1x'_2\cdots x_{n'}$. Given a letter $x=(a,b)$ in
$\mathscr A\times B$, we say that $b$ is the colour of $x$. Then
\begin{align*}
\alpha=\std_B(w)\;\text{ and }&\;\alpha'=\std_B(w')\\[3pt]
&\Longleftrightarrow\left\{\begin{aligned}&\text{$\sigma$ is the
standardization of $w$, $\sigma'$ is the standardization of $w'$,}\\
&\text{$b_i$ is the colour of $x_i$ and $b'_j$ is the colour of
$x'_j$,}\end{aligned}\right.\\[3pt]
&\Longleftrightarrow\left\{\begin{aligned}&\text{$\exists\rho\in
X_{(n,n')}$ such that $\rho(\sigma\times\sigma')$ is the standardization
of $ww'$,}\\&\text{$b_i$ is the colour of $x_i$ and $b'_j$ is the colour
of $x'_j$,}\end{aligned}\right.\\[3pt]
&\Longleftrightarrow\left\{\begin{aligned}&\text{$\exists\rho\in
X_{(n,n')}$ such that}\\&\std_B(ww')=\rho(\sigma\times\sigma')\cdot
(b_1,b_2,\ldots,b_n,b'_1,b'_2,\ldots,b'_{n'};e_{n+n'}).\end{aligned}\right.
\end{align*}
This proves that $\Phi$ is a morphism of algebras. The injectivity of
$\Phi$ is an obvious consequence of the fact that $\mathscr A$ was
chosen infinite.
\item Let $n$ be a positive integer and let $\alpha=(b_1,b_2,\ldots,
b_n;\sigma)$ and $\alpha'$ be two elements in $B\wr\mathfrak S_n$ such
that $\alpha\underset A\sim\alpha'$. Then there exists a simple
transposition $s_i\in\mathfrak S_n$ such that $\alpha'=\alpha\cdot s_i$,
the index $i\in\{1,2,\ldots,n-1\}$ enjoying moreover the property that
the map $j\mapsto b_j$ is not constant on the interval $[\sigma(i),
\sigma(i+1)]$ or that the inequality $|\sigma(i+1)-\sigma(i)|>1$ holds.

Then in each word $w=x_1x_2\cdots x_n$ of length $n$ on the alphabet
$\mathscr A\times B$ whose $B$-standardization is $\alpha$, the letters
$x_i$ and $x_{i+1}$ differ. The word $w'=x_1x_2\cdots x_{i-1}x_{i+1}
x_ix_{i+2}\cdots x_n$ obtained from $w$ by exchanging the letters $x_i$
and $x_{i+1}$ has thus $\alpha\cdot s_i=\alpha'$ for $B$-standardization,
and the map $w\mapsto w'$ is a bijective correspondence
$$\left\{w\Biggm|\begin{aligned}w\ &\text{word on }\mathscr A\times B
\text{ such}\\&\text{that }\std_B(w)=\alpha\end{aligned}\right\}
\stackrel\simeq\longrightarrow
\left\{w'\Biggm|\begin{aligned}w'&\text{ word on }\mathscr A\times B
\text{ such}\\&\text{that }\std_B(w')=\alpha'\end{aligned}\right\}.$$
Therefore $\Phi(\alpha)$ and $\Phi(\alpha')$ have the same image in
$\mathbb K[[\mathscr A\times B]]$, for $w$ and $w'$ have the same
weight. By Proposition~\ref{pr:CombMRBialg}~\ref{it:PrCMRBc}, this
implies that $\Phi\bigl(\mathscr D(B)^\circ\bigr)\subseteq I$.

The morphism $\Phi$ defines therefore a map $\overline\Phi$ from
$\mathscr F(\mathbb KB)/\mathscr D(B)^\circ$ to $\mathbb K[[\mathscr
A\times B]]$. Assertion~\ref{it:ThWRb} will then be proved as soon as
the injectivity of $\overline\Phi$ is established.

We associate a $B$-composition $C(\mu)$ to each $(\mathscr A\times
B)$-weight $\mu$ as follows: we list in increasing order $(a_1,b_1)
<(a_2,b_2)<\cdots<(a_k,b_k)$ the elements $(a,b)$ in the support of
the multiset $\mu$, and we then define $C(\mu)$ as the sequence
$((\mu(a_1,b_1),b_1),(\mu(a_2,b_2),b_2),\ldots,(\mu(a_k,b_k),b_k))$.
One checks that $R(\std_B(w))\preccurlyeq C(\wt(w))$ for any word $w$
on the alphabet $\mathscr A\times B$.

Let $z$ be a non-zero element in $\mathscr F(\mathbb KB)/\mathscr
D(B)^\circ$. By Proposition~\ref{pr:CombMRBialg}, $z$ has an antecedent
in $\mathscr F(\mathbb KB)$ of the form $\sum_{j\in J}a_j\alpha_j$,
where $J$ is a finite index set, $a_j\in\mathbb K\setminus\{0\}$, and
the elements $\alpha_j\in B\wr\mathfrak S_n$ are such that all
$B$-compositions $R(\alpha_j)$ are different. We may then find $j_0\in J$
such that $R(\alpha_{j_0})$ is a minimal element of the set
$\{R(\alpha_j)\mid j\in J\}$ with respect to the refinement order
$\preccurlyeq$, and we may find a word $w$ on the alphabet $\mathscr
A\times B$ such that $\std_B(w)=\alpha_{j_0}$ and $C(\wt(w))=
R(\alpha_{j_0})$. Then $\wt(w)$ appears in $\overline\Phi(z)$ with the
coefficient $a_{j_0}\neq0$, which entails that $\overline\Phi(z)\neq0$.

Therefore $\overline\Phi$ is injective, which completes the proof.
\end{enumerate}
\end{proof}

Assertion~\ref{it:ThWRa} of Theorem~\ref{th:WordReal} says that we can
find a realization of the algebra $\mathscr F(\mathbb KB)$ in terms of
free (non-commutative) quasisymmetric functions. Assertion~\ref{it:ThWRb}
says that the quotient map from $\mathscr F(\mathbb KB)$ onto
$\mathscr F(\mathbb KB)/\mathscr D(B)^\circ\cong\mathscr D(B)^\vee$ is
obtained in this realization by making commutative all words
$w\in\langle\mathscr A\times B\rangle$. This can be translated into
the commutative diagram
\begin{equation}
\raisebox{23pt}{\xymatrix{\mathscr F(\mathbb KB)\,\ar@{^{(}->}[r]
^\Phi\ar@{->>}[d]&\mathbb K\llangle\mathscr A\times B\rrangle
\ar@{->>}[d]\\\mathscr F(\mathbb KB)/\mathscr D(B)^\circ
\ar@{^{(}->}[r]^{\overline\Phi}&\mathbb K[[\mathscr A\times B]].}}
\label{eq:DiagWordReal}
\end{equation}
One can find a similar description of all the algebras that appear in
the diagram~(\ref{eq:DiagThetaCoplac}); for instance, the quotient map
from $\mathscr F(\mathbb KB)$ onto $\mathscr F(\mathbb KB)/\mathscr
Q(B)^\circ=\mathscr P(B)$ amounts to look at the words $w\in\langle
\mathscr A\times B\rangle$ modulo the Knuth relation $\underset K\sim$
of Section~\ref{ss:RSOCorr}.

\subsection{Poirier's quasisymmetric functions}
\label{ss:PoirQuasisymFonc}
Let $\QSym(B)$ denote the image in $\mathbb K[[\mathscr A\times B]]$ of
the map $\overline\Phi$ in the diagram~(\ref{eq:DiagWordReal}). In this
section, we describe $\QSym(B)$ explicitly and compare it with Poirier's
algebra of quasisymmetric functions.

By Proposition~\ref{pr:CombMRBialg}, the class modulo $\mathscr
D(B)^\circ$ of an element $\alpha\in B\wr\mathfrak S_n$ is determined
by its receding composition $R(\alpha)$. A stronger assertion holds:
it is possible to find a combinatorial description of $\overline\Phi
\bigl(\alpha+\mathscr D(B)^\circ\bigr)$ based on the sole data of
$R(\alpha)$.

Indeed let $\mathbf c=((c_1,b_1),(c_2,b_2),\ldots,(c_k,b_k))$
be a $B$-composition of size say $n$, set $t_i=c_1+c_2+\cdots+c_i$
for each $i$, and set
$$(\tilde b_1,\tilde b_2,\ldots,\tilde b_n)=(\underbrace{b_1,b_1,\ldots,
b_1}_{\text{$c_1$ times}},\underbrace{b_2,b_2,\ldots,b_2}_{\text{$c_2$
times}},\ldots,\underbrace{b_k,b_k,\ldots,b_k}_{\text{$c_k$ times}}).$$
From $\mathbf c$, we construct the set $S_{\mathbf c}$ of all $n$-uples
$(x_1,x_2,\ldots,x_n)\in(\mathscr A\times B)^n$ satisfying the three
following conditions: the sequence $(x_1,x_2,\ldots,x_n)$ is
non-decreasing; $x_{t_i}<x_{t_i+1}$ for each $i\in\{1,2,\ldots,k-1\}$;
the second component of $x_i\in\mathscr A\times B$ is $\tilde b_i$.
In other words, a $n$-uple $(x_1,x_2,\ldots,x_n)$ belongs to
$S_{\mathbf c}$ if and only if each $x_i$ can be written $(a_i,\tilde
b_i)$, where $(a_1,a_2,\ldots,a_n)$ is a non-decreasing sequence of
elements of $\mathscr A$ such that
$$\forall i\in\{1,2,\ldots,k-1\},\quad b_i\geq b_{i+1}\Longrightarrow
a_{t_i}<a_{t_i+1}.$$
By analogy with Formula~(2) on p.~324 in~\cite{Poirier98}, we define
the formal series in $\mathbb K[[\mathscr A\times B]]$
$$F_{\mathbf c}=\sum_{(x_1,x_2,\ldots,x_n)\in S_{\mathbf c}}x_1x_2\cdots
x_n.$$
For instance if $B$ is the set $\{\,\Bar\ ,\Bar{\Bar\ }\,\}$ with the order
$\Bar\ <\Bar{\Bar\ }$, then
\begin{xalignat*}3
F_{((2,\Bar\ \,))}&=\sum_{\substack{(x,y)\in\mathscr A^2\\[2pt]x\leq y}}
\Bar x\Bar y,&
F_{((2,\Bar{\Bar\ }\,))}&=\sum_{\substack{(x,y)\in\mathscr A^2\\[2pt]
x\leq y}}\Bar{\Bar x}\Bar{\Bar y},&
F_{((1,\Bar\ \,),(1,\Bar\ \,))}&=\sum_{\substack{(x,y)\in\mathscr A^2
\\[2pt]x<y}}\Bar x\Bar y,\\[4pt]
F_{((1,\Bar\ \,),(1,\Bar{\Bar\ }\,))}&=\sum_{\substack{(x,y)\in\mathscr
A^2\\[2pt]x\leq y}}\Bar x\Bar{\Bar y},&
F_{((1,\Bar{\Bar\ }\,),(1,\Bar\ \,))}&=\sum_{\substack{(x,y)\in\mathscr
A^2\\[2pt]x<y}}\Bar{\Bar x}\Bar y,&
F_{((1,\Bar{\Bar\ }\,),(1,\Bar{\Bar\ }\,))}&=\sum_{\substack{(x,y)\in
\mathscr A^2\\[2pt]x<y}}\Bar{\Bar x}\Bar{\Bar y}.
\end{xalignat*}

The following result is a rewriting of Lemma~11 in \cite{Poirier98};
it implies that the elements $F_{\mathbf c}$ form a basis of the
$\mathbb K$-module $\QSym(B)$, where $\mathbf c$ is a $B$-composition.
\begin{proposition}
\label{pr:ImgPhiBarIsQRib}
For each element $\alpha\in B\wr\mathfrak S_n$, there holds
$F_{R(\alpha)}=\overline\Phi\bigl(\alpha+\mathscr D(B)^\circ\bigr)$.
\end{proposition}
\begin{proof}
We take an element $\alpha\in B\wr\mathfrak S_n$, we write
$$\alpha=(\tilde b_1,\tilde b_2,\ldots,\tilde b_n;\sigma)\qquad
\text{and}\qquad R(\alpha)=((c_1,b_1),(c_2,b_2),\ldots,(c_k,b_k)),$$
and we set $t_i=c_1+c_2+\cdots+c_i$ for each $i$. The definition of
$R(\alpha)$ implies that
$$(\tilde b_1,\tilde b_2,\ldots,\tilde b_n)=(\underbrace{b_1,b_1,\ldots,
b_1}_{\text{$c_1$ times}},\underbrace{b_2,b_2,\ldots,b_2}_{\text{$c_2$
times}},\ldots,\underbrace{b_k,b_k,\ldots,b_k}_{\text{$c_k$ times}}),$$
that the permutation $\sigma^{-1}$ is increasing on each interval
$[t_{i-1}+1,t_i]$, and that
$$\forall i\in\{1,2,\ldots,k-1\},\quad b_i=b_{i+1}\Longrightarrow
\sigma(t_i)>\sigma(t_i+1).$$

Each sequence $(x_1,x_2,\ldots,x_n)\in S_{R(\alpha)}$ yields a word
$w=x_{\sigma(1)}x_{\sigma(2)}\cdots x_{\sigma(n)}$ with letters in
$\mathscr A\times B$. The definition of $S_{R(\alpha)}$ is so shaped
that the standardization of $w$ is $\sigma$; it follows that the
$B$-standardization of $w$ is $\sigma\cdot\bigl(\tilde b_{\sigma(1)},
\tilde b_{\sigma(2)},\ldots,\tilde b_{\sigma(n)};e_n\bigr)=\alpha$.
Conversely, each word $w$ with letters in $\mathscr A\times B$ whose
$B$-standardization is $\alpha$ can be written
$w=x_{\sigma(1)}x_{\sigma(2)}\cdots x_{\sigma(n)}$, where the
sequence $(x_1,x_2,\ldots,x_n)$ belongs to $S_{R(\alpha)}$.

We conclude that the image of
$$\Phi(\alpha)=\sum_{\substack{w\in\langle\mathscr A\times B\rangle\\[2pt]
\std_B(w)=\alpha}}w$$
under the canonical map from $\mathbb K\llangle\mathscr A\times B\rrangle$
to $\mathbb K[[\mathscr A\times B]]$ is equal to
$$\sum_{(x_1,x_2,\ldots,x_n)\in S_{R(\alpha)}}x_{\sigma(1)}x_{\sigma(2)}
\cdots x_{\sigma(n)}=\sum_{(x_1,x_2,\ldots,x_n)\in S_{R(\alpha)}}x_1x_2
\cdots x_n=F_{R(\alpha)}.$$
The proposition follows.
\end{proof}

Now let us enumerate the elements of $B$ in increasing order: $\bar b_1$,
$\bar b_2$, \dots, $\bar b_l$, where $l$ is the cardinality of $B$, and
let us review the definitions of a combinatorial nature that are needed
to introduce Poirier's theory of coloured quasisymmetric functions.
Since Poirier made a slight mistake (in \cite{Poirier98}, Lemma~8 does
not always agree with Formulas~(1) and (2) on p.~324), we will follow
Novelli and Thibon's presentation \cite{Novelli-Thibon04}.

A $l$-partite number is an element of $\mathbb N^l$; we view it
as a column matrix. Given a positive integer~$k$, a $l$-vector
composition of length $k$ is a $k$-uple of non-zero $l$-partite
numbers; it can be viewed as a sequence of column matrices, or
more simply as a matrix with non-negative integral entries in $l$
rows and $k$ columns which has at least one non-zero element in each
column.

Each $l$-vector composition $\mathbf I$ produces a formal power series
in $\mathbb K[[\mathscr A\times B]]$ called a monomial quasisymmetric
function of level $l$ and defined by
$$M_{\mathbf I}=\sum_{\substack{(a_1,a_2,\ldots,a_k)\in\mathscr
A^k\\[2pt]a_1<a_2<\cdots<a_k}}\left(\;\prod_{i=1}^l\ \prod_{j=1}^k
(a_j,\bar b_i)^{m_{ij}}\right),$$
where $(m_{ij})$ is the matrix that represents $\mathbf I$.
For instance in the case where $B$ is the set $\{\,\Bar\ ,\Bar{\Bar\ }\,\}$
with the order $\Bar\ <\Bar{\Bar\ }$, the monomial quasisymmetric
functions of level $l=2$ and of degree $2$ are
\begin{xalignat*}3
M_{\left(\begin{smallmatrix}2\\0\end{smallmatrix}\right)}&=
\sum_{x\in\mathscr A}\Bar x^2,&
M_{\left(\begin{smallmatrix}1\\1\end{smallmatrix}\right)}&=
\sum_{x\in\mathscr A}\Bar x\Bar{\Bar x},&
M_{\left(\begin{smallmatrix}0\\2\end{smallmatrix}\right)}&=
\sum_{x\in\mathscr A}\Bar{\Bar x}^2,\\[4pt]
M_{\left(\begin{smallmatrix}1&1\\0&0\end{smallmatrix}\right)}&=
\sum_{\substack{(x,y)\mathscr A^2\\[2pt]x<y}}\Bar x\Bar y,&
M_{\left(\begin{smallmatrix}1&0\\0&1\end{smallmatrix}\right)}&=
\sum_{\substack{(x,y)\mathscr A^2\\[2pt]x<y}}\Bar x\Bar{\Bar y},&
M_{\left(\begin{smallmatrix}0&1\\1&0\end{smallmatrix}\right)}&=
\sum_{\substack{(x,y)\mathscr A^2\\[2pt]x<y}}\Bar{\Bar x}\Bar y,\\[4pt]
M_{\left(\begin{smallmatrix}0&0\\1&1\end{smallmatrix}\right)}&=
\sum_{\substack{(x,y)\mathscr A^2\\[2pt]x<y}}\Bar{\Bar x}\Bar{\Bar y}.
\end{xalignat*}

Let $\mathbf I$ be a $l$-vector composition, represented by the matrix
$(m_{ij})$. We form the list of all pairs $(m_{ij},\bar b_i)$, reading
columnwise the entries of $(m_{ij})$ from top to bottom and from left
to right. Erasing in this list all the pairs whose first component
$m_{ij}$ is zero, we obtain a $B$-composition, which we call the
sequential reading of $\mathbf I$ and which we denote by
$\seqr(\mathbf I)$. For instance with our favorite set $B=\{\,\Bar\
,\Bar{\Bar\ }\,\}$ with the order $\Bar\ <\Bar{\Bar\ }$, the $l$-vector
compositions represented by the matrices
$$\begin{pmatrix}1&0&4\\3&2&1\end{pmatrix}\qquad\text{and}\qquad
\begin{pmatrix}1&0&0&4\\0&3&2&1\end{pmatrix}$$
have both sequential reading $((1,\Bar\ \,),(3,\Bar{\Bar\ }\,),
(2,\Bar{\Bar\ }\,),(4,\Bar\ \,),(1,\Bar{\Bar\ }\,))$. We see therefore
that the map $\mathbf I\mapsto\seqr(\mathbf I)$ is not injective.

This definition allows us to express each formal power series
$F_{\mathbf c}$ as a linear combination of monomial quasisymmetric
functions.
\begin{proposition}
\label{pr:QRibbonAsLCMonom}
For each $B$-composition $\mathbf c$, there holds
$$F_{\mathbf c}=\sum_{\substack{\mathbf I\text{ $l$-vector
composition}\\[2pt]\mathbf c\preccurlyeq\seqr(\mathbf
I)}}M_{\mathbf I}.$$
\end{proposition}
\begin{proof}
Let $S$ be the set of all non-decreasing finite sequences
$(x_1,x_2,\ldots,x_n)$ of elements of $\mathscr A\times B$. We define
a map $\psi$ from $S$ to the set of all $l$-vector compositions by
the following recipe. Let $(x_1,x_2,\ldots,x_n)$ in $S$; write $x_i=
(a_i,b_i)$ for each $i$; let $\tilde a_1$, $\tilde a_2$, \dots,
$\tilde a_k$ be the (distinct) elements of $\{a_i\mid1\leq i\leq n\}$
enumerated in increasing order. Then $\psi(x_1,x_2,\ldots,x_n)$ is
the $l$-vector composition of length $k$ represented by the matrix
$(m_{ij})$, where each $m_{ij}$ counts the number of times that the
element $(\tilde a_j,\bar b_i)$ appears in the sequence
$(x_1,x_2,\ldots,x_n)$.

Given a $l$-vector composition $\mathbf I$, we set $T_{\mathbf I}=
\psi^{-1}\bigl(\{\mathbf I\}\bigr)$. Then by definition
$$M_{\mathbf I}=\sum_{(x_1,x_2,\ldots,x_n)\in T_{\mathbf I}}
x_1x_2\cdots x_n.$$

Now let $\mathbf c$ be a $B$-composition. Routine arguments show that
$$\forall\mathbf x\in S,\quad\mathbf x\in S_{\mathbf c}
\Longleftrightarrow\mathbf c\preccurlyeq\seqr(\psi(\mathbf x)).$$
In other words, $S_{\mathbf c}$ is the disjoint union of the sets
$T_{\mathbf I}$, where $\mathbf I$ is an $l$-vector composition such that
$\mathbf c\preccurlyeq\seqr(\mathbf I)$. It follows that
\begin{align*}
F_{\mathbf c}
&=\sum_{(x_1,x_2,\ldots,x_n)\in S_{\mathbf c}}x_1x_2\cdots x_n\\[4pt]
&=\sum_{\substack{\mathbf I\text{ $l$-vector composition}\\[2pt]\mathbf
c\preccurlyeq\seqr(\mathbf I)}}\left(\sum_{(x_1,x_2,\ldots,x_n)\in
T_{\mathbf I}}x_1x_2\cdots x_n\right)\\[4pt]
&=\sum_{\substack{\mathbf I\text{ $l$-vector composition}\\[2pt]\mathbf
c\preccurlyeq\seqr(\mathbf I)}}M_{\mathbf I},
\end{align*}
which proves the proposition.
\end{proof}

Paraphrasing a construction of Poirier, Novelli and Thibon endow the
set of $l$-vector compositions with a partial order $\leq$ and define
for each $l$-vector composition $\mathbf I$ the formal power series
$$F_{\mathbf I}=\sum_{\substack{\mathbf J\text{ $l$-vector composition}
\\[2pt]\mathbf I\leq\mathbf J}}M_{\mathbf J},$$
which they call a quasi-ribbon function of level $l$. On the other side,
one can show quite easily that for each $B$-composition $\mathbf c$,
there exists a unique $l$-vector composition $\mathbf K(\mathbf c)$
such that
$$\left\{\mathbf J\Biggm|\begin{aligned}\text{$\mathbf J$ $l$-vector
}&\text{composition}\\\mathbf c\preccurlyeq&\;\seqr(\mathbf
J)\end{aligned}\right\}=\left\{\mathbf J\Biggm|\begin{aligned}
\text{$\mathbf J$ $l$-vector}&\text{ composition}\\\mathbf K&(\mathbf
c)\leq\mathbf J\end{aligned}\right\}.$$
With these notations, Proposition~\ref{pr:QRibbonAsLCMonom} asserts that
the formal power series $F_{\mathbf c}$ coincides with the quasi-ribbon
function $F_{\mathbf K(\mathbf c)}$.

Let us denote by $\NTQSym l$ the submodule of $\mathbb K[[\mathscr
A\times B]]$ spanned by the monomial quasisymmetric functions of level $l$.
Novelli and Thibon claim in \cite{Novelli-Thibon04} that $\NTQSym l$ is
a subalgebra of $\mathbb K[[\mathscr A\times B]]$, and moreover that
$\NTQSym l$ has the structure of a graded bialgebra, whose dual can be
identified to the Novelli-Thibon bialgebra $\NT(\mathbb KB)$. In this
context, Propositions~\ref{pr:ImgPhiBarIsQRib} and
\ref{pr:QRibbonAsLCMonom} imply that $\QSym(B)$ is a subalgebra of
$\NTQSym l$. It is amusing to note here that the graded algebra
$\QSym(B)$, which is isomorphic to the dual of the graded bialgebra
$\mathscr D(B)$, can also be viewed as a quotient of $\NTQSym l$, since
$\mathscr D(B)$ is a graded subbialgebra of $\NT(\mathbb KB)$.

We conclude this paper by mentioning that Aval, F.~Bergeron and N.~Bergeron
recently observed that coloured quasisymmetric functions of level $l=2$
appear in a completely different context. We refer the reader to their
paper \cite{Aval-Bergeron-Bergeron04} for additional details.

\bigskip
\noindent\parbox[t]{241pt}{\noindent
Pierre Baumann\\
Institut de Recherche Mathématique Avancée\\
Université Louis Pasteur et CNRS\\
7, rue René Descartes\\
67084 Strasbourg Cedex\\
France\\[6pt]
E-mail: \texttt{baumann@math.u-strasbg.fr}}
\parbox[t]{192pt}{\noindent
Christophe Hohlweg\\
The Fields Institute\\
222 College Street\\
Toronto Ontario M5T 3J1\\
Canada\\[6pt]
E-mail: \texttt{chohlweg@fields.utoronto.ca}}
\end{document}